\documentclass[12pt]{elsarticle}

\usepackage{amssymb}
\usepackage{amsmath}

\usepackage{rotating}
\usepackage{hyperref}

\usepackage[margin=1.75cm,nohead]{geometry} 

\journal{arXiv}

\begin{document}



\begin{frontmatter}


\title{On the zeros of Riemann's Xi Function}

\author{Akhila Raman }

\address{Email: akhila.raman@berkeley.edu. }



\begin{abstract}

We consider Riemann's Xi function $\xi(s)$ which is evaluated at $s = \frac{1}{2} + \sigma + i \omega$, given by $\xi(\frac{1}{2} + \sigma + i \omega)= E_{p\omega}(\omega)$, where $\sigma, \omega$ are real and compute its inverse Fourier transform given by $E_p(t)$. We study the properties of $E_p(t)$ and a promising new method is presented which could be used to show that the Fourier Transform of $E_p(t)$ given by $E_{p\omega}(\omega) = \xi(\frac{1}{2} + \sigma + i \omega)$ does not have zeros for finite and real $\omega$ when $0 < |\sigma| < \frac{1}{2}$, corresponding to the critical strip excluding the critical line.
\end{abstract}


\end{frontmatter}
   







\newcommand{\nocontentsline}[3]{}
\newcommand{\tocless}[2]{\bgroup\let\addcontentsline=\nocontentsline#1{#2}\egroup}
\tocless\section{\label{sec:level2} \textbf{Introduction} \protect\\  \lowercase{} }

It is well known that Riemann's Zeta function given by $\zeta(s) = \sum\limits_{m=1}^{\infty} \frac{1}{m^{s}}$ converges in the half-plane where the real part of $s$ is greater than 1. Riemann proved that  $\zeta(s)$  has an analytic continuation to the whole s-plane apart from a simple pole at $s = 1$ and that $\zeta(s)$  satisfies a symmetric functional equation given by $\xi(s) = \xi(1-s) =  \frac{1}{2} s (s-1) \pi^{-\frac{s}{2}} \Gamma(\frac{s}{2}) \zeta(s) $ where $\Gamma(s) = \int_{0}^{\infty} e^{-u} u^{s-1} du$ is the Gamma function.{\citep{FWE}} {\citep{JBC}} We can see that if Riemann's Xi function has a zero in the critical strip, then Riemann's Zeta function also has a zero at the same location. Riemann made his conjecture in his 1859 paper, that all of the non-trivial zeros of $\zeta(s)$ lie on the critical line with real part of $s = \frac{1}{2}$, which is called the Riemann Hypothesis.{\citep{BR}} \\

Hardy and Littlewood later proved that infinitely many of the zeros of $\zeta(s)$ are on the critical line with real part of $s=\frac{1}{2}$.{\citep{GH}} It is well known that $\zeta(s)$ does not have non-trivial zeros when real part of $s = \frac{1}{2} + \sigma + i \omega$, given by $\frac{1}{2} + \sigma \geq 1$ and $\frac{1}{2} + \sigma \leq 0$. In this paper, \textbf{critical strip} $0 < Re[s] < 1 $ corresponds to $0 \leq |\sigma| < \frac{1}{2}$.\\

In this paper, a \textbf{new method} is discussed and a specific solution is presented to prove Riemann's Hypothesis. If the specific solution presented in this paper is incorrect, it is \textbf{hoped} that the new method discussed in this paper will lead to a correct solution by other researchers.\\

In Section~\ref{sec:Section_2} to Section~\ref{sec:Section_A_1_6}, we prove Riemann's hypothesis by taking the analytic continuation of Riemann's Zeta Function derived from Riemann's Xi function  $\xi(\frac{1}{2} + \sigma + i \omega)= E_{p\omega}(\omega)$ and compute inverse Fourier transform of $E_{p\omega}(\omega)$ given by $E_p(t)$ and show that its Fourier transform $E_{p\omega}(\omega)$ does not have zeros for finite and real $\omega$ when $0 < |\sigma| < \frac{1}{2}$, corresponding to the critical strip \textbf{excluding} the critical line.\\

In Section~\ref{sec:Section_3}, it is shown that the new method is \textbf{not} applicable to Hurwitz zeta function and related functions and  \textbf{does not} contradict the existence of their non-trivial zeros away from the critical line with real part of $s = \frac{1}{2}$. \\ 


We present an \textbf{outline} of the new method below.


.\\ \tocless\subsection{\label{sec:Section_1_1} \textbf{Step 1: Inverse Fourier Transform of $\xi(\frac{1}{2} + i \omega)$ } \protect\\  \lowercase{} }

Let us start with Riemann's Xi Function $\xi(s)$ evaluated at $s = \frac{1}{2} + i \omega$ given by $\xi(\frac{1}{2} + i \omega)= \Xi(\omega) =  E_{0\omega}(\omega)$, where $\omega$ is real. Its inverse Fourier Transform is given by $ E_0(t)= \frac{1}{2 \pi}  \int_{-\infty}^{\infty} E_{0\omega}(\omega) e^{i\omega t} d\omega $, where $\omega, t$ are real, as follows \href{https://web.archive.org/web/20240218091913/https://www.ams.org/notices/200303/fea-conrey-web.pdf#page=5}{(link)}.{\citep{ECT}} \href{https://www.ocf.berkeley.edu/~araman/files/math_z/titchmarsh_pp254_255.pdf}{(Titchmarsh pp254-255)} 
We take the term $e^{\frac{t}{2}}$ out of the bracket and rearrange the terms as follows.
\begin{align}\label{sec_intro_eq_1}
E_{0}(t) = \Phi(t) = 2 \sum_{n=1}^{\infty}  [ 2 n^{4}  \pi^{2} e^{\frac{9t}{2}}    - 3 n^{2} \pi  e^{\frac{5t}{2}} ]  e^{- \pi n^{2} e^{2t}} =     \sum_{n=1}^{\infty}  [ 4 \pi^{2} n^{4} e^{4t}    - 6 \pi n^{2}   e^{2t} ]  e^{- \pi n^{2} e^{2t}} e^{\frac{t}{2}}   
\end{align}
We see that $E_{0}(t)=E_{0}(-t)$ is a real and \textbf{even} function of $t$, given that  $E_{0\omega}(\omega) = E_{0\omega}(-\omega)$ because $\xi(s)=\xi(1-s)$  \href{https://www.ocf.berkeley.edu/~araman/files/math_z/Ellison_p147-152.pdf#page=6}{(link)} and hence $\xi(\frac{1}{2} + i \omega)=\xi(\frac{1}{2} - i \omega)$ when evaluated at $s = \frac{1}{2} + i \omega$.(Details in ~\ref{sec:appendix_C_7})  \\

The inverse Fourier Transform of   $\xi(\frac{1}{2}+ \sigma + i \omega) = E_{p\omega}(\omega)$ is given by the real function $E_p(t)$.  We can write $E_p(t)$ as follows for $0 < |\sigma| < \frac{1}{2}$ and this is shown in detail in ~\ref{sec:appendix_A}. 
\begin{align}\label{sec_intro_eq_2}
E_{p}(t) = E_{0}(t) e^{-\sigma t} =    \sum_{n=1}^{\infty}  [ 4 \pi^{2} n^{4} e^{4t}    - 6 \pi n^{2}   e^{2t} ]  e^{- \pi n^{2} e^{2t}} e^{\frac{t}{2}} e^{-\sigma t} 
\end{align}
We can see that $E_{p}(t)$ is an analytic function for real $t$, given that the sum and product of exponential functions are analytic for real $t$ and hence infinitely differentiable for real $t$.

.\\ \tocless\subsection{\label{sec:level2} \textbf{Step 2: On the zeros of a related function $G(\omega, t_0)$ } \protect\\  \lowercase{} }

\textbf{Statement 1}: Let us assume that Riemann's Xi function $\xi(\frac{1}{2} + \sigma + i \omega)= E_{p\omega}(\omega)$ has a zero at $\omega = \omega_{0}$ where $\omega_{0}$ is real and $0 < |\sigma| < \frac{1}{2}$, corresponding to the critical strip excluding the critical line.  We will prove that this assumption leads to a \textbf{contradiction}. \\

Let us consider $0 < \sigma < \frac{1}{2}$ at first. Let us consider a new function  $g(t, t_0) =   f(t, t_0) e^{-\sigma t}  u(-t) +  f(t, t_0)  e^{\sigma t}  u(t)  $, where $f(t, t_0) = e^{- \sigma t_0}   E_0(t + t_0 )  e^{-\sigma t}  + e^{\sigma t_0}  E_0(t - t_0 ) e^{-\sigma t}  + 2 \cosh{( \sigma t_0)} E_p(t) $ and $f_1(t, t_0)=e^{\sigma t_0} E_{p}(t + t_0)$ and $f_2(t, t_0)= f_1(t, -t_0)  $ and $t_0$ is real  and $g(t, t_0)$ is a real function of variable $t$ and $u(t)$ is Heaviside unit step function. We can see that $g(t, t_0) h(t) =  f(t, t_0)$ where $h(t)=[ e^{ \sigma t} u(-t) + e^{ -  \sigma t} u(t) ]$ . \\

In Section~\ref{sec:Section_2_1}, we show that the Fourier transform of the \textbf{even function} $g_{even}(t, t_0)=\frac{1}{2} [g(t, t_0)+g(-t, t_0) ] $ given by $ G_{R}(\omega, t_0)$  must have \textbf{at least one zero} at $\omega = \omega_{z}(t_0) \neq 0$, for every value of $t_0 \in \Re$, where $G_{R}(\omega, t_0)$ crosses the zero line to the opposite sign, to satisfy Statement 1, where $\omega_{z}(t_0)$ is real.

.\\ \tocless\subsection{\label{sec:level2} \textbf{Step 3: On the zeros of the function $G_{R}(\omega, t_0)$  } \protect\\  \lowercase{} }
In Section~\ref{sec:Section_2_1_a}, we compute the Fourier transform of the function $g(t, t_0)$ and compute its real part given by $G_{R}(\omega, t_0)$ and we can write as follows. 
\begin{align} \label{sec:sec_1_2_eq_1} 
G_{R}(\omega, t_0) =  \int_{-\infty}^{0}  E_0(\tau + t_0 )  ( e^{- \sigma t_0} e^{-2 \sigma \tau} + e^{\sigma t_0})]  \cos{(\omega \tau)} d\tau \notag\\
+ \int_{-\infty}^{0}    E_0(\tau - t_0  )  ( e^{\sigma t_0} e^{-2 \sigma \tau} + e^{- \sigma t_0})   \cos{(\omega \tau)} d\tau  
+ 2 \cosh{( \sigma t_0)} \int_{-\infty}^{0}  E_0(\tau) ( e^{-2 \sigma \tau} + 1) \cos{(\omega \tau )} d\tau 
\end{align}
We require $G_{R}(\omega, t_0) =0$ for $\omega=\omega_{z}(t_0)$ for every value of $t_0$, to satisfy \textbf{Statement 1}. In general $\omega_{z}(t_0) \neq \omega_0$. Hence we can see that $P(t_0)=G_{R}(\omega_{z}(t_0), t_0) = 0$.

.\\ \tocless\subsection{\label{sec:level2} \textbf{Step 4: Zero Crossing function $\omega_z(t_0)$ is an \textbf{even} function of variable $t_0$ } \protect\\  \lowercase{} }

In Section~\ref{sec:Section_A_1_2}, we show the result in Eq.~\ref{sec:sec_1_2_eq_1_1} and that $\omega_{z}(t_0)=\omega_z( -t_0)$. It is shown that  $P(t_0) = G_{R}(\omega_z(t_0), t_0)=P_{odd}(t_0) +P_{odd}( -t_0) = 0$ and that $P_{odd}(t_0) $ is an \textbf{odd} function of $t_0$, as follows. 
\begin{align}  \label{sec:sec_1_2_eq_1_1} 
P_{odd}(t_0) =  \int_{-\infty}^{0}  E_0(\tau + t_0 )  ( e^{- \sigma t_0} e^{-2 \sigma \tau} + e^{\sigma t_0})]  \cos{(\omega_{z}(t_0) \tau)} d\tau +  \cosh{( \sigma t_0)} \int_{-\infty}^{0}  E_0(\tau) ( e^{-2 \sigma \tau} + 1) \cos{(\omega_{z}(t_0) \tau )} d\tau
\end{align}

.\\ \tocless\subsection{\label{sec:level2} \textbf{Step 5: Final Step  } \protect\\  \lowercase{} }

In Section~\ref{sec:Section_A_1_8_1}, it is shown that $\omega_z(t_0)$ is a \textbf{continuous}  function of variable $t_0$, for $t_0 \in \Re$. In Section~\ref{sec:Section_A_1_6}, it is shown that $E_0(t)$ is \textbf{strictly decreasing} for $t > 0$.\\

In Section~\ref{sec:Section_A_1_3}, we set $t_0=t_{0c}$, such that $\omega_z(t_{0c}) t_{0c}=  \pi  $ and substitute in the equation for $P_{odd}(t_0)$ in Eq.~\ref{sec:sec_1_2_eq_1_1}  and show the result in Eq.~\ref{sec:sec_1_2_eq_2}. 
\begin{align} \label{sec:sec_1_2_eq_2}    
A(t_{0c}) = \int_{0}^{t_{0c}}     E_0(\tau ) (  \sinh{(\sigma  t_{0c} )} - \sinh{( 2 \sigma  \tau - \sigma  t_{0c} )} )  \cos{(\omega_z(t_{0c}) \tau)} d\tau  = 0 
\end{align} 
We show that \textbf{each} of the terms in the integrand in Eq.~\ref{sec:sec_1_2_eq_2} are \textbf{positive}, in the interval $0 < \tau < t_{0c}$ where $t_{0c} >0$. \\

Hence the result in Eq.~\ref{sec:sec_1_2_eq_2}  leads to a \textbf{contradiction} for $0 < \sigma < \frac{1}{2}$.\\

We show this result for $0 < \sigma < \frac{1}{2}$ and then use the property $\xi(\frac{1}{2} + \sigma + i \omega) = \xi(\frac{1}{2} - \sigma - i \omega)$ to show the result for $-\frac{1}{2} < \sigma < 0$. Hence we produce a \textbf{contradiction} of  \textbf{Statement 1} that the Fourier Transform of the function $E_p(t) = E_0(t) e^{-\sigma t} $ has a zero at $\omega = \omega_{0}$ for  $0 < |\sigma| < \frac{1}{2}$.

.\\ \tocless\section{\label{sec:Section_2} An Approach towards Riemann's Hypothesis  \protect\\  \lowercase{} }

\textbf{Theorem 1}:   Riemann's Xi function $\xi(\frac{1}{2} + \sigma + i \omega)= E_{p\omega}(\omega)$ does not have zeros for any real value of $\omega$, for $0 < |\sigma| < \frac{1}{2}$, corresponding to the critical strip excluding the critical line, given that $E_0(t) = E_0(-t)$ is an even function of variable $t$, where $E_p(t) =  \frac{1}{2 \pi} \int_{-\infty}^{\infty} E_{p\omega}(\omega) e^{i \omega t} d\omega $, $E_p(t)= E_0(t) e^{-\sigma t}$ and $E_{0}(t) =    \sum_{n=1}^{\infty}  [ 4 \pi^{2} n^{4} e^{4t}    - 6 \pi n^{2}   e^{2t} ]  e^{- \pi n^{2} e^{2t}} e^{\frac{t}{2}} $.\\

\textbf{Proof}: We assume that Riemann Hypothesis is false and prove its truth using proof by contradiction.\\

\textbf{Statement 1}: Let us assume that Riemann's Xi function $\xi(\frac{1}{2} + \sigma + i \omega)= E_{p\omega}(\omega)$ has a zero at $\omega = \omega_{0}$ where $\omega_{0}$ is real and $0 < |\sigma| < \frac{1}{2}$, corresponding to the critical strip excluding the critical line.  We will prove that this assumption leads to a \textbf{contradiction}. \\

We will prove it for $0 < \sigma < \frac{1}{2}$ first and then use the property $\xi(\frac{1}{2} + \sigma + i \omega) = \xi(\frac{1}{2} - \sigma - i \omega)$ to show the result for $-\frac{1}{2} < \sigma < 0$ and hence show the result for  $0 < |\sigma| < \frac{1}{2}$.\\

We know that $\omega_{0} \neq 0$, because $\zeta(s)$ has no zeros on the real axis between 0 and 1, when $s=\frac{1}{2}+\sigma + i \omega$ is real, $\omega=0$ and $0 \leq |\sigma| < \frac{1}{2}$. {\citep{ECT}} \href{https://www.ocf.berkeley.edu/~araman/files/math_z/Titchmarsh_pp30_31.pdf}{(Titchmarsh pp30-31)}. This is shown in detail in first two paragraphs in ~\ref{sec:appendix_C_1}.

.\\ \tocless\subsection{\label{sec:Section_2_1} \textbf{New function $g(t, t_0)$} \protect\\  \lowercase{} }

We consider the function $E_{p}(t)=  E_0(t)  e^{-\sigma t}$ whose Fourier transform $E_{p\omega}(\omega)$ has a zero at  $\omega = \omega_{0}$, using Statement 1. We consider the function $f_1(t, t_0)= e^{\sigma t_0} E_{p}(t + t_0)$ where $t_0 \in \Re$ is real, and its Fourier transform is given by $F_1(\omega, t_0)  = E_{p\omega}(\omega)  e^{\sigma t_0}  e^{i \omega t_0} $ using linearity and time shift properties of the Fourier transform ( \href{https://archive.is/CLVqX}{link}). We consider $f_2(t, t_0) = f_1(t,-t_0) $ and its Fourier transform is given by $F_2(\omega, t_0)  = F_1(\omega, -t_0)   =   E_{p\omega}(\omega)    e^{-\sigma t_0} e^{-i \omega t_0}$. (We call this \textbf{Result 2.1.a})\\ 

We consider $f(t, t_0)=  e^{- \sigma t_0} f_1(t, t_0) +  e^{\sigma t_0}  f_2(t, t_0) + 2 \cosh{( \sigma t_0)} E_p(t)$ and the Fourier Transform of this function $F(\omega, t_0)  = e^{- \sigma t_0} F_1(\omega, t_0)  +  e^{\sigma t_0}  F_2(\omega, t_0) + 2 \cosh{( \sigma t_0)} E_{p\omega}(\omega) =  E_{p\omega}(\omega) [ e^{- \sigma t_0} e^{\sigma t_0}  e^{i \omega t_0} + e^{\sigma t_0} e^{-\sigma t_0}  e^{-i \omega t_0} + 2 \cosh{( \sigma t_0)} ] = E_{p\omega}(\omega) [  e^{i \omega t_0} +  e^{-i \omega t_0} + 2 \cosh{( \sigma t_0)} ]$, using Result 2.1.a and $e^{\sigma t_0} e^{-\sigma t_0} = 1$. Hence $ F(\omega, t_0)  =  E_{p\omega}(\omega) [  2 \cos{( \omega t_0)} + 2 \cosh{( \sigma t_0)} ] $ also has a zero at the \textbf{same} $\omega = \omega_{0}$, using Statement 1. ( \textbf{Result 2.1.b}).\\

We can write the above equations as follows. We use $e^{\sigma t_0} e^{-\sigma t_0} = 1$ below.  We see that $f_2(t, t_0) = f_1(t, -t_0)$ and $f(t, t_0)=f(t, -t_0)$ (\textbf{Result 2.1.c}) as shown below. 
\begin{align} \label{sec:sec_2_1_eq_0}    
E_{p}(t)=  E_0(t) e^{- \sigma t}  \notag\\
f_1(t, t_0)= e^{\sigma t_0} E_{p}(t + t_0)  = e^{\sigma t_0}   E_0(t + t_0 )  e^{-\sigma t} e^{-\sigma t_0} \notag\\
 f_2(t, t_0) =  f_1(t, -t_0)  =   e^{-\sigma t_0}  E_0(t - t_0 )  e^{-\sigma t} e^{\sigma t_0} \notag\\
f(t, t_0)=  e^{- \sigma t_0} f_1(t, t_0) +  e^{\sigma t_0}  f_2(t, t_0) + 2 \cosh{( \sigma t_0)} E_p(t) \notag\\
f(t, t_0)  =  e^{- \sigma t_0}   E_0(t + t_0 )  e^{-\sigma t}  + e^{\sigma t_0}  E_0(t - t_0 ) e^{-\sigma t}  + 2 \cosh{( \sigma t_0)} E_p(t)
\end{align} 
We consider a new function  $g(t, t_0) = f(t, t_0) e^{-\sigma t} u(-t) + f(t, t_0) e^{\sigma t}  u(t)  $ where $g(t, t_0)$ is a real function of variable $t$ and $u(t)$ is Heaviside unit step function. We see that $g(t, t_0)h(t)= f(t, t_0)$ where $h(t)= [ e^{ \sigma t} u(-t) + e^{ - \sigma t} u(t) ] $. We note that $f(0, t_0)=  4  E_0(t) e^{-\sigma t} $ and $g(0, t_0) = f(0, t_0) e^{-\sigma t} u(-t) + f(0, t_0) e^{\sigma t}  u(t)  $ and hence $f(t, t_0)$ and $g(t, t_0)$ are \textbf{non-zero} at $t_0=0$.
\begin{align} \label{sec:sec_2_1_eq_0a0}    
g(t, t_0) = [ f(t, t_0) e^{-\sigma t}] u(-t) + [ f(t, t_0) e^{\sigma t}] u(t)  \notag\\
g(t, t_0)h(t)= f(t, t_0), \quad h(t)=  e^{ \sigma t} u(-t) + e^{ - \sigma t} u(t)  
\end{align}
We can show that $E_p(t),  h(t)$ are absolutely integrable functions and go to zero as $|t| \to  \infty$. Hence their respective Fourier transforms given by $E_{p\omega}(\omega), H(\omega)$ are finite for real $\omega$ and go to zero as $|\omega| \to \infty$, as per Riemann Lebesgue Lemma \href{https://archive.is/HV6zJ}{(link)}. We can show that  $E_0(t)$ and $E_0(t) e^{-2 \sigma t}$ are absolutely \textbf{integrable} functions. These results are shown in ~\ref{sec:appendix_C_1}. \\

In Section~\ref{sec:Section_2_1_a} and Section~\ref{sec:Section_A_1_2}, it is shown that $g(t, t_0)=  e^{- \sigma t_0} g_1(t, t_0) +  e^{\sigma t_0}  g_2(t, t_0) + 2 \cosh{( \sigma t_0)} g_3(t)$ is a Fourier transformable function and its Fourier transform given by $G(\omega, t_0) = e^{- \sigma t_0} G_1(\omega, t_0)  + e^{\sigma t_0} G_2(\omega, t_0) +  2 \cosh{( \sigma t_0)} G_3(\omega) $ converges. (Eq.~\ref{sec:sec_2_1_eq_1_2} and Eq.~\ref{sec:sec_2_1_eq_9a}) \\

We take the Fourier transform of the equation $g(t, t_0)  h(t) = f(t, t_0)$ where $h(t) = e^{\sigma t} u(-t) + e^{- \sigma t} u(t)  $ and get $\frac{1}{2\pi} [ G(\omega, t_0) \ast H(\omega)] = F(\omega, t_0)  =  E_{p\omega}(\omega) [  2 \cos{( \omega t_0)} + 2 \cosh{( \sigma t_0)} ] =  F_{R}(\omega, t_0) + i F_{I}(\omega, t_0) $ as per \textbf{convolution theorem} \href{https://archive.is/AjVin}{(link)}, and Result 2.1.b, where $\ast$ denotes convolution operation given by $F(\omega, t_0) = \frac{1}{2\pi} \int_{-\infty}^{\infty} G(\omega', t_0) H(\omega - \omega') d\omega'$. (\textbf{Result 2.1.d})\\

We see that $H(\omega)= H_{R}(\omega) = [ \frac{1}{  \sigma - i \omega} +  \frac{1}{  \sigma + i \omega}   ]  =  \frac{2 \sigma}{(\sigma^{2} + \omega^{2})}  $ is real and is the Fourier transform of the function $h(t)$ \href{https://web.archive.org/web/20240218075832/https://engineering.purdue.edu/~mikedz/ee301/FourierTransformTable.pdf#page=3}{(link)}.  $G(\omega, t_0)=  G_{R}(\omega, t_0) + i G_{I}(\omega, t_0)$ is the Fourier transform of the function $g(t, t_0)$. We can write $g(t, t_0)= g_{even}(t, t_0) + g_{odd}(t, t_0)$ where $g_{even}(t, t_0)$ is an even function and $g_{odd}(t, t_0)$ is an odd function of variable $t$. \\

If Statement 1 is true, then we require the Fourier transform of the function $f(t, t_0)$ given by $F(\omega, t_0)$ to have a zero at $\omega = \omega_{0}$ for \textbf{every value} of $t_0$, using Result 2.1.b. This implies that the \textbf{real} part of the Fourier transform of the \textbf{even function} $g_{even}(t, t_0)=\frac{1}{2} [g(t, t_0)+g(-t, t_0) ] $  given by $G_{R}(\omega, t_0)$(Details in ~\ref{sec:appendix_I_3} and \href{https://archive.is/CLVqX}{link}) must have \textbf{at least one zero} at $\omega = \omega_{z}(t_0) \neq 0$ where $\omega_{z}(t_0)$ is real, where $G_{R}(\omega, t_0)$ crosses the zero line to the opposite sign, explained below. We note that $\omega_{z}(t_0)$ can be different from $\omega_0$ in general. \\

Because $H(\omega) = \frac{2 \sigma}{(\sigma^{2} + \omega^{2})}$ is real and does not have zeros for any finite value of $\omega$, \textbf{if} $G_{R}(\omega, t_0)$ does not have at least one zero for some $\omega  = \omega_{z}(t_0) \neq 0$, where $G_{R}(\omega, t_0)$ crosses the zero line to the opposite sign, \textbf{then} the \textbf{real part} of $F(\omega, t_0)$ given by $F_{R}(\omega, t_0)= \frac{1}{2 \pi}  [ G_{R}(\omega, t_0) \ast H(\omega)]$, obtained by the convolution of $H(\omega)$ and $G_{R}(\omega, t_0)$, \textbf{cannot} possibly have zeros for any non-zero finite value of $\omega$, which contradicts Result 2.1.b and  \textbf{Statement 1}. This is shown in detail in Lemma 1. \\

The proof for Lemma 1 below is shown for \textbf{a fixed value} of $t_0=t_{0f}$, for $t_0 \in \Re$ (\textbf{Interval A}). The proof continues to hold for our choice of \textbf{each} \textbf{fixed value} of  $t_0$ in interval A.\\

\textbf{Lemma 1:} Let $t_0 \in \Re$ be fixed value and $\xi(\frac{1}{2} + \sigma + i \omega)= E_{p\omega}(\omega)$ has a zero at $\omega=\omega_0$ using Statement 1. Then the real part of the Fourier transform of the \textbf{even function} $g_{even}(t, t_0) $  given by $G_{R}(\omega, t_0)$  must have \textbf{at least one zero} at $\omega = \omega_{z}(t_0) \neq 0$, where $G_{R}(\omega, t_0)$ crosses the zero line to the opposite sign and $\omega_{z}(t_0)$ is real.\\

\textbf{Proof}: If $E_{p\omega}(\omega)$ has a zero at $\omega=\omega_0$  to satisfy Statement 1, then $F(\omega, t_0)$  has a zero at $\omega=\omega_0$, using Result 2.1.b and  its real part given by $F_{R}(\omega, t_0)$ has a zero at $\omega=\omega_0$ (\textbf{Result 2.1.e}). We see that $\omega_{0} \neq 0$ in para 5 of Section~\ref{sec:Section_2}.\\

We do not have a closed form solution for $G_{R}(\omega, t_{0})$ and do not know the exact location of its zeros at $\omega = \omega_z( t_0) $. For a specific choice of $t_0$, \textbf{only one} of the 2 cases is possible: \\ \textbf{Case A:} $G_{R}(\omega, t_{0})$ does not have a zero crossing for any choice of $\omega \neq 0$ or \\
\textbf{Case B: }$G_{R}(\omega, t_{0})$ has at least one zero crossing for a specific $\omega \neq 0$. \\ 
\textbf{If} Statement 1 is true, \textbf{then} Case B is the \textbf{only} possibility and Case A is \textbf{ruled out}, as shown below.\\

We want to show the \textbf{Result 2.1.h} that $G_{R}(\omega, t_{0})$ \textbf{must have at least one} zero crossing at \textbf{some value} of $\omega = \omega_z( t_0) \neq 0$ (\textbf{Case B}), to satisfy \textbf{Statement 1}, for this choice of fixed $t_0$. \\

To show Result 2.1.h, we \textbf{assume the opposite Case A}, that $G_{R}(\omega, t_0)$ \textbf{does not} have at least one zero for \textbf{any} value of $\omega  \neq 0$, where $G_{R}(\omega, t_0)$  crosses the zero line to the opposite sign (zero crossing) and will show that $F_{R}(\omega, t_0)$ does not have at least one zero at finite $\omega \neq 0$ for this case, which \textbf{contradicts} Result 2.1.e and Statement 1 and hence we \textbf{rule out} Case A and arrive at Case B, in the following Proof of Lemma 1. (Case B is the same as Result 2.1.h).\\

This \textbf{does not} mean that, proof of Lemma 1 will work \textbf{only if} $G_{R}(\omega, t_{0})$  does not have a zero  crossing for any value of $\omega \neq 0$, for any choice of $t_0$. The device \textbf{Proof by Contradiction} is used here to \textbf{rule out} Case A and arrive at Case B. (Details of other cases in Section~\ref{sec:Section_2_1_a1} and Section~\ref{sec:Section_3_1_a0}) \\

It is noted that, for \textbf{Case B}, we \textbf{do not} use Eq.~\ref{app_D_5_eq_1} to Eq.~\ref{app_D_5_eq_2} and related arguments, because \\ Case B is the desired Result 2.1.h. (\textbf{\textbf{Result 2.1.f}}) \\

The arguments above and following proof continue to hold for our choice of \textbf{each} \textbf{fixed value} of  $t_0$ in interval A.\\

Given that $H(\omega)$ is real, using Result 2.1.d, we write the convolution theorem only for the real parts as follows.
\begin{align} \label{app_D_5_eq_1}   
F_{R}(\omega, t_0) = \frac{1}{2 \pi}  \int_{-\infty}^{\infty} G_R(\omega', t_0) H(\omega - \omega') d\omega' 
\end{align}
We can show that the above integral converges for real $\omega$, given that the integrand is absolutely integrable because $G_R(\omega', t_0) H(\omega - \omega')$ has a minimum fall-off rate of $\frac{1}{\omega^2}$ as $|\omega| \to \infty$ because the first derivative of $h(t)$ is discontinuous at $t=0$.(Details in ~\ref{sec:appendix_C_2} and ~\ref{sec:appendix_C_5z})\\

We substitute $H(\omega) = \frac{2 \sigma}{(\sigma^{2} + \omega^{2})}$ in Eq.~\ref{app_D_5_eq_1}  and we get
\begin{align} \label{app_D_51_eq_1_1}   
F_{R}(\omega, t_0) = \frac{\sigma}{\pi}  \int_{-\infty}^{\infty} G_R(\omega', t_0) \frac{1}{(\sigma^{2} + (\omega - \omega')^{2})}  d\omega'
\end{align}
We can split the integral in Eq.~\ref{app_D_51_eq_1_1} using $\int_{-\infty}^{\infty}= \int_{-\infty}^{0} + \int_{0}^{\infty}$, as follows.
\begin{align} \label{app_D_5_eq_1_2}   
F_{R}(\omega, t_0)  = \frac{\sigma}{ \pi}  [ \int_{-\infty}^{0} G_R(\omega', t_0) \frac{1}{(\sigma^{2} + (\omega - \omega')^{2})}  d\omega' + \int_{0}^{\infty} G_R(\omega', t_0) \frac{1}{(\sigma^{2} + (\omega - \omega')^{2})}  d\omega' ]
\end{align}
We see that $G_{R}(-\omega, t_0)= G_{R}(\omega, t_0)$ because $g(t, t_0)$ is a real function of variable $t$. (Details in ~\ref{sec:appendix_I_2} and \href{https://archive.is/CLVqX}{link}) We substitute $\omega' = -\omega''$ and $d\omega' = -d\omega''$ in the first integral in Eq.~\ref{app_D_5_eq_1_2} and substituting $\omega'' = \omega'$ in the result, we can write as follows.
\begin{align}\label{app_D_5_eq_2}   
F_{R}(\omega, t_0) = \frac{\sigma}{\pi}   \int_{0}^{\infty} G_R(\omega', t_0) [  \frac{1}{(\sigma^{2} + (\omega + \omega')^{2})} + \frac{1}{(\sigma^{2} + (\omega - \omega')^{2})} ]  d\omega'  
\end{align}
We note that $G_{R}(\omega',  t_{0})$ is a function of $\omega'$ \textbf{only}, for a given choice of $t_0$ and the integrand in Eq.~\ref{app_D_5_eq_2} is integrated over the variable $\omega'$ \textbf{only}.\\

We see that $G_R(\omega', t_0)$ converges for real $\omega'$ (Eq.~\ref{sec:sec_2_1_eq_1_2} and Eq.~\ref{sec:sec_2_1_eq_9}). As $\omega' \to \infty$, the term $ \frac{1}{(\sigma^{2} + (\omega + \omega')^{2})} + \frac{1}{(\sigma^{2} + (\omega - \omega')^{2})}$ and the integrand in Eq.~\ref{app_D_5_eq_2} goes to zero. For finite $\omega \geq 0$ and finite $\omega' \geq 0$, we can see that the term $ \frac{1}{(\sigma^{2} + (\omega + \omega')^{2})} + \frac{1}{(\sigma^{2} + (\omega - \omega')^{2})} = \frac{2(\sigma^{2} + \omega^{2} + \omega'^{2} ) }{(\sigma^{2} + (\omega + \omega')^{2}) (\sigma^{2} + (\omega - \omega')^{2})} > 0$, for $0 < \sigma < \frac{1}{2}$. We see that $G_R(\omega', t_0)$ is \textbf{not} an all zero function of variable $\omega'$ (Details in Section~\ref{sec:Section_2_1_a0}). (\textbf{Result 2.1.g})\\

$\bullet$ \textbf{\textbf{Case 1:} $G_R(\omega', t_0) \geq 0$ for all finite $\omega' \geq 0$} \\

We see that $F_{R}(\omega, t_0) > 0$ for all finite $\omega > 0$, using Result 2.1.g. We see that $F_{R}(-\omega, t_0)= F_{R}(\omega, t_0)$ because $f(t, t_0)$ is a real function ( ~\ref{sec:appendix_I_2}) and \href{https://archive.is/CLVqX}{link} ). Hence $F_{R}(\omega, t_0) > 0$ for all finite $\omega < 0$.\\

 This \textbf{contradicts}  Statement 1 and Result 2.1.e which requires $F_{R}(\omega, t_0)$ to have at least one zero at finite $\omega \neq 0$. Therefore $G_R(\omega', t_0)$ must have\textbf{ at least one zero} at $\omega' =  \omega_{z}(t_0) > 0$ where it crosses the zero line and becomes negative, where $\omega_{z}(t_0)$ is finite. \\

$\bullet$ \textbf{\textbf{Case 2:} $G_R(\omega', t_0) \leq 0$ for all finite  $\omega' \geq 0$} \\

We see that  $F_{R}(\omega, t_0) < 0$ for all finite $\omega > 0$, using Result 2.1.g. We see that $F_{R}(-\omega, t_0)= F_{R}(\omega, t_0)$ because $f(t, t_0)$ is a real function ( ~\ref{sec:appendix_I_2}) and \href{https://archive.is/CLVqX}{link} ). Hence $F_{R}(\omega, t_0) < 0$ for all finite $\omega < 0$.\\

This \textbf{contradicts} Statement 1 and Result 2.1.e which requires $F_{R}(\omega, t_0)$ to have at least one zero at finite $\omega \neq 0$. Therefore $G_R(\omega', t_0)$ must have\textbf{ at least one zero} at $\omega' =  \omega_{z}(t_0) > 0$ , where it crosses the zero line and becomes positive, where $\omega_{z}(t_0)$ is finite. \\

We have shown that, $G_{R}(\omega, t_0)$ must have\textbf{ at least one zero} at finite $\omega =  \omega_{z}(t_0) \neq 0$ where it crosses the zero line to the opposite sign, to satisfy \textbf{Statement 1}, for specific choices of fixed $t_0$. We call this \textbf{Result 2.1.h}. \\

\textbf{Points to be Noted:}\\

$\bullet$ The arguments above and the proof continue to hold for our choice of \textbf{each} \textbf{fixed value} of  $t_0$ in Interval A ($t_0 \in \Re$). \\

$\bullet$ We consider the case where zeros of  $G_{R}(\omega, t_0)$  come from $t_0$ alone, for specific choice of $t_0$. For this case, $G_{R}(\omega, t_0)=0$ for all $\omega$, for that specific choice of $t_0$. But we see that $G_{R}(\omega, t_0)$ is \textbf{not} an all zero function of variable $\omega$, as shown in Section~\ref{sec:Section_2_1_a0}. It is shown in Note 1 in Section~\ref{sec:Section_C_6_0z} that the zeros of $G_{R}(\omega, t_0)$  are \textbf{isolated}. This \textbf{rules out} the case where zeros of  $G_{R}(\omega, t_0)$  come from $t_0$ alone. This means that, if Statement 1 is true, then $G_{R}(\omega, t_0)$  must have at least one zero crossing at $\omega= \omega_z(t_0)$, \textbf{due to a term} related to $\omega$.\\

$\bullet$ Note that $\omega_z(t_0)$ is a \textbf{dependent} function of $t_0$, given that $G_{R}(\omega, t_{0})$ may be different for different choices of $t_0$ and hence the zero crossing point $\omega_z(t_0)$ may be different for different choices of $t_0$. \\

This means that, we fix $t_0 = t_{0f}$ in Pass 1 and use Lemma 1 to show that $\omega_z(t_{0f}) \neq 0$ for \textbf{that specific} choice of $t_0$ and $G_{R}(\omega, t_{0f})$. In Pass 2, we use a \textbf{different choice} of $ t_0 = t_{0f}^{'}$ and use Lemma 1 to show that $\omega_z(t_{0f}^{'}) \neq 0$ for that specific \textbf{different} choice of $t_0$ and $G_{R}(\omega, t_{0f}^{'})$. In general,  $G_{R}(\omega, t_{0f}) \neq G_{R}(\omega, t_{0f}^{'})$ and hence the zero crossing point $\omega_z(t_{0f}) \neq \omega_z(t_{0f}^{'})$. This argument holds for Pass 1, 2 to Pass N, for  \textbf{each} value of $t_0$ in Interval A. \\

In general, $\omega_z(t_{0f}) \neq \omega_z(t_{0f}^{'})$ in Lemma 1. But the specific case of $\omega_z( t_0)$ equals a constant for all $t_0$ is allowed in Lemma 1 and this constant $\omega_z(t_0)$ is a \textbf{continuous function} of  $t_0$ and there is \textbf{no need} for using Implicit Function Theorem (IFT) in Section~\ref{sec:Section_A_1_8_1}. \\

$\bullet$ We do not have a closed form solution for $G_{R}(\omega, t_{0})$ and do not know the exact location of its zeros at $\omega = \omega_z( t_0) $, for each fixed choice of $t_0$. This problem \textbf{exists} in the case of $\xi(\frac{1}{2} + i \omega)= \Xi(\omega) =  E_{0\omega}(\omega)$ \textbf{also} (see Eq.~\ref{sec:App_H_eq_3_1} in ~\ref{sec:appendix_H}) given that we \textbf{do not} have a closed form solution for $ \Xi(\omega) $ and \textbf{do not} know the exact location of its zeros. \textbf{Despite} this problem, Hardy and Littlewood were able to prove that infinitely many of the zeros of $\xi(s)$ and $\zeta(s)$ are on the critical line with real part of $s=\frac{1}{2}$.{\citep{GH}} \\

Similarly, it is shown in proof of Lemma 1 that $G_{R}(\omega, t_0)$ must have\textbf{ at least one zero} at finite $\omega =  \omega_{z}(t_0) \neq 0$ where it crosses the zero line to the opposite sign, to satisfy \textbf{Statement 1} and this is sufficient to prove Theorem 1, \textbf{despite} this problem with $G_{R}(\omega, t_{0})$.\\

$\bullet$ In the rest of the sections, we consider only the \textbf{first} zero crossing to the right of origin, where $G_{R}(\omega, t_0)$  crosses the zero line to the opposite sign. Hence $0 < \omega_z(t_0) < \infty$, for $t_0 \in \Re$, to satisfy \textbf{Statement 1}. 

.\\ \tocless\subsubsection{\label{sec:Section_2_1_a1} \textbf{Discussion of Lemma 1 } \protect\\  \lowercase{} }

\textbf{Result 2.1.h:} $G_{R}(\omega, t_0)$ must have\textbf{ at least one zero} at finite $\omega =  \omega_{z}(t_0) \neq 0$ where it crosses the zero line to the opposite sign, to satisfy \textbf{Statement 1}.\\

For each \textbf{fixed} value of $t_0$, only 2 cases are possible for $G_{R}(\omega, t_{0})$. \textbf{Case A:} $G_{R}(\omega, t_{0})$ does not have a zero crossing for any choice of $\omega \neq 0$.  \textbf{Case B:} $G_{R}(\omega, t_{0})$ has at least one zero crossing for a specific $\omega \neq 0$. Proof of Lemma 1 assumes Case A and uses \textbf{Proof by Contradiction} to rule out Case A and arrive at Case B, for each choice of fixed $t_0$. This does not mean that Proof of Lemma 1 does not work for Case B. For Case B, we \textbf{do not} use Proof of Lemma 1 and jump to the end of the proof because we already have the desired Result 2.1.h which is the same as Case B.\\

The logic used is this proof is as follows: \textbf{If} Statement 1 is true(RH is false), \textbf{then} Result 2.1.h is true (Case B), for \textbf{each} \textbf{fixed} value of $t_0$ in interval A ($t_0 \in \Re$) and hence Case A is \textbf{ruled out} and only Case B is possible for $G_{R}(\omega, t_0)$. Then we proceed with Result 2.1.h to Section 2.3, 2.4 and Section 3, to produce a \textbf{contradiction} of Statement 1 in Eq.~\ref{sec:sec_az3_1_3_0_0_eq_3_g1_0}  and thus prove the truth of RH.

.\\ \tocless\subsubsection{\label{sec:Section_3_1_a0}  \textbf{Discussion of all possible cases of  $G_{R}(\omega, t_{0})$} \protect\\  \lowercase{} }

We analyze all possible cases of  $G_{R}(\omega, t_{0})$ below. We can arrive at Result 2.1.h in \textbf{Lemma 1} in Section~\ref{sec:Section_2_1}, for \textbf{each} \textbf{fixed} value of $t_0$ in interval A ($t_0 \in \Re$ ), using Proof of Lemma 1 for Case C and Case D or using Case E, as explained below. \\

It is noted that $F_{R}(\omega, t_0)$ and $G_{R}(\omega, t_0)$ may have more zeros than $F(\omega, t_0)$ and $G(\omega, t_0)$ respectively. That \textbf{does not} affect the proof of Lemma 1, as explained below.\\

We consider 3 possible cases of $G_{R}(\omega, t_{0})$ below.\\

$\bullet$  \textbf{Case C:} We consider the case that $G_{R}(\omega, t_{0})$ \textbf{does not} have a zero crossing, for any value of  $\omega \neq 0$, for \textbf{each and every} choice of $ t_{0}$ in Interval A and we assume Statement 1 and use Proof of Lemma 1 for each and every choice of $ t_{0}$, to show that it leads to a \textbf{contradiction} of Statement 1, and hence prove Result 2.1.h, for each and every choice of $ t_{0}$. \\

Hence Case C is \textbf{ruled out}, \textbf{if} Statement 1 is true. (\textbf{Result 2.1.2.c})\\

$\bullet$ \textbf{Case D:} We consider the case  $G_{R}(\omega, t_0)$ has at least one zero crossing at finite $\omega = \omega_z( t_0) \neq 0$,  for \textbf{specific} choices of $t_0=t_{0}^{'}$, but \textbf{not} for all possible choices of $ t_{0}$. \\

For Case D, this means that $G_{R}(\omega, t_{0})$ has \textbf{at least one zero crossing} at $\omega = \omega_z( t_0) $, for \textbf{specific} choices of $t_0=t_{0}^{'}$, which is the desired \textbf{Result 2.1.h} and hence we \textbf{do not}  go through the arguments in this proof and we can jump to end of Proof of Lemma 1 (using Result 2.1.f in Proof of Lemma 1). In this case, we \textbf{have not} assumed Statement 1 and yet arrived at Result 2.1.h, for \textbf{specific} choices of $t_0=t_{0}^{'}$. \\

For Case D, there is \textbf{at least one} choice of $t_0=t_{0f}$ for which $G_{R}(\omega, t_{0})$ \textbf{does not} have a zero crossing, for any value of  $\omega \neq 0$. For this choice of $t_0=t_{0f}$, we assume Statement 1 and use Proof of Lemma 1 to show that it leads to a \textbf{contradiction} of Statement 1, and hence prove Result 2.1.h. \\ 

Hence Case D is \textbf{ruled out}, \textbf{if} Statement 1 is true. (\textbf{Result 2.1.2.d})\\

$\bullet$ \textbf{Case E:} We consider the case  $G_{R}(\omega, t_{0})$ has at least one zero crossing at finite $\omega = \omega_z( t_0) \neq 0$,  for \textbf{each and every} value of $ t_{0}$ in Interval A. This is the same as Result 2.1.h, for \textbf{each and every} value of $ t_{0}$ in Interval A. \\ 

There are \textbf{only 3} possible cases for $G_{R}(\omega, t_{0})$ given by Case C,D and E.  \textbf{If} Statement 1 is true, Case C and Case D have been \textbf{ruled out} using Result 2.1.2.c and Result 2.1.2.d. This leaves only Case E, which is the same as Result 2.1.h. Then we proceed with Result 2.1.h to Section 2.3, 2.4 and Section 3, to produce a \textbf{contradiction} of Statement 1 in Eq.~\ref{sec:sec_az3_1_3_0_0_eq_3_g1_0}.\\

Thus we have produced a \textbf{contradiction} of  \textbf{Statement 1} that the Fourier Transform of the function $E_p(t) = E_0(t) e^{-\sigma t} $ has a zero at $\omega = \omega_{0}$ for  $0 < |\sigma| < \frac{1}{2}$ and hence prove the truth of Riemann's Hypothesis.\\

$\bullet$ We can also consider Case E, \textbf{without} assuming Statement 1. This is the same as Result 2.1.h, for \textbf{each and every} value of $ t_{0}$ in Interval A. Then we proceed with Result 2.1.h to Section 2.3, 2.4 and Section 3, to produce a \textbf{contradiction} of this Case E in Eq.~\ref{sec:sec_az3_1_3_0_0_eq_3_g1_0}. Hence this Case E is \textbf{ruled out}. (\textbf{Result 2.1.2.e})\\

There are \textbf{only 3} possible cases for $G_{R}(\omega, t_{0})$ given by Case C,D and E. We have ruled out Case E using Result 2.1.2.e. This leaves only Case C and Case D for $G_{R}(\omega, t_{0})$. \textbf{If} Statement 1 is true, Case C and Case D have been \textbf{ruled out} using Result 2.1.2.c and Result 2.1.2.d. Given that Case C,D and E are the \textbf{only 3} possible cases for $G_{R}(\omega, t_{0})$, this means \textbf{Statement 1  is false}.

.\\ \tocless\subsection{\label{sec:Section_2_1_a0}  \textbf{$G_{R}(\omega', t_0)$  is not an all zero function of variable $\omega'$} \protect \\  \lowercase{} }

If $G_{R}(\omega', t_0)=0$  for all values of real $\omega'$, for a \textbf{specific} choice of $t_0 \in \Re$ (\textbf{Statement 2}), then $F_{R}(\omega, t_0)$ in Eq.~\ref{app_D_5_eq_1}  is an all zero function of $\omega$, for real $\omega$. Hence $2 f_{even}(t, t_0)= f(t, t_0)+ f(-t, t_0)=0$ for all values of real $t$, given that the Fourier transform of $f_{even}(t, t_0)$ is given by $ F_{R}(\omega, t_0)$, using symmetry properties of Fourier transform(~\ref{sec:appendix_I_3}) and \href{https://archive.is/CLVqX}{link}). Hence $f(t, t_0)=f_{even}(t, t_0)+f_{odd}(t, t_0)=f_{odd}(t, t_0)$ is an \textbf{odd function} of variable $t$ and hence $f(t, t_0)  = -f(-t, t_0)$. (\textbf{Result 2.2}) \\ 

We copy $f(t, t_0)$ from  Eq.~\ref{sec:sec_2_1_eq_0} and write as follows, using $E_p(t)=E_0(t) e^{- \sigma t} $. 
\begin{align} \label{sec_2_2_eq_0a}  
f(t, t_0)  =  e^{- \sigma t_0}   E_0(t + t_0 )  e^{-\sigma t}  + e^{\sigma t_0}  E_0(t - t_0 ) e^{-\sigma t}  + 2 \cosh{( \sigma t_0)} E_0(t) e^{- \sigma t} 
\end{align}
For $t_0 \in \Re$, it is shown that Result 2.2 is false. We will compute $f(t, t_0)$ in Eq.~\ref{sec_2_2_eq_0a} at $t=0$ and show that it does not equal zero. \\

We see that $f(0,  t_0) = e^{- \sigma t_0}   E_0(t_0)   + e^{\sigma t_0} E_0( -t_0 ) + 2 \cosh{( \sigma t_0)} E_0(0)   =   2 \cosh{( \sigma t_0)} E_0( t_0) + 2 \cosh{( \sigma t_0)} E_0(0) =  2 \cosh{(  \sigma t_0)} [ E_0( t_0) + E_0(0) ] $. We use the fact that $E_0(t_0)=E_0(-t_0)$ (Details in ~\ref{sec:appendix_C_7}).\\
 
If Result 2.2 is true, then we require $f(0,  t_0) =   2 \cosh{( \sigma t_0)} [ E_0( t_0) + E_0(0) ] =0$. For our choice of $0 < \sigma < \frac{1}{2}$ and $t_0 \in \Re$, we see that $E_0( t_0)> 0$ and $E_0(0) >0$ using ~\ref{sec:appendix_C_5a} and hence $f(0,  t_0) > 0$ for $t_0 \in \Re$. Hence Result 2.2 is false and Statement 2 is false, for any specific choice of $t_0 \in \Re$ and $G_{R}(\omega', t_0)$  is \textbf{not} an all zero function of variable $\omega'$.\\ 

Hence $G_{R}(\omega', t_0)$  is \textbf{not} an all zero function of variable $\omega'$, for any choice of  $t_0 \in \Re$.

.\\ \tocless\subsection{\label{sec:Section_2_1_a} \textbf{ On the zeros of a related function $G(\omega, t_0)$ } \protect\\  \lowercase{} }

In this section, we compute the Fourier transform of the function $g_{even}(t, t_0)= \frac{1}{2} [ g(t, t_0) + g(-t, t_0) ]$ given by $ G_{R}(\omega, t_0)$(using ~\ref{sec:appendix_I_3}). We require $G_{R}(\omega, t_0) =0$ for $\omega=\omega_{z}(t_0)$ for \textbf{every value} of $t_0$, to satisfy \textbf{Statement 1}, using Lemma 1 in Section~\ref{sec:Section_2_1}.\\

From Eq.~\ref{sec:sec_2_1_eq_0}, we see that $f(t, t_0)=  e^{- \sigma t_0} f_1(t, t_0) +  e^{\sigma t_0}  f_2(t, t_0) + 2 \cosh{( \sigma t_0)} E_p(t) $ and $g(t, t_0)h(t)= f(t, t_0)$ where $h(t)= [ e^{ \sigma t} u(-t) + e^{ - \sigma t} u(t) ]$. We get $g(t, t_0)= \frac{f(t, t_0)}{h(t)} =  e^{- \sigma t_0} \frac{f_1(t, t_0)}{h(t)} +  e^{\sigma t_0}  \frac{f_2(t, t_0)}{h(t)} +  \frac{2 \cosh{( \sigma t_0)} E_p(t)}{h(t)}$. We \textbf{define} $g_1(t, t_0)= \frac{f_1(t, t_0)}{h(t)}$ and $g_2(t, t_0)= \frac{f_2(t, t_0)}{h(t)}$ and $g_3(t)=\frac{ E_p(t)}{h(t)}$ and get $g(t, t_0)=  e^{- \sigma t_0} g_1(t, t_0) +  e^{\sigma t_0}  g_2(t, t_0) + 2 \cosh{( \sigma t_0)} g_3(t)$. (\textbf{Result 2.3.a})\\

We get $g_3(t) =  E_p(t) e^{-\sigma t}  u(-t) +   E_p(t) e^{\sigma t}  u(t)$. We compute the Fourier transform of the function $g_3(t)$ given by $G_3(\omega)= G_{3_R}(\omega) + i G_{3_I}(\omega) $ as follows. We use $E_p(t) = E_0(t) e^{-\sigma t} $. 
\begin{align}\label{sec:sec_2_0_eq_1}    
G_3(\omega) = \int_{-\infty}^{\infty} g_3(t) e^{-i \omega t} dt = \int_{-\infty}^{0} g_3(t) e^{-i \omega t} dt 
+ \int_{0}^{\infty} g_3(t) e^{-i \omega t} dt  \notag\\
G_3(\omega) =   \int_{-\infty}^{0} E_0(t)  e^{-\sigma t} e^{-\sigma t}  e^{-i \omega t} dt 
+ \int_{0}^{\infty}  E_0(t) e^{-\sigma t} e^{\sigma t} e^{-i \omega t} dt   \notag\\  
G_3(\omega) = \int_{-\infty}^{0}  E_0(t) e^{-2 \sigma t}  e^{-i \omega t} dt 
+ \int_{0}^{\infty}   E_0(t) e^{-i \omega t} dt   
\end{align}
We substitute $t = -t$ in the second integral in Eq.~\ref{sec:sec_2_0_eq_1}. We use $-\int_{0}^{-\infty} =  \int_{-\infty}^{0} $.
\begin{align}\label{sec:sec_2_0_eq_1_1} 
G_3(\omega) = \int_{-\infty}^{0}  E_0(t) e^{-2 \sigma t}  e^{-i \omega t} dt 
+ \int_{-\infty}^{0}   E_0(-t) e^{i \omega t} dt 
\end{align}
We use $E_0(t)=E_0(-t)$ from ~\ref{sec:appendix_C_7} and write Eq.~\ref{sec:sec_2_0_eq_1_1}  as follows. The integral in Eq.~\ref{sec:sec_2_0_eq_1_2} converges, given that $E_0(t) e^{-2 \sigma t}$ is an absolutely \textbf{integrable} function. (Details in ~\ref{sec:appendix_C_1}) 
\begin{align}\label{sec:sec_2_0_eq_1_2}   
G_3(\omega) = \int_{-\infty}^{0}  E_0(t) e^{-2 \sigma t}  e^{-i \omega t} dt + \int_{-\infty}^{0}   E_0(t) e^{i \omega t} dt
= G_{3R}(\omega) + i G_{3_I}(\omega)
\end{align}
The above equations can be expanded as follows using the identity $e^{i \omega t } = \cos(\omega t) + i \sin(\omega t)$. Comparing the \textbf{real parts} of $G_{3}(\omega)$, we have
\begin{align}\label{sec:sec_2_0_eq_5}  
G_{3_R}(\omega) =  \int_{-\infty}^{0}  E_0(t) e^{-2 \sigma t}  \cos{(\omega t)} dt 
+ \int_{-\infty}^{0}   E_0(t)   \cos{(\omega t)} dt  = \int_{-\infty}^{0}  E_0(t) ( e^{-2 \sigma t} + 1) \cos{(\omega t)} dt 
\end{align}

.\\ \tocless\subsubsection{\label{sec:Section_2_1_a_z0} \textbf{ Fourier transform of $g_1(t, t_0)$ } \protect\\  \lowercase{} }

Using $g_1(t, t_0)= \frac{f_1(t, t_0)}{h(t)}$ in Result 2.3.a, we get $g_1(t, t_0) = f_1(t, t_0) e^{-\sigma t}  u(-t) +  f_1(t, t_0) e^{\sigma t}  u(t)$, where $f_1(t, t_0)= E_0(t + t_0 )  e^{-\sigma t} $, using Eq.~\ref{sec:sec_2_1_eq_0}.
Then we compute the Fourier transform of the function $g_1(t, t_0)$ given by $G_1(\omega, t_0)= G_{1_R}(\omega, t_0) + i G_{1_I}(\omega, t_0) $ as follows. 
\begin{align}\label{sec:sec_2_1_eq_1}    
G_1(\omega, t_0) = \int_{-\infty}^{\infty} g_1(t, t_0) e^{-i \omega t} dt = \int_{-\infty}^{0} g_1(t, t_0) e^{-i \omega t} dt 
+ \int_{0}^{\infty} g_1(t, t_0) e^{-i \omega t} dt  \notag\\
G_1(\omega, t_0) = \int_{-\infty}^{0}   E_0(t + t_0 )  e^{-\sigma t} e^{-\sigma t}  e^{-i \omega t} dt  
+ \int_{0}^{\infty}    E_0(t + t_0 )  e^{-\sigma t} e^{\sigma t} e^{-i \omega t} dt   \notag\\
G_1(\omega, t_0) = \int_{-\infty}^{0}   E_0(t + t_0 )  e^{-2 \sigma t}  e^{-i \omega t} dt 
+ \int_{0}^{\infty}     E_0(t + t_0 )   e^{-i \omega t} dt   
\end{align}
We substitute $t = -t$ in the second integral in Eq.~\ref{sec:sec_2_1_eq_1}. We use $-\int_{0}^{-\infty} =  \int_{-\infty}^{0} $.
\begin{align}\label{sec:sec_2_1_eq_1_1} 
G_1(\omega, t_0) = \int_{-\infty}^{0}   E_0(t + t_0 )  e^{-2 \sigma t}  e^{-i \omega t} dt  
+ \int_{-\infty}^{0}    E_0(-t + t_0 )  e^{i \omega t} dt 
\end{align}
We use $E_0(t)=E_0(-t)$ from ~\ref{sec:appendix_C_7} and get $E_0(-t + t_0)= E_0(t - t_0 )$ and write Eq.~\ref{sec:sec_2_1_eq_1_1}  as follows.The integral in Eq.~\ref{sec:sec_2_1_eq_1_2} converges, given that $E_0(t) e^{-2 \sigma t}$ is an absolutely \textbf{integrable} function and its $t_0$ shifted versions are absolutely \textbf{integrable}. (Details in ~\ref{sec:appendix_C_1}) 
\begin{align}\label{sec:sec_2_1_eq_1_2}   
G_1(\omega, t_0) = \int_{-\infty}^{0}   E_0(t + t_0 )  e^{-2 \sigma t}  e^{-i \omega t} dt   
+ \int_{-\infty}^{0}    E_0(t - t_0 )   e^{i \omega t} dt 
= G_{1R}(\omega, t_0) + i G_{1_I}(\omega, t_0)
\end{align}
The above equations can be expanded as follows using the identity $e^{i \omega t } = \cos(\omega t) + i \sin(\omega t)$. Comparing the \textbf{real parts} of $G_{1}(\omega, t_0)$, we have
\begin{align}\label{sec:sec_2_1_eq_5}  
G_{1_R}(\omega, t_0) = \int_{-\infty}^{0}  E_0(t + t_0 )  e^{-2 \sigma t}  \cos{(\omega t)} dt 
+ \int_{-\infty}^{0}    E_0(t - t_0 )  \cos{(\omega t)} dt 
\end{align}

.\\ \tocless\subsection{\label{sec:Section_A_1_2} \textbf{ Zero crossing function $\omega_z(t_0)$ is an \textbf{even} function of variable $t_0$ } \protect\\  \lowercase{} }

Using $g_2(t, t_0)= \frac{f_2(t, t_0)}{h(t)}$ in Result 2.3.a, we get $g_2(t, t_0) = f_2(t, t_0) e^{-\sigma t}  u(-t) +  f_2(t, t_0) e^{\sigma t}  u(t)$ (\textbf{Result 2.4.a}), where $f_2(t, t_0)= f_1(t,-t_0) = E_0(t - t_0 )  e^{-\sigma t} $, using Eq.~\ref{sec:sec_2_1_eq_0} and Result 2.1.c in Section~\ref{sec:Section_2_1}. We compute the Fourier transform of the function $g_2(t, t_0)$ given by $G_2(\omega, t_0)=G_1(\omega, -t_0)$, given that $f_2(t, t_0) = f_1(t, -t_0)$ using Result 2.1.c in Section~\ref{sec:Section_2_1}  and hence $g_2(t, t_0) = g_1(t, -t_0)$ using Result 2.4.a and hence $G_{2}(\omega, t_0) = G_{1}(\omega, -t_0)$.  \\

We compute the Fourier transform of the function $g(t, t_0)=  e^{- \sigma t_0} g_1(t, t_0) +  e^{\sigma t_0}  g_2(t, t_0) + 2 \cosh{( \sigma t_0)} g_3(t)$ in Result 2.3.a, given by $G(\omega, t_0)$, and compute its real part $G_{R}(\omega, t_0)$ using the procedure in Section~\ref{sec:Section_2_1_a}, similar to Eq.~\ref{sec:sec_2_0_eq_5} and Eq.~\ref{sec:sec_2_1_eq_5} and we can write as follows in Eq.~\ref{sec:sec_2_1_eq_9a}. We substitute $t=\tau$ in the equation for $G_{3_R}(\omega)$ and $G_{1_R}(\omega, t_0) $ below, copied from Eq.~\ref{sec:sec_2_0_eq_5} and Eq.~\ref{sec:sec_2_1_eq_5}.
\begin{align}\label{sec:sec_2_1_eq_9a}  
G_{R}(\omega, t_0) =  e^{- \sigma t_0} G_{1_R}(\omega, t_0) + e^{\sigma t_0}  G_{2_R}(\omega, t_0) + 2 \cosh{( \sigma t_0)} G_{3_R}(\omega) \notag\\ =  e^{- \sigma t_0} G_{1_R}(\omega, t_0) + e^{\sigma t_0}  G_{1_R}(\omega, -t_0) + 2 \cosh{( \sigma t_0)} G_{3_R}(\omega) \notag\\
G_{3_R}(\omega) =   \int_{-\infty}^{0}  E_0(\tau) ( e^{-2 \sigma \tau} + 1) \cos{(\omega \tau)} d\tau \notag\\
G_{1_R}(\omega, t_0) = [\int_{-\infty}^{0}   E_0(\tau + t_0  ) e^{-2 \sigma \tau}  \cos{(\omega \tau)} d\tau 
+ \int_{-\infty}^{0}    E_0(\tau - t_0 )    \cos{(\omega \tau)} d\tau ] \notag\\
G_{R}(\omega, t_0) = e^{- \sigma t_0}  [\int_{-\infty}^{0}   E_0(\tau + t_0  ) e^{-2 \sigma \tau}  \cos{(\omega \tau)} d\tau 
+ \int_{-\infty}^{0}   E_0(\tau - t_0 )   \cos{(\omega \tau)} d\tau ] \notag\\
+ e^{\sigma t_0}  [\int_{-\infty}^{0}   E_0(\tau - t_0  )  e^{-2 \sigma \tau}  \cos{(\omega \tau)} d\tau 
+ \int_{-\infty}^{0}    E_0(\tau + t_0 )  \cos{(\omega \tau)} d\tau ] \notag\\ + 2 \cosh{( \sigma t_0)} \int_{-\infty}^{0}  E_0(\tau) ( e^{-2 \sigma \tau} + 1) \cos{(\omega \tau )} d\tau 
\end{align}
We combine first and fourth integrals and then combine third and second integrals in $G_{R}(\omega, t_0) $ in Eq.~\ref{sec:sec_2_1_eq_9a}, as follows. 
\begin{align}\label{sec:sec_2_1_eq_9}  
G_{R}(\omega, t_0) =  \int_{-\infty}^{0}  E_0(\tau + t_0 )  ( e^{- \sigma t_0} e^{-2 \sigma \tau} + e^{\sigma t_0})]  \cos{(\omega \tau)} d\tau \notag\\
+ \int_{-\infty}^{0}    E_0(\tau - t_0  )  ( e^{\sigma t_0} e^{-2 \sigma \tau} + e^{- \sigma t_0})   \cos{(\omega \tau)} d\tau 
+ 2 \cosh{( \sigma t_0)} \int_{-\infty}^{0}  E_0(\tau) ( e^{-2 \sigma \tau} + 1) \cos{(\omega \tau )} d\tau 
\end{align}
We require $G_{R}(\omega, t_0) =0$ for $\omega=\omega_{z}(t_0)$ for every value of $t_0$, to satisfy \textbf{Statement 1}, using Lemma 1 in Section~\ref{sec:Section_2_1}. Hence we write $P(t_0) =  G_{R}(\omega_{z}(t_0), t_0) =  0$ as follows. 
\begin{align}\label{sec:sec_2_1_eq_10_a} 
P(t_0) =  \int_{-\infty}^{0}  E_0(\tau + t_0 )  ( e^{- \sigma t_0} e^{-2 \sigma \tau} + e^{\sigma t_0})]  \cos{(\omega_{z}(t_0) \tau)} d\tau \notag\\
+ \int_{-\infty}^{0}    E_0(\tau - t_0  )  ( e^{\sigma t_0} e^{-2 \sigma \tau} + e^{- \sigma t_0})   \cos{(\omega_{z}(t_0) \tau)} d\tau ] 
+ 2 \cosh{( \sigma t_0)} \int_{-\infty}^{0}  E_0(\tau) ( e^{-2 \sigma \tau} + 1) \cos{(\omega_{z}(t_0) \tau )} d\tau    = 0 
\end{align}
We use the fact that $f(t, t_0)  =  e^{- \sigma t_0}  E_0(t + t_0 )  e^{-\sigma t}  + e^{\sigma t_0}  E_0(t - t_0 )  e^{-\sigma t} + 2 \cosh{( \sigma t_0)} E_p(t) = f(t, -t_0)$ in Eq.~\ref{sec:sec_2_1_eq_0}, is \textbf{unchanged} by the substitution $t_0=-t_0$. \textbf{If} $f(t, t_0)=f(t, -t_0)$ is unchanged by the substitution $t_0=-t_0$, \textbf{then} $g(t, t_0)=g(t, -t_0)$ is unchanged by the substitution $t_0=-t_0$, using the fact that $g(t, t_0)h(t)= f(t, t_0)$ and $h(t)= [ e^{ \sigma t} u(-t) + e^{ - \sigma t} u(t) ] $ in Eq.~\ref{sec:sec_2_1_eq_0}. \\

Hence the Fourier transform of  $g(t, t_0)$ given by $G(\omega, t_0)=G(\omega,  -t_0)$ and its real part given by $G_{R}(\omega, t_0)=G_{R}(\omega, -t_0)$ is \textbf{unchanged} by the substitution $t_0=-t_0$ and the zero crossing in $G_{R}(\omega,  -t_0)$ given by $\omega_{z}( -t_0)$ is the \textbf{same} as the zero crossing in $G_{R}(\omega, t_0)$ given by $\omega_{z}(t_0)$  and we get $\omega_{z}(t_0)=\omega_{z}( -t_0)$ and hence $\omega_{z}(t_0)$ is an \textbf{even} function of variable $t_0$. (\textbf{Result 2.4.b})\\

We can write Eq.~\ref{sec:sec_2_1_eq_10_a} as follows, where $P_{odd}(t_0)$ is an \textbf{odd} function of variable $t_0$.  We use $\omega_{z}(t_0)=\omega_{z}( -t_0)$ and $D(t_0)= D(-t_0)= \cosh{( \sigma t_0)} \int_{-\infty}^{0}  E_0(\tau) ( e^{-2 \sigma \tau} + 1) \cos{(\omega_{z}(t_0) \tau )} d\tau$.(\textbf{Result 2.4.c})
\begin{align}\label{sec:sec_2_1_eq_10}   
P(t_0)= P_{odd}(t_0) +  P_{odd}( -t_0) = 0 \notag\\
P_{odd}(t_0) =  \int_{-\infty}^{0}  E_0(\tau + t_0 )  ( e^{- \sigma t_0} e^{-2 \sigma \tau} + e^{\sigma t_0})]  \cos{(\omega_{z}(t_0) \tau)} d\tau +  D(t_0)
\end{align}

.\\ \tocless\section{\label{sec:Section_A_1_3} \textbf{ Final Step} \protect\\  \lowercase{} }

We expand $ P_{odd}(t_0)$ in  Eq.~\ref{sec:sec_2_1_eq_10}  as follows, using the substitution $\tau + t_0  = \tau'$. We get $\tau = \tau' - t_0$ and $d\tau = d\tau'$ and substitute back $\tau' = \tau$ in the second line below. We use $e^{- \sigma t_0} e^{2 \sigma t_0} = e^{\sigma t_0}$ and  use the identity $   \cos{(\omega_z(t_0)  (\tau -  t_0)} =\cos{(\omega_z(t_0)  t_0)} \cos{(\omega_z(t_0) \tau)} + \sin{(\omega_z(t_0)  t_0)} \sin{(\omega_z(t_0) \tau)}$ below.
\begin{align}\label{sec:sec_az3_1_3_eq_1}   
P_{odd}(t_0) =  \int_{-\infty}^{t_0}    E_0(\tau' )  [ e^{- \sigma t_0} e^{-2 \sigma \tau'} e^{2 \sigma t_0}  + e^{\sigma t_0}   ]  \cos{(\omega_{z}(t_0) (\tau^{'} - t_0)} d\tau' +  D(t_0) \notag\\
P_{odd}(t_0) =   [ \cos{(\omega_z(t_0) t_0)}  \int_{-\infty}^{t_0}    E_0(\tau )  ( e^{\sigma t_0}  e^{-2 \sigma \tau} + e^{\sigma t_0}) \cos{(\omega_z(t_0) \tau)} d\tau  \notag\\  + \sin{(\omega_z(t_0) t_0)}  \int_{-\infty}^{t_0}   E_0(\tau ) ( e^{\sigma t_0} e^{-2 \sigma \tau} + e^{\sigma t_0}) \sin{(\omega_z(t_0) \tau)} d\tau ]  +  D(t_0)
\end{align}
In Lemma 1 in Section~\ref{sec:Section_2_1}, it is shown that $0 < \omega_z(t_0) < \infty$, for $t_0 \in \Re$. \\ 

In Section~\ref{sec:Section_A_1_8_1}, it is shown that $\omega_z(t_0)$ is a \textbf{continuous}  function of variable $t_0$, for $t_0 \in \Re$.\\

In Section~\ref{sec:Section_A_1_6}, it is shown that $E_0(t)$ is \textbf{strictly decreasing} for $t > 0$.\\

Given that $\omega_z(t_0)$ is a continuous function of $t_0$ and given that $t_0$ is a continuous function, we see that the \textbf{product} of two continuous functions $\omega_z(t_0) t_0$ is a \textbf{continuous }function and is positive for $t_0>0$ because $0 < \omega_z(t_0) < \infty$. \\

It is shown in Section~\ref{sec:Section_C_6_0} that $\omega_z(t_0) t_0$ starts from zero and increases to a larger and larger positive value, as $t_0$ increases to a larger and larger finite value without bounds and that $\omega_z(t_0) t_0$ does not go to zero and does not approach a constant, as $t_0$ increases.
Hence $\omega_z(t_0) t_0 =   \pi$, can be reached for specific values of $t_0$, as finite $t_0$ increases without bounds. As $t_0$ increases from $0$, to a larger and larger finite value without bounds, the continuous function $\omega_z(t_0) t_0$ starts from $0$ and  will pass through $\pi$, using Intermediate Value Theorem,  for specific value of $t_0$.\\ 

We set $t_0=t_{0c} > 0$ and  such that $\omega_z(t_{0c}) t_{0c}=   \pi$, in Eq.~\ref{sec:sec_az3_1_3_eq_1} as follows. We use the fact that $\cos{( \omega_z(t_{0c}) t_{0c} )} = -1$, $\sin{( \omega_z(t_{0c}) t_{0c} )} = 0$ and  $\omega_z(-t_{0c}) = \omega_z(t_{0c})$ shown in Result 2.4.b in Section~\ref{sec:Section_A_1_2}. 
\begin{align}\label{sec:sec_az3_1_3_0_eq_2a}   
P_{odd}(t_{0c})=  - \int_{-\infty}^{t_{0c}}    E_0(\tau)  ( e^{  \sigma  t_{0c}} e^{-  2 \sigma  \tau} + e^{  \sigma  t_{0c}} ) \cos{( \omega_z(t_{0c}) \tau)} d\tau  + D(t_{0c})
\end{align}
We compute $P_{odd}(-t_0)$ in Eq.~\ref{sec:sec_az3_1_3_eq_1} as follows. We use $\omega_z( -t_{0}) = \omega_z( t_{0})$ in Result 2.4.b and $D(t_0)=D(-t_0)$ in Result 2.4.c in Section~\ref{sec:Section_A_1_2}.
\begin{align}\label{sec:sec_az3_1_3_eq_2b}   
P_{odd}(-t_0) =   [ \cos{(\omega_z(t_0) t_0)}  \int_{-\infty}^{-t_0}    E_0(\tau )  ( e^{-\sigma t_0}  e^{-2 \sigma \tau} + e^{-\sigma t_0}) \cos{(\omega_z(t_0) \tau)} d\tau  \notag\\  - \sin{(\omega_z(t_0) t_0)}  \int_{-\infty}^{-t_0}   E_0(\tau ) ( e^{-\sigma t_0} e^{-2 \sigma \tau} + e^{-\sigma t_0}) \sin{(\omega_z(t_0) \tau)} d\tau ]  +  D(t_0)
\end{align}
We set $t_0=t_{0c} > 0$  such that $\omega_z(t_{0c}) t_{0c}=  \pi $, in Eq.~\ref{sec:sec_az3_1_3_eq_2b} as follows. We use $\cos{( \omega_z(t_{0c}) t_{0c} )} = -1$, $\sin{( \omega_z(t_{0c}) t_{0c} )} = 0$. 
\begin{align}\label{sec:sec_az3_1_3_0_eq_2c}   
P_{odd}( -t_{0c})= - \int_{-\infty}^{-t_{0c}}    E_0(\tau)  ( e^{ - \sigma  t_{0c}} e^{-  2 \sigma  \tau} +  e^{  -\sigma  t_{0c}} ) \cos{( \omega_z(t_{0c}) \tau)} d\tau  + D(t_{0c})
\end{align}
We compute $ [ P_{odd}(t_0) + P_{odd}( -t_0) ] = 0$ in Eq.~\ref{sec:sec_2_1_eq_10}, at $t_0=t_{0c}$ using Eq.~\ref{sec:sec_az3_1_3_0_eq_2a} and Eq.~\ref{sec:sec_az3_1_3_0_eq_2c}. 
\begin{align}\label{sec:sec_az3_1_3_0_eq_2} 
- \int_{-\infty}^{t_{0c}}    E_0(\tau)  ( e^{  \sigma  t_{0c}} e^{-  2 \sigma  \tau} + e^{  \sigma  t_{0c}} ) \cos{( \omega_z(t_{0c}) \tau)} d\tau  \notag\\
 - \int_{-\infty}^{-t_{0c}}    E_0(\tau)  ( e^{ - \sigma  t_{0c}} e^{-  2 \sigma  \tau} +  e^{  -\sigma  t_{0c}} ) \cos{( \omega_z(t_{0c}) \tau)} d\tau  + 2 D(t_{0c})  = 0
\end{align}
We split the integrals in the left hand side of  Eq.~\ref{sec:sec_az3_1_3_0_eq_2} using $ \int_{-\infty}^{t_{0c}} = \int_{-\infty}^{0}  +  \int_{0}^{t_{0c}}$ as follows.
\begin{align}\label{sec:sec_ab_1_3_0_eq_2_b}   
- \int_{-\infty}^{0}    E_0(\tau)  ( e^{  \sigma  t_{0c}} e^{-  2 \sigma  \tau} + e^{  \sigma  t_{0c}} ) \cos{( \omega_z(t_{0c}) \tau)} d\tau  
- \int_{0}^{t_{0c}}    E_0(\tau)  ( e^{  \sigma  t_{0c}} e^{-  2 \sigma  \tau} + e^{  \sigma  t_{0c}} ) \cos{( \omega_z(t_{0c}) \tau)} d\tau  \notag\\
 - \int_{-\infty}^{0}    E_0(\tau)  ( e^{ - \sigma  t_{0c}} e^{-  2 \sigma  \tau} +  e^{  -\sigma  t_{0c}} ) \cos{( \omega_z(t_{0c}) \tau)} d\tau 
 - \int_{0}^{-t_{0c}}    E_0(\tau)  ( e^{ - \sigma  t_{0c}} e^{-  2 \sigma  \tau} +  e^{  -\sigma  t_{0c}} ) \cos{( \omega_z(t_{0c}) \tau)} d\tau  \notag\\ + 2 D(t_{0c})  = 0
\end{align}
We keep the second and fourth integrals in the left hand side and consolidate the rest of the integrals in Eq.~\ref{sec:sec_ab_1_3_0_eq_2_b} to the right hand side as follows, using $2 \cosh{( \sigma t_{0c})} = e^{  \sigma  t_{0c}} + e^{ - \sigma  t_{0c}}$.
\begin{align}\label{sec:sec_az3_1_3_0_0_eq_3_b} 
- \int_{0}^{t_{0c}}    E_0(\tau)  ( e^{  \sigma  t_{0c}} e^{-  2 \sigma  \tau} + e^{  \sigma  t_{0c}} ) \cos{( \omega_z(t_{0c}) \tau)} d\tau  
 - \int_{0}^{-t_{0c}}    E_0(\tau)  ( e^{ - \sigma  t_{0c}} e^{-  2 \sigma  \tau} +  e^{  -\sigma  t_{0c}} ) \cos{( \omega_z(t_{0c}) \tau)} d\tau \notag\\
 =  2 \cosh{( \sigma t_{0c})}  \int_{-\infty}^{0}     E_0(\tau )  e^{-  2 \sigma  \tau} \cos{( \omega_z(t_{0c}) \tau)} d\tau 
+ 2 \cosh{( \sigma t_{0c})}  \int_{-\infty}^{0}     E_0(\tau )   \cos{( \omega_z(t_{0c}) \tau)} d\tau - 2 D(t_{0c})
\end{align}
We denote the right hand side of Eq.~\ref{sec:sec_az3_1_3_0_0_eq_3_b} as $RHS$. We can substitute $\tau=-\tau$ in the second integral in the left hand side (LHS) of  Eq.~\ref{sec:sec_az3_1_3_0_0_eq_3_b}  as follows.  We use $E_0(-\tau) = E_0(\tau)$. 
\begin{align}\label{sec:sec_az3_1_3_0_0_eq_3_d}   
LHS = - \int_{0}^{t_{0c}}    E_0(\tau)  ( e^{  \sigma  t_{0c}} e^{-  2 \sigma  \tau} + e^{  \sigma  t_{0c}} ) \cos{( \omega_z(t_{0c}) \tau)} d\tau  
 + \int_{0}^{t_{0c}}    E_0(\tau)  ( e^{ - \sigma  t_{0c}} e^{  2 \sigma  \tau} +  e^{  -\sigma  t_{0c}} ) \cos{( \omega_z(t_{0c}) \tau)} d\tau 
  \end{align}
We combine the integrals in Eq.~\ref{sec:sec_az3_1_3_0_0_eq_3_d} as follows.
\begin{align}\label{sec:sec_az3_1_3_0_0_eq_3_d1}   
LHS = \int_{0}^{t_{0c}}    E_0(\tau)  ( -e^{  \sigma  t_{0c}} e^{-  2 \sigma  \tau} - e^{  \sigma  t_{0c}} + e^{ - \sigma  t_{0c}} e^{  2 \sigma  \tau} +  e^{  -\sigma  t_{0c}} ) \cos{( \omega_z(t_{0c}) \tau)} d\tau 
\end{align}
We substitute $ 2 \sinh{( \sigma  \tau)} =  e^{  \sigma  \tau}  -  e^{ - \sigma  \tau} $ and $ 2 \sinh{( \sigma  t_{0c} )} =  e^{  \sigma  t_{0c}}  -  e^{- \sigma  t_{0c}} $ in Eq.~\ref{sec:sec_az3_1_3_0_0_eq_3_d1} as follows.
  \begin{align}\label{sec:sec_az3_1_3_0_0_eq_3_d2}   
LHS = - 2 \int_{0}^{t_{0c}}     E_0(\tau ) ( \sinh{( \sigma  t_{0c} - 2 \sigma  \tau)} + \sinh{(\sigma  t_{0c} )}) \cos{(\omega_z(t_{0c}) \tau)} d\tau 
  \end{align}
We substitute $\tau=-\tau$ in the integrals in the right hand side of Eq.~\ref{sec:sec_az3_1_3_0_0_eq_3_b} as follows.  We use $E_0(\tau) =E_0(-\tau)$ and $-\int_{\infty}^{0}= \int_{0}^{\infty}$.
\begin{align}\label{sec:sec_az3_1_3_0_0_eq_3_g}   
RHS  =  -2 \cosh{( \sigma t_{0c})}  \int_{\infty}^{0}     E_0(\tau )  e^{  2 \sigma  \tau} \cos{( \omega_z(t_{0c}) \tau)} d\tau 
- 2 \cosh{( \sigma t_{0c})}  \int_{\infty}^{0}     E_0(\tau )   \cos{( \omega_z(t_{0c}) \tau)} d\tau - 2 D(t_{0c}) \notag\\
RHS  =  2 \cosh{( \sigma t_{0c})}  \int_{0}^{\infty}     E_0(\tau )  e^{  2 \sigma  \tau} \cos{( \omega_z(t_{0c}) \tau)} d\tau 
+ 2 \cosh{( \sigma t_{0c})}  \int_{0}^{\infty}     E_0(\tau )   \cos{( \omega_z(t_{0c}) \tau)} d\tau - 2 D(t_{0c}) \notag\\
RHS  =  2 \cosh{( \sigma t_{0c})}  \int_{0}^{\infty}     E_0(\tau ) [ e^{  2 \sigma  \tau} + 1 ]\cos{( \omega_z(t_{0c}) \tau)} d\tau - 2 D(t_{0c}) 
\end{align}
We equate Eq.~\ref{sec:sec_az3_1_3_0_0_eq_3_d2}  and Eq.~\ref{sec:sec_az3_1_3_0_0_eq_3_g} after dividing by a factor of 2,  as follows.  
\begin{align}\label{sec:sec_az3_1_3_0_0_eq_3_g1_0a}   
-  \int_{0}^{t_{0c}}     E_0(\tau ) ( \sinh{( \sigma  t_{0c} - 2 \sigma  \tau)} + \sinh{(\sigma  t_{0c} )}) \cos{(\omega_z(t_{0c}) \tau)} d\tau \notag\\
=   \cosh{( \sigma t_{0c})}  \int_{0}^{\infty}     E_0(\tau ) [ e^{  2 \sigma  \tau} + 1 ] \cos{( \omega_z(t_{0c}) \tau)} d\tau  -  D(t_{0c})
\end{align}
We use $D(t_0)=  \cosh{( \sigma t_0)} \int_{-\infty}^{0}  E_0(\tau) ( e^{-2 \sigma \tau} + 1) \cos{(\omega_{z}(t_0) \tau )} d\tau$ using Result 2.4.c in Section 2.4 and we substitute $\tau=-\tau$ in the integrand and $t_0=t_{0c}$ and we get $D(t_{0c})=  \cosh{( \sigma t_{0c})} \int_{0}^{\infty}  E_0(\tau) ( e^{2 \sigma \tau} + 1) \cos{(\omega_{z}(t_{0c}) \tau )} d\tau$ using $E_0(\tau) =E_0(-\tau)$ and $-\int_{\infty}^{0}= \int_{0}^{\infty}$. Hence we can cancel RHS terms in Eq.~\ref{sec:sec_az3_1_3_0_0_eq_3_g1_0a} as follows. We multiply the equation by a factor of $-1$ in the second line below.
\begin{align}\label{sec:sec_az3_1_3_0_0_eq_3_g1_0}   
-A(t_{0c}) = -  \int_{0}^{t_{0c}}     E_0(\tau ) ( \sinh{( \sigma  t_{0c} - 2 \sigma  \tau)} + \sinh{(\sigma  t_{0c} )}) \cos{(\omega_z(t_{0c}) \tau)} d\tau  = 0 \notag\\
A(t_{0c}) = \int_{0}^{t_{0c}}     E_0(\tau ) (  \sinh{(\sigma  t_{0c} )} - \sinh{( 2 \sigma  \tau - \sigma  t_{0c} )} )  \cos{(\omega_z(t_{0c}) \tau)} d\tau  = 0 
\end{align}

\textbf{Eq.~\ref{sec:sec_az3_1_3_0_0_eq_3_g1_0} leads to contradiction:}\\

We will show that $A(t_{0c}) > 0$ and Eq.~\ref{sec:sec_az3_1_3_0_0_eq_3_g1_0} leads to a contradiction for $0 < \sigma < \frac{1}{2}$.\\

We see that $E_0(\tau) > 0$ for real $\tau$ (~\ref{sec:appendix_C_5a}). We see that $E_0(\tau)$ is a \textbf{strictly decreasing} function for $\tau > 0$ as shown in Section~\ref{sec:Section_A_1_6}. (\textbf{Result 3.0.a}) \\

The first derivative of $ B(\tau, t_{0c})=\sinh{(\sigma  t_{0c} )} - \sinh{( 2 \sigma  \tau - \sigma  t_{0c} )}$ in Eq.~\ref{sec:sec_az3_1_3_0_0_eq_3_g1_0}, given by $\frac{\partial B(\tau, t_{0c})}{\partial \tau} = - 2 \sigma \cosh{( 2 \sigma  \tau - \sigma  t_{0c} )} < 0$, for $0 < \sigma < \frac{1}{2}$ and hence $ B(\tau, t_{0c})=\sinh{(\sigma  t_{0c} )} - \sinh{( 2 \sigma  \tau - \sigma  t_{0c} )}$ is a \textbf{strictly decreasing function} of $\tau$ in the interval $0 < \tau < t_{0c}$. (\textbf{Result 3.0.b}) \\

We see that $ B(\tau, t_{0c}) = 2 \sinh{(\sigma  t_{0c} )} > 0$ at $\tau=0$ and $B(\tau, t_{0c})= 0$ at $\tau=t_{0c}$. Given that $B(\tau, t_{0c})$ is a \textbf{strictly decreasing function} of $\tau$ in the interval $0 < \tau < t_{0c}$, we see that $B(\tau, t_{0c}) >  0$ in the interval $0 < \tau < t_{0c}$. (\textbf{Result 3.0.c}) \\

Using Result 3.0.a, Result 3.0.b and Result 3.0.c, we see that the term $C(\tau, t_{0c}) = E_0(\tau ) B(\tau, t_{0c})$ in Eq.~\ref{sec:sec_az3_1_3_0_0_eq_3_g1_0}, is a \textbf{strictly decreasing function} of $\tau$, given that the product of two strictly decreasing functions, which are positive in the interval $0 < \tau < t_{0c}$, is also a strictly decreasing function, because $C(\tau, t_{0c}) - C(\tau + d\tau, t_{0c}) = E_0(\tau ) B(\tau, t_{0c}) - E_0(\tau + d\tau ) B(\tau + d\tau, t_{0c}) > 0$, given that $E_0(\tau )  > E_0(\tau + d\tau )$ and $B(\tau, t_{0c}) > B(\tau + d\tau, t_{0c})$, for $d\tau > 0$, in the interval $0 < \tau < t_{0c}$ and $0 < \tau + d\tau < t_{0c}$.  (\textbf{Result 3.0.d})\\

Using Result 3.0.a, Result 3.0.b, Result 3.0.c and Result 3.0.d and $C(\tau, t_{0c}) = E_0(\tau ) (\sinh{(\sigma  t_{0c} )} - \sinh{( 2 \sigma  \tau - \sigma  t_{0c} )})$, we can split Eq.~\ref{sec:sec_az3_1_3_0_0_eq_3_g1_0}  to 2 integrals, as follows. 
\begin{align}\label{sec:sec_3_3_eq_1}   
A(t_{0c}) = \int_{0}^{t_{0c}}  C(\tau, t_{0c})  \cos{(\omega_z(t_{0c}) \tau)} d\tau = \int_{0}^{\frac{t_{0c}}{2}}  C(\tau, t_{0c})  \cos{(\omega_z(t_{0c}) \tau)} d\tau +  \int_{\frac{t_{0c}}{2}}^{t_{0c}}  C(\tau, t_{0c})  \cos{(\omega_z(t_{0c}) \tau)} d\tau
\end{align}
We substitute $t_{0c} - \tau = \tau^{'}$ and hence $\tau = t_{0c} - \tau^{'}$ and $d\tau = - d\tau^{'}$ in the second integral in the right hand side of Eq.~\ref{sec:sec_3_3_eq_1} and substitute back $\tau^{'} = \tau$ in the second line below and write as follows.
\begin{align}\label{sec:sec_3_3_eq_2}   
A(t_{0c}) =  \int_{0}^{\frac{t_{0c}}{2}} C(\tau, t_{0c})  \cos{(\omega_z(t_{0c}) \tau)}  d\tau -  \int_{\frac{t_{0c}}{2}}^{0}  C( t_{0c} - \tau^{'}, t_{0c})  \cos{(\omega_z(t_{0c}) (t_{0c} - \tau^{'}) )} d\tau^{'} \notag\\
A(t_{0c}) =  \int_{0}^{\frac{t_{0c}}{2}} C(\tau, t_{0c})  \cos{(\omega_z(t_{0c}) \tau)}  d\tau -  \int_{\frac{t_{0c}}{2}}^{0}  C( t_{0c} - \tau, t_{0c})  \cos{(\omega_z(t_{0c}) (t_{0c} - \tau) )} d\tau
\end{align}
We expand above equation using $ \cos{(\omega_z(t_{0c}) (t_{0c} - \tau) )} =  \cos{(\omega_z(t_{0c}) t_{0c} )}  \cos{(\omega_z(t_{0c}) \tau )} \\ +  \sin{(\omega_z(t_{0c}) t_{0c} )}  \sin{(\omega_z(t_{0c})  \tau )}$ and use $\omega_z(t_{0c}) t_{0c} = \pi$ and hence $\sin{(\omega_z(t_{0c}) t_{0c} )}=0$ and \\
$\cos{(\omega_z(t_{0c}) t_{0c} )}  = -1$. We use $-  \int_{\frac{t_{0c}}{2}}^{0} = \int_{0}^{\frac{t_{0c}}{2}}$ below.
\begin{align}\label{sec:sec_3_3_eq_3}   
A(t_{0c}) =  \int_{0}^{\frac{t_{0c}}{2}} C(\tau, t_{0c})  \cos{(\omega_z(t_{0c}) \tau)}  d\tau - \int_{0}^{\frac{t_{0c}}{2}}  C(t_{0c} - \tau, t_{0c})  \cos{(\omega_z(t_{0c}) \tau )} d\tau \notag\\
A(t_{0c}) =  \int_{0}^{\frac{t_{0c}}{2}} [ C(\tau, t_{0c}) - C(t_{0c} - \tau, t_{0c})  ] \cos{(\omega_z(t_{0c}) \tau)}  d\tau 
\end{align}
At $\tau=0$, the term $C(\tau, t_{0c}) = E_0(\tau ) B(\tau, t_{0c})  > 0$, given that $E_0(0) > 0$ and $ B(\tau, t_{0c}) = \sinh{(\sigma  t_{0c} )} - \sinh{( 2 \sigma  \tau - \sigma  t_{0c} )}=2 \sinh{(\sigma  t_{0c} )} > 0$ at $\tau=0$. We see that $C(t_{0c} - \tau, t_{0c})=0$ at  $\tau=0$, given that $C(t_{0c} - \tau, t_{0c})= E_0(t_{0c} - \tau ) B(t_{0c} - \tau, t_{0c})$ and $ B( t_{0c} -\tau, t_{0c}) = \sinh{(\sigma  t_{0c} )} - \sinh{( 2 \sigma  (t_{0c} -\tau) - \sigma  t_{0c} )}=0$ at $\tau=0$. Hence $ I(\tau, t_{0c})= C(\tau, t_{0c}) - C(t_{0c} - \tau , t_{0c}) > 0$ at $\tau=0$ in Eq.~\ref{sec:sec_3_3_eq_3}. At $\tau=\frac{t_{0c}}{2}$, we see that $ I(\frac{t_{0c}}{2}, t_{0c}) = C(\frac{t_{0c}}{2}, t_{0c}) - C(\frac{t_{0c}}{2} , t_{0c}) = 0 $.\\

Using \textbf{strictly decreasing} property of $C(\tau, t_{0c})$ in Result 3.0.d, we see that $C(\tau, t_{0c}) - C(t_{0c} - \tau, t_{0c})  > 0$ in the interval $0 < \tau < \frac{t_{0c}}{2} $. The term $\cos{(\omega_z(t_{0c}) \tau )} > 0$ in the interval $0 < \tau < \frac{t_{0c}}{2} $ using $\omega_z(t_{0c}) \frac{t_{0c}}{2} = \frac{\pi}{2} $. Hence the integral for $A(t_{0c})$ in Eq.~\ref{sec:sec_3_3_eq_3}  and in LHS of Eq.~\ref{sec:sec_az3_1_3_0_0_eq_3_g1_0} is $>0$.(\textbf{Result 3.0.e}) \\

We require $A(t_{0c})= 0$ in Eq.~\ref{sec:sec_az3_1_3_0_0_eq_3_g1_0}. Hence this leads to a \textbf{contradiction}, for $0 < \sigma < \frac{1}{2}$.\\

For $\sigma=0$, both sides of Eq.~\ref{sec:sec_az3_1_3_0_0_eq_3_g1_0}  is zero, given the term $ \sinh{(\sigma  t_{0c} )} - \sinh{( 2 \sigma  \tau - \sigma  t_{0c} )}= 0$ and \textbf{does not} lead to a contradiction. \\

We have shown this result for $0 < \sigma < \frac{1}{2}$. \textbf{If} the Fourier transform of $E_p(t)= E_0(t) e^{-\sigma t}$ given by $E_{p\omega}(\omega)=E_{pR\omega}(\omega) + i E_{pI\omega}(\omega) $ has a zero at $\omega=\omega_0$, to satisfy Statement 1, \textbf{then} the real part $E_{pR\omega}(\omega)$ and the imaginary part $E_{pI\omega}(\omega)$ \textbf{also} have a zero at $\omega=\omega_0$.(\textbf{Result 3.b}) \\

Given that $E_p(t) =  E_0(t) e^{-\sigma t}$  is real, its Fourier transform $E_{p\omega}(\omega) = \xi(\frac{1}{2} + \sigma + i \omega) $ has symmetry properties and hence $E_{pR\omega}(-\omega)=E_{pR\omega}(\omega)$  and $E_{pI\omega}(-\omega)= - E_{pI\omega}(\omega)$ (\href{https://archive.is/CLVqX}{Symmetry property}) and hence $E_{p\omega}(-\omega)=\xi(\frac{1}{2} + \sigma - i \omega)=E_{pR\omega}(-\omega) + i E_{pI\omega}(-\omega)=E_{pR\omega}(\omega) - i E_{pI\omega}(\omega)$ \textbf{also} has a zero at $\omega=\omega_0$ to satisfy Statement 1, using Result 3.b. (\textbf{Result 3.c})\\

Using the property $\xi(s)=\xi(1-s)$, we get  $\xi(\frac{1}{2} + \sigma - i \omega) = \xi(\frac{1}{2} - \sigma + i \omega)$ at $s=\frac{1}{2} + \sigma - i \omega$ and hence $E_{q\omega}(\omega)=\xi(\frac{1}{2} - \sigma + i \omega)$ \textbf{also} has a zero at $\omega=\omega_0$  to satisfy Statement 1, using Result 3.c. We see that $E_{q\omega}(\omega)=\xi(\frac{1}{2} - \sigma + i \omega)$ is obtained by replacing $\sigma$ in $E_{p\omega}(\omega)= \xi(\frac{1}{2} + \sigma + i \omega) $ by $-\sigma$. Hence the results in above sections hold for $-\frac{1}{2} < \sigma < 0$ and hence for  $0 < |\sigma| < \frac{1}{2}$.\\

Hence we have produced a \textbf{contradiction} of  \textbf{Statement 1} that the Fourier Transform of the function $E_p(t) = E_0(t) e^{-\sigma t} $ has a zero at $\omega = \omega_{0}$ for  $0 < |\sigma| < \frac{1}{2}$.\\

Hence the assumption in \textbf{Statement 1} that  Riemann's Xi Function given by $\xi(\frac{1}{2} + \sigma + i \omega)= E_{p\omega}(\omega)$ has a zero at $\omega= \omega_{0}$, where $\omega_{0}$ is real, leads to a \textbf{contradiction} for the region $0 < |\sigma| < \frac{1}{2}$ which corresponds to the critical strip excluding the critical line. Hence $\zeta(s)$ does not have non-trivial zeros in the critical strip excluding the critical line and we have proved Riemann's Hypothesis.

.\\ \tocless\section{\label{sec:Section_A_1_8_1} \textbf{ $\omega_z(t_0)$ is a continuous  function of $t_0$ } \protect\\  \lowercase{} }

It is shown in \textbf{Lemma 1} in Section~\ref{sec:Section_2_1} that  $G_{R}(\omega, t_0)=0$ at $\omega=\omega_z(t_0)$ where it crosses the zero line to the opposite sign, if Statement 1 is true, and that $\omega_z(t_0)$ is \textbf{finite and non-zero} for all $t_0 \in \Re$ and that $\omega_z(t_0)$ is an even function of variable $t_0$ (Result 2.4.b in Section ~\ref{sec:Section_A_1_2}). For a given $t_0$, $\omega_z(t_0)$ can have more than one value, corresponding to multiple zero crossings in $G_{R}(\omega, t_0)$, but we consider only the first zero crossing to the right of origin in the section below, where $G_{R}(\omega, t_0)$  crosses the zero line to the opposite sign, as detailed in \textbf{Lemma 1} in Section~\ref{sec:Section_2_1}. \\

We consider the Fourier transform of the even part of $g(t, t_0)$ given by $G_{R}(\omega, t_0) $ in the section below and show that, under this Fourier transformation, as we change $t_0$, the zero crossing in $G_{R}(\omega, t_0) $  given by $\omega_z(t_0)$ is a continuous  function of $t_0$, for $t_0 \in \Re$ . This is shown in the steps below using \textbf{Implicit Function Theorem}. \\

$\bullet$  It is shown in Section~\ref{sec:Section_A_1_8_1_2_0} that $G_{R}(\omega, t_0)$ and $G_{R, 2r}(\omega, t_0)= \frac{\partial^{2r} G_{R}(\omega, t_0)}{\partial \omega^{2r}}$ are partially differentiable at least twice with respect to $\omega$, for some value of $r \in W$ (element of set of whole numbers including zero.)\\

$\bullet$ It is shown in Section~\ref{sec:Section_A_1_8_1_2b} that $G_{R, 2r}(\omega, t_0) $ is partially differentiable at least twice with respect to $t_0$.\\

$\bullet$ In Section~\ref{sec:Section_order_1}, it is shown that, \textbf{if} $G_{R}(\omega, t_0)=0$ at $\omega=\pm \omega_z(t_0)$, for each fixed choice of $t_0 \in \Re$ and $(2r+1)$ is the highest order of the zero at $\omega=\pm \omega_z(t_0)$  for some value of $r \in W$ (element of set of whole numbers including zero), \textbf{then} $ G_{R, 2r}(\omega, t_0) = \frac{\partial^{2r} G_{R}(\omega, t_0)}{\partial \omega^{2r}}= 0$ at $\omega= \pm \omega_z(t_0)$ and $\frac{\partial G_{R, 2r}(\omega, t_0)}{\partial \omega}=\frac{ \partial^{2r+1} G_{R}(\omega, t_0)}{\partial \omega^{2r+1}} \neq 0$ at $\omega= \pm \omega_z(t_0)$.\\

$\bullet$ It is shown in Section~\ref{sec:Section_A_1_8_1_c_0} that  $\omega_z(t_0)$ is a \textbf{continuous}  function of $t_0$ , for $t_0 \in \Re$ , using \textbf{Implicit Function Theorem} in $\Re^2$. 

.\\ \tocless\subsection{\label{sec:Section_A_1_8_1_2_0} \textbf{ $G_{R}(\omega, t_0) $ and $ G_{R, 2r}(\omega, t_0)$ are  partially differentiable twice as a function of $\omega$} \protect\\  \lowercase{} }

$G_{R}(\omega, t_0) $ in Eq.~\ref{sec:sec_2_1_eq_9a}  is copied below.
\begin{align}\label{sec:sec_a_1_8_1_1z}  
G_{R}(\omega, t_0) = e^{- \sigma t_0}  [\int_{-\infty}^{0}   E_0(\tau + t_0  ) e^{-2 \sigma \tau}  \cos{(\omega \tau)} d\tau 
+ \int_{-\infty}^{0}   E_0(\tau - t_0 )   \cos{(\omega \tau)} d\tau ] \notag\\
+ e^{\sigma t_0}  [\int_{-\infty}^{0}   E_0(\tau - t_0  )  e^{-2 \sigma \tau}  \cos{(\omega \tau)} d\tau 
+ \int_{-\infty}^{0}    E_0(\tau + t_0 )  \cos{(\omega \tau)} d\tau ] \notag\\ + 2 \cosh{( \sigma t_0)} \int_{-\infty}^{0}  E_0(\tau) ( e^{-2 \sigma \tau} + 1) \cos{(\omega \tau )} d\tau 
\end{align}
We write it succinctly as follows.  
\begin{align}\label{sec:sec_a_1_8_1_1}  
G_{R}(\omega, t_0) = \int_{-\infty}^{0} A(\tau, t_0) \cos{(\omega \tau)} d\tau \notag\\
A(\tau, t_0) = e^{- \sigma t_0} A_1(\tau, t_0) + e^{\sigma t_0}  A_1(\tau, -t_0) +  2 \cosh{( \sigma t_0)}  E_0(\tau) ( e^{-2 \sigma \tau} + 1) \notag\\
A_1(\tau, t_0) =E_0(\tau + t_0  ) e^{-2 \sigma \tau} +   E_0(\tau - t_0 )
\end{align}
We see that $E_{0}(\tau)$ in Eq.~\ref{sec_intro_eq_1} and its $t_0$  shifted versions are analytic functions of $\tau$ and $t_0$, given that the sum and product of exponential functions are analytic and hence infinitely differentiable.(\textbf{Result 4.1})\\ 

In Eq.~\ref{sec:sec_a_1_8_1_1}, $G_{R}(\omega, t_0) $ is partially differentiable at least twice with respect to $\omega$ and the integrals converge in Eq.~\ref{sec:sec_a_1_8_1_1} to Eq.~\ref{sec:sec_a_1_8_1_1a}  for $0 < \sigma < \frac{1}{2}$, because the term $ \tau^r E_{0}(\tau + M t_0) e^{-2 \sigma \tau}$ has \textbf{exponential} asymptotic fall-off rate as $|\tau| \to \infty$, for $r \in W$  (Details in Section~\ref{sec:appendix_C_5b}). The integrands in Eq.~\ref{sec:sec_a_1_8_1_1} to Eq.~\ref{sec:sec_a_1_8_1_1a} are analytic functions of variables $\omega$, $t_0$(using Result 4.1 in this section and given that the terms $ \cos{(\omega \tau)},  \sin{(\omega \tau)}$ and $ e^{-2 \sigma \tau}$ are analytic functions). The integrands are finite and have \textbf{exponential} asymptotic fall-off rate (Details in Section~\ref{sec:appendix_C_5b}) and absolutely integrable and we can find a suitable dominating function with exponential asymptotic fall-off rate which is absolutely integrable.(Details in Section~\ref{sec:appendix_C_5z1}) The partial derivative of the integrands in the equation for $G_{R}(\omega, t_0)$ in Eq.~\ref{sec:sec_a_1_8_1_1}, with respect to $\omega$ exists, as shown in the integrands  in Eq.~\ref{sec:sec_a_1_8_1_1az} and Eq.~\ref{sec:sec_a_1_8_1_1a}. Hence we can interchange the order of partial differentiation and integration  in Eq.~\ref{sec:sec_a_1_8_1_1az} and Eq.~\ref{sec:sec_a_1_8_1_1a} using theorem of differentiability of functions defined by Lebesgue integrals and theorem of dominated convergence, recursively as follows.(\href{https://archive.is/QjvGe}{theorem}) (\textbf{Result 4.1.A})
\begin{align}\label{sec:sec_a_1_8_1_1az}    
\frac{\partial G_{R}(\omega, t_0)}{\partial \omega} = - \int_{-\infty}^{0} \tau A(\tau, t_0) \sin{(\omega \tau)} d\tau
\end{align}
We write the second derivative $\frac{\partial^2 G_{R}(\omega, t_0)}{\partial \omega^2}$ as follows.
\begin{align}\label{sec:sec_a_1_8_1_1a}    
\frac{\partial^2 G_{R}(\omega, t_0)}{\partial \omega^2} = - \int_{-\infty}^{0} \tau^{2} A(\tau, t_0) \cos{(\omega \tau)} d\tau
\end{align}
We can use the arguments in the above paras and derive the $(2r)^{th}$ derivative of $G_{R}(\omega, t_0)$, for $r \in W$, which is differentiable at least twice, as follows.
\begin{align}\label{sec:sec_a_1_8_1_1ab}    
 G_{R, 2r}(\omega, t_0) = \frac{\partial^{2r} G_{R}(\omega, t_0)}{\partial \omega^{2r}} = (-1)^{r} \int_{-\infty}^{0} \tau^{2r} A(\tau, t_0) \cos{(\omega \tau)} d\tau
\end{align}
We can prove Eq.~\ref{sec:sec_a_1_8_1_1ab} using induction. We use Eq.~\ref{sec:sec_a_1_8_1_1ab} as Induction Hypothesis. We take the first and second derivative of Eq.~\ref{sec:sec_a_1_8_1_1ab} and we interchange the order of differentiation and integration, using the arguments used to derive Eq.~\ref{sec:sec_a_1_8_1_1a} as follows.
\begin{align}\label{sec:sec_a_1_8_1_1ab_1}
 \frac{\partial^{2r+1} G_{R}(\omega, t_0)}{\partial \omega^{2r+1}} = (-1)^{r+1}  \int_{-\infty}^{0} \tau^{2r+1} A(\tau, t_0) \sin{(\omega \tau)} d\tau
\end{align}
\begin{align}\label{sec:sec_a_1_8_1_1ab_1z}
 \frac{\partial^{2r+2} G_{R}(\omega, t_0)}{\partial \omega^{2r+2}} = (-1)^{r+1} \int_{-\infty}^{0} \tau^{2r+2} A(\tau, t_0) \cos{(\omega \tau)} d\tau
\end{align}
We see that the equation in Eq.~\ref{sec:sec_a_1_8_1_1ab_1z} is the same as the equation obtained by setting $r=r+1$ in Eq.~\ref{sec:sec_a_1_8_1_1ab}. Thus we have proved Eq.~\ref{sec:sec_a_1_8_1_1ab} using mathematical induction.

.\\ \tocless\subsection{\label{sec:appendix_C_5b} \textbf{  Exponential Fall off rate of  $ B(t)= t^r E_{0}(t + M t_0) e^{-2 \sigma t}$ for $r \in W$} \protect\\  \lowercase{} }

In this section, it is shown that the term $ B(t)= t^r E_{0}(t + M t_0) e^{-2 \sigma t}$, for integer $M$, has exponential asymptotic fall-off rate as $|t| \to \infty$, for $r \in W$. (\textbf{Result 4.2.1}).\\

We consider $C(t)= t^r   E_0(t+t_a) e^{-2 \sigma t}$ for real $t_a$. We see that $C(t-t_a)=(t-t_a)^r   e^{2 \sigma t_a} E_0(t) e^{-2 \sigma t}$. We see that $E_0(t)e^{-2 \sigma t}$ is an absolutely integrable function, for $0 \leq |\sigma| < \frac{1}{2}$ given that it has exponential fall-off rates as $|t| \to \infty$. (Details in ~\ref{sec:appendix_C_5} and ~\ref{sec:appendix_C_5z}).\\

Hence $C(t-t_a)=(t-t_a)^r   e^{2 \sigma t_a} E_0(t) e^{-2 \sigma t}$ also has exponential fall-off rates as $|t| \to \infty$, for  $r \in W$ and finite $t_a$ and is an absolutely integrable function.\\

Hence $C(t)= t^r   E_0(t+t_a) e^{-2 \sigma t}$ has exponential fall-off rates as $|t| \to \infty$, for finite $t_a$ and is an absolutely integrable function. We set $t_a=M t_0$ and see that $B(t)$ in Result 4.2.1, has \textbf{exponential fall-off rates} as $|t| \to \infty$, for finite $t_0$ and is an absolutely integrable function.

.\\ \tocless\subsection{\label{sec:appendix_C_5z1} \textbf{ Dominating function } \protect \\ \lowercase{} }

We consider $x(t)=E_0(t) e^{-2 \sigma t}$ which has asymptotic exponential fall-off rate of  $o[e^{-0.5|t|}]$.(Details in ~\ref{sec:appendix_C_5}) We see that $x(t+t_a)$ also has the same asymptotic exponential fall-off rate, for finite shift of $t_a=M t_0$ for integer $M$ and $y(t, t_a)=t^{r} x(t+t_a) e^{2 \sigma t_a}$ also has the same asymptotic exponential fall-off rate, for $r \in W$. We consider the intervals $0 \leq t_0 \leq t_{0_{max}}$ and $0 \leq t_a \leq t_{a_{max}}$ where $t_{0_{max}}, t_{a_{max}}$ are finite.\\

We consider finite $t_d >> t_{a_{max}}$ where $y(t, t_a)=t^{r} x(t+t_a) e^{2 \sigma t_a}$ falls off at the rate of $o[e^{0.5 t}]$ for $t << -t_d$. We consider $f(t, t_a,\omega)= y(t, t_a)  \cos{(\omega t)}$ and we get $\frac{\partial f(t, t_a,\omega)}{\partial \omega}= - t y(t, t_a)  \sin{(\omega t)} $ which falls off at the rate of $o[e^{0.5 t}]$ for $t << -t_d$. We see that $ x(t) = E_0(t) e^{-2 \sigma t}>0$ and finite, for $|t| < \infty$ (Details in ~\ref{sec:appendix_C_1}). Let $0 < f_{max} < \infty$ be the maximum value of $|\frac{\partial f(t, t_a,\omega)}{\partial \omega}|$ in the interval $-\infty < t < \infty$. \\ 

 We can find a suitable \textbf{dominating function} $D(t)= e^{-K |t|} f_{max} e^{K t_d}>0$ with a fall off rate of $O[e^{-K |t|}]$ where $0 < K< 0.5$ and hence $D(t)$ has a slower fall off rate than $\frac{\partial f(t, t_a,\omega)}{\partial \omega}$ and $D(t)=f_{max}$ at  $t = -t_d$ and hence $D(t)> |\frac{\partial f(t, t_a,\omega)}{\partial \omega}|$ for $-\infty < t \leq 0$ and hence $|\frac{\partial f(t, t_a,\omega)}{\partial \omega}| \leq D(t)$ in the interval $(-\infty, 0]$ and $\int_{-\infty}^{0} |D(t)| dt = \int_{-\infty}^{0}  e^{K t} f_{max} e^{K t_d} dt =\frac{1}{K} f_{max} e^{K t_d} [e^{K t}]_{-\infty}^{0}=  \frac{1}{K}  f_{max} e^{K t_d}  $ is finite.(\textbf{Result 4.3.1})\\
 
The terms in Eq.~\ref{sec:sec_a_1_8_1_1a} are given by $B(t)=  t^{r} e^{-2 \sigma t}  E_0(t+ M t_0)$ using Eq.~\ref{sec:sec_a_1_8_1_1} and Result 4.2.1 in Section~\ref{sec:appendix_C_5b}. We set $t_a=M t_0$ and get $B(t)= t^{r}  e^{-2 \sigma t} E_0(t+t_a) $. Hence $y(t, t_a)= t^{r} x(t+t_a) e^{2 \sigma t_a}= t^{r} E_0(t+t_a) e^{-2 \sigma t}$ in the second para, corresponds to $B(t)= t^{r}  e^{-2 \sigma t}  E_0(t+t_a) $ and Result 4.3.1 holds for this term.  We see that Result 4.3.1 holds for the  terms in Eq.~\ref{sec:sec_a_1_8_1_1a} and Eq.~\ref{sec:sec_a_1_8_1_1} using arguments in above paragraphs and setting $M=1, -1$  and setting $\sigma=0$ as needed.\\
 
 As $t_{0_{max}}, t_{a_{max}}$  increase to a larger and larger \textbf{finite value} without bounds, we consider larger intervals $0 \leq t_0 \leq t_{0_{max}}$ and $0 \leq t_a \leq t_{a_{max}}$ and $f_{max}$ and $t_d$ also increase correspondingly and the results in above paragraphs are valid in these intervals.\\ %
 
Similarly, we consider $f(t, t_a,\omega)= y(t, t_a)  \sin{(\omega t)}=t^{r} E_0(t+t_a) e^{-2 \sigma t}  \sin{(\omega t)}= t^{r} E_0(t+M  t_0) e^{-2 \sigma t}  \sin{(\omega t)}$ and we see that $\frac{\partial f(t, t_a,\omega)}{\partial t_0}$ which falls off at the rate of $o[e^{0.5 t}]$ for $t << -t_d$, using  Eq.~\ref{sec:sec_a_1_8_1_5} and $E_0(t)=E_0(-t)$ and due to the term $e^{-\pi n^{2} e^{-2t}}$ and we can use arguments in above paragraphs to get a result similar to Result 4.3.1 for the terms in Eq.~\ref{sec:sec_a_1_8_1_2}. We can use these arguments to get a result similar to Result 4.3.1 for the second derivative terms $\frac{\partial^{2} f(t, t_a,\omega)}{\partial t_0^{2}}$ in Eq.~\ref{sec:sec_a_1_8_1_3_a}. 

.\\ \tocless\subsection{\label{sec:Section_A_1_8_1_2b} \textbf{ $G_{R, 2r}(\omega, t_0)$ are partially differentiable twice as a function of $t_0$, $r \in W$} \protect\\  \lowercase{} }

We copy Eq.~\ref{sec:sec_a_1_8_1_1ab} and Eq.~\ref{sec:sec_a_1_8_1_1} as follows. 
\begin{align}\label{sec:sec_a_1_8_1_2z0}  
 G_{R, 2r}(\omega, t_0) = \frac{\partial^{2r} G_{R}(\omega, t_0)}{\partial \omega^{2r}} =  (-1)^{r} \int_{-\infty}^{0} \tau^{2r} A(\tau, t_0) \cos{(\omega \tau)} d\tau \notag\\
A(\tau, t_0) = e^{- \sigma t_0} A_1(\tau, t_0) + e^{\sigma t_0}  A_1(\tau, -t_0) +  2 \cosh{( \sigma t_0)}  E_0(\tau) ( e^{-2 \sigma \tau} + 1) \notag\\
A_1(\tau, t_0) =E_0(\tau + t_0  ) e^{-2 \sigma \tau} +   E_0(\tau - t_0 )
\end{align} 
In Eq.~\ref{sec:sec_a_1_8_1_2z0}, $G_{R, 2r}(\omega, t_0)$ is partially differentiable at least twice as a function of $t_0$ and the integrals converge in Eq.~\ref{sec:sec_a_1_8_1_2} shown as follows. Using Result 4.1.A in Section~\ref{sec:Section_A_1_8_1_2_0}, we can interchange the order of partial differentiation and integration  in Eq.~\ref{sec:sec_a_1_8_1_2z0} as follows.
\begin{align}\label{sec:sec_a_1_8_1_2} 
\frac{\partial G_{R, 2r}(\omega, t_0)}{\partial t_0} =    (-1)^{r} \int_{-\infty}^{0} \tau^{2r} \frac{\partial A(\tau, t_0)}{\partial t_0} \cos{(\omega \tau)} d\tau \notag\\
\frac{\partial A(\tau, t_0)}{\partial t_0} = - \sigma e^{- \sigma t_0} A_1(\tau, t_0) + e^{- \sigma t_0} \frac{\partial A_1(\tau, t_0)}{\partial t_0}  +  \sigma e^{\sigma t_0} A_1(\tau, -t_0) + e^{\sigma t_0} \frac{\partial A_1(\tau, -t_0)}{\partial t_0}  \notag\\
+  2 \sigma \sinh{( \sigma t_0)}  E_0(\tau) ( e^{-2 \sigma \tau} + 1) 
\end{align}
We show that the integrals in Eq.~\ref{sec:sec_a_1_8_1_2} converge, as follows. We see that the integrals for the first and third terms in Eq.~\ref{sec:sec_a_1_8_1_2} converge because the terms $ \tau^{2r}  E_{0}(\tau + M t_0) e^{-2 \sigma \tau}$ and $ \tau^{2r}  E_{0}(\tau + M t_0) $ in Eq.~\ref{sec:sec_a_1_8_1_2z0}, for integer $M$,  have exponential asymptotic fall-off rate as $|\tau| \to \infty$(Details in Section~\ref{sec:appendix_C_5b}). \\

We consider the integrand $\frac{\partial A_1(\tau, t_0)}{\partial t_0}$  in Eq.~\ref{sec:sec_a_1_8_1_2} first. We consider the term $E_0(\tau+ M t_0)$ for integer $M$ in Eq.~\ref{sec:sec_a_1_8_1_2z0} and can show that the corresponding integrals converge in Eq.~\ref{sec:sec_a_1_8_1_2}, as follows. We take the factor of $2$ out of the summation in $E_0(\tau)$ in Eq.~\ref{sec_intro_eq_1} copied below.
\begin{align}\label{sec:sec_a_1_8_1_4} 
E_{0}(\tau) = 2 \sum_{n=1}^{\infty}  [ 2 \pi^{2} n^{4} e^{4 \tau}    - 3 \pi n^{2}   e^{2 \tau} ]  e^{- \pi n^{2} e^{2 \tau}} e^{\frac{\tau}{2}}  \notag\\
E_0(\tau+ M t_0) = 2 \sum_{n=1}^{\infty}  [ 2 \pi^{2} n^{4} e^{4 \tau} e^{4 M t_0}   - 3 \pi n^{2}   e^{2 \tau} e^{2 M t_0} ]  e^{- \pi n^{2} e^{2 \tau} e^{2 M t_0}} e^{\frac{\tau}{2}} e^{\frac{M t_0}{2}} 
\end{align}
We can show that $ \frac{\partial}{\partial t_0}  E_0(\tau+ M t_0)  = M \frac{\partial}{\partial \tau}  E_0(\tau+ M t_0)  $ as follows, given that the equation for $E_0(\tau+ M t_0)$ in Eq.~\ref{sec:sec_a_1_8_1_4} has terms of the form $e^{\tau+ M t_0}$ and the equation is \textbf{invariant }if we interchange the variables $\tau$ and $M t_0$. (\textbf{Result 4.4.A})
\begin{align}\label{sec:sec_a_1_8_1_5} 
\frac{\partial}{\partial \tau}  E_0(\tau+ M t_0) = 2 \sum_{n=1}^{\infty}   e^{- \pi n^{2} e^{2 \tau} e^{2 M t_0}} e^{\frac{\tau}{2}} e^{\frac{ M t_0}{2}} [ 8 \pi^{2} n^{4} e^{4 \tau} e^{4 M t_0}   - 6 \pi n^{2}   e^{2 \tau} e^{2 M t_0} \notag\\ + ( \frac{1}{2} - 2 \pi n^{2} e^{2 \tau} e^{2 M t_0} ) (2 \pi^{2} n^{4} e^{4 \tau} e^{4 M t_0}   - 3 \pi n^{2}   e^{2 \tau} e^{2 M t_0})] \notag\\
\frac{\partial}{\partial t_0}  E_0(\tau+ M t_0)  = 2 \sum_{n=1}^{\infty}   e^{- \pi n^{2} e^{2 \tau} e^{2 M t_0 }} e^{\frac{\tau}{2}} e^{\frac{ M t_0}{2}} [ 8 M \pi^{2} n^{4} e^{4 \tau} e^{4 M t_0 }   - 6 M \pi n^{2}   e^{2 \tau} e^{2 M t_0} \notag\\+ ( \frac{M}{2} - 2 M \pi n^{2} e^{2 \tau} e^{2 M t_0} ) (2 \pi^{2} n^{4} e^{4 \tau} e^{4 M t_0}   - 3 \pi n^{2}   e^{2 \tau} e^{2 M t_0})]= M \frac{\partial}{\partial \tau}  E_0(\tau+ M t_0)  
\end{align}
We can write the integral corresponding to the term $E_0(\tau+ M t_0) e^{-2 \sigma \tau}$ in Eq.~\ref{sec:sec_a_1_8_1_2z0}, in the equation for $\frac{\partial A_1(\tau, t_0)}{\partial t_0}$  in Eq.~\ref{sec:sec_a_1_8_1_2}, using $M=1$ in Result 4.4.A,   as follows. We use the fact that $\int_{-\infty}^{0}  \frac{d A(\tau)}{d\tau} B(\tau) d\tau = \int_{-\infty}^{0}  \frac{d (A(\tau) B(\tau)) }{d\tau} d\tau - \int_{-\infty}^{0}  A(\tau) \frac{d B(\tau)}{d\tau}  d\tau$.
\begin{align}\label{sec:sec_a_1_8_1_6}   
\int_{-\infty}^{0}    \frac{\partial ( E_{0}(\tau +   t_0)   )}{\partial t_0} \tau^{2r} e^{-2 \sigma \tau} \cos{(\omega \tau)} d\tau =  \int_{-\infty}^{0} \frac{\partial  ( E_{0}(\tau +  t_0)  }{\partial \tau} \tau^{2r} e^{-2 \sigma \tau}  \cos{(\omega \tau)} d\tau  \notag\\
=  \int_{-\infty}^{0} \frac{\partial  ( E_{0}(\tau +   t_0) \tau^{2r}  e^{-2 \sigma \tau}  \cos{(\omega \tau)})}{\partial \tau} d\tau -\int_{-\infty}^{0}  E_{0}(\tau +   t_0)  \frac{\partial  ( \tau^{2r}  e^{-2 \sigma \tau}  \cos{(\omega \tau)}}{\partial \tau} d\tau \notag\\
= [ E_{0}(\tau +   t_0) \tau^{2r}  e^{-2 \sigma \tau}  \cos{(\omega \tau)} ]_{-\infty}^{0} +  \omega \int_{-\infty}^{0}  E_{0}(\tau +   t_0)  ) \tau^{2r} e^{-2 \sigma \tau}  \sin{(\omega \tau)} d\tau \notag\\+  2 \sigma \int_{-\infty}^{0}  E_{0}(\tau +   t_0)  ) \tau^{2r} e^{-2 \sigma \tau}  \cos{(\omega \tau)} d\tau -  2 r \int_{-\infty}^{0}  E_{0}(\tau +   t_0)  ) \tau^{2r-1} e^{-2 \sigma \tau}  \cos{(\omega \tau)} d\tau
\end{align}
We see that the integrals in Eq.~\ref{sec:sec_a_1_8_1_6} converge because the integrands are absolutely integrable because the terms $ E_{0}(\tau +   t_0)   \tau^{2r} e^{-2 \sigma \tau}  \sin{(\omega \tau)}$, $E_{0}(\tau +   t_0)   \tau^{2r} e^{-2 \sigma \tau}  \cos{(\omega \tau)} $ and $E_{0}(\tau +   t_0)  ) \tau^{2r-1} e^{-2 \sigma \tau}  \cos{(\omega \tau)} $ have exponential asymptotic fall-off rate as $|\tau| \to \infty$(Details in Section~\ref{sec:appendix_C_5b}). The term $[ E_{0}(\tau +   t_0)  \tau^{2r} e^{-2 \sigma \tau}  \cos{(\omega \tau)} ]_{-\infty}^{0} $ is finite, given that $ \tau^{2r} E_{0}(\tau)   e^{-2 \sigma \tau}$ and its shifted versions go to zero as $t \to -\infty$(Details in ~\ref{sec:appendix_C_5} ). Hence the integral $\int_{-\infty}^{0}    \frac{\partial  ( E_{0}(\tau +  t_0) \tau^{2r} e^{-2 \sigma \tau}  )}{\partial t_0} \cos{(\omega \tau)} d\tau$ in Eq.~\ref{sec:sec_a_1_8_1_6} and  in Eq.~\ref{sec:sec_a_1_8_1_2} corresponding to the term $ E_{0}(\tau +   t_0) e^{-2 \sigma \tau}$, converges.\\

We set $\sigma=0$ and $M=-1$ in the term $ E_{0}(\tau + M t_0)   e^{-2 \sigma \tau}$  and see that  the integral \\$\int_{-\infty}^{0}    \frac{\partial  ( E_{0}(\tau -  t_0)   )}{\partial t_0} \tau^{2r} \cos{(\omega \tau)} d\tau$ in Eq.~\ref{sec:sec_a_1_8_1_2} corresponding to the term  $ E_{0}(\tau -  t_0)$  also converges, using Result 4.4.A and the procedure used in Eq.~\ref{sec:sec_a_1_8_1_4} to Eq.~\ref{sec:sec_a_1_8_1_6}. Hence the integral for the second term $\frac{\partial A_1(\tau, t_0)}{\partial t_0}$  in Eq.~\ref{sec:sec_a_1_8_1_2} converges.\\

We set $M=-1$ in the term $ E_{0}(\tau + M t_0)   e^{-2 \sigma \tau}$ in Eq.~\ref{sec:sec_a_1_8_1_4} to Eq.~\ref{sec:sec_a_1_8_1_6} and see that  the integral $\int_{-\infty}^{0}    \frac{\partial  ( E_{0}(\tau -  t_0) e^{-2 \sigma \tau} ) }{\partial  t_0} \tau^{2r} \cos{(\omega \tau)} d\tau$ in Eq.~\ref{sec:sec_a_1_8_1_2} corresponding to the term $ E_{0}(\tau -  t_0)   e^{-2 \sigma \tau}$  also converges, using Result 4.4.A. \\

We set  $\sigma=0$ and set $M=1$ in the term $ E_{0}(\tau + M t_0)   e^{-2 \sigma \tau}$ and see that  the integral \\$\int_{-\infty}^{0}    \frac{\partial  ( E_{0}(\tau + t_0)   )}{\partial t_0} \tau^{2r} \cos{(\omega \tau)} d\tau$ in Eq.~\ref{sec:sec_a_1_8_1_2} corresponding to the term $ E_{0}(\tau + t_0) $  also converges, using Result 4.4.A and the procedure used in Eq.~\ref{sec:sec_a_1_8_1_4} to Eq.~\ref{sec:sec_a_1_8_1_6}.  Hence the integral for the fourth term $\frac{\partial A_1(\tau, -t_0)}{\partial t_0}$  in Eq.~\ref{sec:sec_a_1_8_1_2} converges.\\



The integral corresponding to the last term $ E_0(\tau) ( e^{-2 \sigma \tau} + 1) $ in Eq.~\ref{sec:sec_a_1_8_1_2} converges. Hence all the integrals in $\frac{\partial G_{R, 2r}(\omega, t_0)}{\partial t_0}$ in Eq.~\ref{sec:sec_a_1_8_1_2} converge.

.\\ \tocless\subsubsection{\label{sec:Section_A_1_8_1_c_01} \textbf{Second Partial Derivative of  $G_{R, 2r}(\omega, t_0) $ with respect to $t_0$} \protect  \lowercase{} }

The second partial derivative of  $G_{R, 2r}(\omega, t_0) $ with respect to $t_0$ is given by $\frac{\partial^2 G_{R, 2r}(\omega, t_0)}{\partial t_0^2} =  \frac{\partial}{\partial t_0} \frac{\partial G_{R, 2r}(\omega, t_0)}{\partial t_0}$ as follows. We use the result in Eq.~\ref{sec:sec_a_1_8_1_2}  and using Result 4.1.A in Section~\ref{sec:Section_A_1_8_1_2_0}, we can interchange the order of partial differentiation and integration  in Eq.~\ref{sec:sec_a_1_8_1_2} as follows.
\begin{align}\label{sec:sec_a_1_8_1_3_a} 
\frac{\partial^2 G_{R, 2r}(\omega, t_0)}{\partial t_0^2} =  (-1)^{r} \int_{-\infty}^{0} \tau^{2r} \frac{\partial^{2} A(\tau, t_0)}{\partial t_0^{2}} \cos{(\omega \tau)} d\tau \notag\\
\frac{\partial^{2} A(\tau, t_0)}{\partial t_0^{2}} =  \sigma^2 e^{- \sigma t_0}  A_1(\tau, t_0)  - 2 \sigma e^{- \sigma t_0} \frac{\partial A_1(\tau, t_0)}{\partial t_0} + e^{- \sigma t_0} \frac{\partial^{2} A_1(\tau, t_0)}{\partial t_0^{2}} \notag\\
+  \sigma^2 e^{\sigma t_0}  A_1(\tau, -t_0)  + 2 \sigma e^{\sigma t_0} \frac{\partial A_1(\tau, -t_0)}{\partial t_0} + e^{\sigma t_0} \frac{\partial^{2} A_1(\tau, -t_0)}{\partial t_0^{2}} \notag\\
+  2 \sigma^{2} \cosh{( \sigma t_0)}  E_0(\tau) ( e^{-2 \sigma \tau} + 1) 
\end{align}
The first two terms and fourth, fifth and seventh terms in $\frac{\partial^{2} A(\tau, t_0)}{\partial t_0^{2}}$ in Eq.~\ref{sec:sec_a_1_8_1_3_a} are similar to those in the equation for  $\frac{\partial G_{R, 2r}(\omega, t_0)}{\partial t_0}$ in Eq.~\ref{sec:sec_a_1_8_1_2} and their corresponding integrals have been shown to converge in Section~\ref{sec:Section_A_1_8_1_2b}. We will show that the integrals for third and sixth terms in Eq.~\ref{sec:sec_a_1_8_1_3_a} converge, as follows. We consider the integrand in the third term in Eq.~\ref{sec:sec_a_1_8_1_3_a}  first. We consider the term $E_0(\tau+M t_0)$ first in Eq.~\ref{sec:sec_a_1_8_1_2z0} and copy Eq.~\ref{sec:sec_a_1_8_1_4} below.
\begin{align}\label{sec:sec_a_1_8_1_4_a} 
E_{0}(\tau) = 2 \sum_{n=1}^{\infty}  [ 2 \pi^{2} n^{4} e^{4 \tau}    - 3 \pi n^{2}   e^{2 \tau} ]  e^{- \pi n^{2} e^{2 \tau}} e^{\frac{\tau}{2}}  \notag \\
E_0(\tau+ M t_0) = 2 \sum_{n=1}^{\infty}  [ 2 \pi^{2} n^{4} e^{4 \tau} e^{4 M t_0}   - 3 \pi n^{2}   e^{2 \tau} e^{2 M t_0} ]  e^{- \pi n^{2} e^{2 \tau} e^{2 M t_0}} e^{\frac{\tau}{2}} e^{\frac{M t_0}{2}} 
\end{align}
We showed in Result 4.4.A that $ \frac{\partial}{\partial t_0}  E_0(\tau+ M t_0)  = M \frac{\partial}{\partial \tau}  E_0(\tau+ M t_0)  $, given that the equation has terms of the form $e^{\tau+M t_0}$ and the equation \textbf{is invariant} if we interchange the variables $\tau$ and $M t_0$.\\

We use $t_0^{'}= M t_0$ in Eq.~\ref{sec:sec_a_1_8_1_4_a} and see that $\frac{\partial^2}{\partial (t_0^{'})^2} E_{0}(\tau+ t_0^{'}) = \frac{\partial^2}{\partial \tau^2}  E_{0}(\tau+ t_0^{'})$ (\textbf{Result 4.4.1.E'}) given that the equation has terms of the form $e^{\tau+t_0^{'}}$ and the equation \textbf{is invariant} if we interchange the variables $\tau$ and $t_0^{'}$. \\

Given that $\frac{\partial}{\partial t_0} = \frac{\partial}{\partial t_0^{'}}  \frac{\partial  t_0^{'}}{\partial t_0}   = M  \frac{\partial}{\partial t_0^{'}} $, we get $\frac{\partial^2}{\partial t_0^2} = \frac{\partial}{\partial t_0} ( \frac{\partial}{\partial t_0}) =  M \frac{\partial}{\partial t_0} ( \frac{\partial}{\partial t_0^{'}}) = M^{2} \frac{\partial}{\partial t_0^{'}} ( \frac{\partial}{\partial t_0^{'}}) = M^{2} \frac{\partial^2}{\partial (t_0^{'})^2}  $,  we substitute it in Result 4.4.1.E' and get $\frac{\partial^2}{\partial t_0^2} E_{0}(\tau+ M t_0) = M^{2} \frac{\partial^2}{\partial \tau^2}  E_{0}(\tau+M t_0)$ .(\textbf{Result 4.4.1.A'})\\

We can write the third term $\frac{\partial^{2} A_1(\tau, t_0)}{\partial t_0^{2}} $ in  Eq.~\ref{sec:sec_a_1_8_1_3_a}, corresponding to the term $E_0(\tau+ t_0) e^{-2 \sigma \tau}$, using $M=1$ in Result $4.4.1.A'$,   as follows. We use the fact that $\int_{-\infty}^{0}  \frac{d A(\tau)}{d\tau} B(\tau) d\tau = \int_{-\infty}^{0}  \frac{d (A(\tau) B(\tau)) }{d\tau} d\tau - \int_{-\infty}^{0}  A(\tau) \frac{d B(\tau)}{d\tau}  d\tau$.
\begin{align}\label{sec:sec_a_1_8_1_6_a}   
\int_{-\infty}^{0}    \frac{\partial^2 ( E_{0}(\tau + t_0)   )}{\partial t_0^2} \tau^{2r}  e^{-2 \sigma \tau} \cos{(\omega \tau)} d\tau =  \int_{-\infty}^{0} \frac{\partial^2  ( E_{0}(\tau + t_0)  )}{\partial \tau^2}\tau^{2r}   e^{-2 \sigma \tau}  \cos{(\omega \tau)} d\tau  \notag\\
= \int_{-\infty}^{0} \frac{\partial  ( \frac{\partial E_{0}(\tau + t_0)}{\partial \tau}    \tau^{2r}  e^{-2 \sigma \tau}  \cos{(\omega \tau)})}{\partial \tau} d\tau - \int_{-\infty}^{0}  \frac{\partial E_{0}(\tau + t_0)}{\partial \tau}   \frac{\partial  ( \tau^{2r}  e^{-2 \sigma \tau}  \cos{(\omega \tau))}}{\partial \tau} d\tau \notag\\
= [ \frac{\partial E_{0}(\tau + t_0)}{\partial \tau}  \tau^{2r}    e^{-2 \sigma \tau}  \cos{(\omega \tau)} ]_{-\infty}^{0} + \omega \int_{-\infty}^{0} \frac{\partial E_{0}(\tau + t_0)}{\partial \tau} \tau^{2r}  e^{-2 \sigma \tau}  \sin{(\omega \tau)} d\tau \notag\\+ 2 \sigma \int_{-\infty}^{0}  \frac{\partial E_{0}(\tau + t_0)}{\partial \tau}  \tau^{2r}    e^{-2 \sigma \tau}  \cos{(\omega \tau)} d\tau - 2 r \int_{-\infty}^{0}  \frac{\partial E_{0}(\tau + t_0)}{\partial \tau}  \tau^{2r-1}    e^{-2 \sigma \tau}  \cos{(\omega \tau)} d\tau
\end{align}
We see that the integrals $\int_{-\infty}^{0}  \frac{\partial E_{0}(\tau + t_0)}{\partial \tau}  \tau^{2r}      e^{-2 \sigma \tau}  \cos{(\omega \tau)} d\tau$ and $\int_{-\infty}^{0}  \frac{\partial E_{0}(\tau + t_0)}{\partial \tau}  \tau^{2r-1}    e^{-2 \sigma \tau}  \cos{(\omega \tau)} d\tau$ in Eq.~\ref{sec:sec_a_1_8_1_6_a}  converge, using Eq.~\ref{sec:sec_a_1_8_1_6} in the previous subsection. We see the term $ [\frac{\partial E_{0}(\tau + t_0)}{\partial \tau}     \tau^{2r}    e^{-2 \sigma \tau}  \cos{(\omega \tau)} ]_{-\infty}^{0} $ also converges, given that $E_0(\tau)=E_0(-\tau)$ and $E_{0}(\tau + t_0)=E_{0}(-\tau -  t_0)$ and we consider $\frac{\partial E_{0}(\tau +  t_0)}{\partial \tau}     \tau^{2r}    e^{-2 \sigma \tau}=\frac{\partial E_{0}(-\tau -  t_0)}{\partial \tau}     \tau^{2r}    e^{-2 \sigma \tau} $ using Eq.~\ref{sec:sec_a_1_8_1_5} and see that the term $e^{- \pi n^{2} e^{-2 \tau}} $ goes to zero faster than the rising term $ \tau^{2r}    e^{-2 \sigma \tau} e^{-6 \tau} e^{-\frac{\tau}{2}}$, as $\tau \to -\infty$. (\textbf{Result 4.4.1.1}) \\

It is shown below that the term $\int_{-\infty}^{0} \frac{\partial E_{0}(\tau +  t_0)}{\partial \tau} \tau^{2r}    e^{-2 \sigma \tau}  \sin{(\omega \tau)} d\tau$  in Eq.~\ref{sec:sec_a_1_8_1_6_a} also converges.
\begin{align}\label{sec:sec_a_1_8_1_6_a0}   
\int_{-\infty}^{0} \frac{\partial  ( E_{0}(\tau +  t_0)  )}{\partial \tau} \tau^{2r}    e^{-2 \sigma \tau}  \sin{(\omega \tau)} d\tau  \notag\\
= \int_{-\infty}^{0} \frac{\partial  ( E_{0}(\tau +  t_0) \tau^{2r}   e^{-2 \sigma \tau}  \sin{(\omega \tau)})}{\partial \tau} d\tau - \int_{-\infty}^{0}  E_{0}(\tau +  t_0)\frac{\partial  ( \tau^{2r} e^{-2 \sigma \tau}  \sin{(\omega \tau)}}{\partial \tau} d\tau \notag\\
= [ E_{0}(\tau +  t_0) \tau^{2r}  e^{-2 \sigma \tau}  \sin{(\omega \tau)} ]_{-\infty}^{0} - \omega \int_{-\infty}^{0}  E_{0}(\tau +  t_0)   \tau^{2r} e^{-2 \sigma \tau}  \cos{(\omega \tau)} d\tau \notag\\+ 2 \sigma \int_{-\infty}^{0}  E_{0}(\tau +  t_0)   \tau^{2r} e^{-2 \sigma \tau}  \sin{(\omega \tau)} d\tau - 2 r \int_{-\infty}^{0}  E_{0}(\tau +  t_0)   \tau^{2r-1} e^{-2 \sigma \tau}  \sin{(\omega \tau)} d\tau
\end{align}
We see that the integrals in Eq.~\ref{sec:sec_a_1_8_1_6_a0} converge because the integrands are absolutely integrable because the terms $ E_{0}(\tau +  t_0)   \tau^{2r} e^{-2 \sigma \tau}  \sin{(\omega \tau)}$, $E_{0}(\tau +  t_0)   \tau^{2r-1} e^{-2 \sigma \tau}  \sin{(\omega \tau)}$ and $E_{0}(\tau + t_0)   \tau^{2r} e^{-2 \sigma \tau}  \cos{(\omega \tau)} $  have exponential asymptotic fall-off rate as $|\tau| \to \infty$(Details in Section~\ref{sec:appendix_C_5b}). The term $[ E_{0}(\tau +  t_0)  \tau^{2r}  e^{-2 \sigma \tau}  \sin{(\omega \tau)} ]_{-\infty}^{0} $ is finite, given that $\tau^{2r} E_{0}(\tau)   e^{-2 \sigma \tau}$ and its shifted versions go to zero as $t \to -\infty$(Details in ~\ref{sec:appendix_C_5} ). Hence the integral $\int_{-\infty}^{0}    \frac{\partial^2  ( E_{0}(\tau +  t_0) \tau^{2r} e^{-2 \sigma \tau}  )}{\partial t_0^2} \cos{(\omega \tau)} d\tau$ in Eq.~\ref{sec:sec_a_1_8_1_6_a} and in  Eq.~\ref{sec:sec_a_1_8_1_3_a} corresponding to the term  $ E_{0}(\tau +  t_0)   e^{-2 \sigma \tau}$, also converges.\\

We set $\sigma=0$ and set $M=-1$ in the term $ E_{0}(\tau + M t_0)   e^{-2 \sigma \tau}$ and see that  the integral \\$\int_{-\infty}^{0}    \frac{\partial^2  ( E_{0}(\tau - t_0)   )}{\partial t_0^2} \tau^{2r} \cos{(\omega \tau)} d\tau$ in  Eq.~\ref{sec:sec_a_1_8_1_3_a} corresponding to the term  $ E_{0}(\tau - t_0)  $  also converges, using Result $4.4.1.A'$ and the procedure used in Eq.~\ref{sec:sec_a_1_8_1_4_a} to Eq.~\ref{sec:sec_a_1_8_1_6_a0}. Hence the third integral corresponding to the term $\frac{\partial^{2} A_1(\tau, t_0)}{\partial t_0^{2}} $ in  Eq.~\ref{sec:sec_a_1_8_1_3_a}, also converges.\\

We set $M=-1$ in the term $ E_{0}(\tau +  M t_0)   e^{-2 \sigma \tau}$ in Eq.~\ref{sec:sec_a_1_8_1_4_a} to Eq.~\ref{sec:sec_a_1_8_1_6_a0} and see that  the integral $\int_{-\infty}^{0}    \frac{\partial^2  ( E_{0}(\tau - t_0 ) \tau^{2r} e^{-2 \sigma \tau} ) }{\partial  t_0^2}  \cos{(\omega \tau)} d\tau$ in  Eq.~\ref{sec:sec_a_1_8_1_3_a} corresponding to the term  $ E_{0}(\tau - t_0)   e^{-2 \sigma \tau}$ also converges, using Result $4.4.1.A'$. \\ 

We set $\sigma=0$ and set $M=1$ in the term $ E_{0}(\tau + M t_0)   e^{-2 \sigma \tau}$ and see that  the integral \\$\int_{-\infty}^{0}    \frac{\partial^2  ( E_{0}(\tau +  t_0)   )}{\partial t_0^2} \tau^{2r} \cos{(\omega \tau)} d\tau$ in  Eq.~\ref{sec:sec_a_1_8_1_3_a} corresponding to the term  $ E_{0}(\tau +  t_0)$  also converges, using Result $4.4.1.A'$ and the procedure used in Eq.~\ref{sec:sec_a_1_8_1_4_a} to Eq.~\ref{sec:sec_a_1_8_1_6_a0}.  Hence the sixth integral corresponding to the term $\frac{\partial^{2} A_1(\tau, -t_0)}{\partial t_0^{2}} $ converges in  Eq.~\ref{sec:sec_a_1_8_1_3_a}.\\

 Hence all the integrals in Eq.~\ref{sec:sec_a_1_8_1_3_a} converge.


.\\ \tocless\subsection{\label{sec:Section_A_1_8_1_c_0} \textbf{Zero Crossings in $G_{R, 2r}(\omega, t_0) $ move continuously  as a function of $t_0$ , for $r \in W$ } \protect\\  \lowercase{} }

\textbf{Result 4.5.1:} It is shown in \textbf{Lemma 1} in Section~\ref{sec:Section_2_1} that  $G_{R}(\omega, t_0)=0$ at $\omega=\omega_z(t_0)$ where it crosses the zero line to the opposite sign, if Statement 1 is true. It is shown in Section~\ref{sec:Section_order_1} that $G_{R, 2r}(\omega, t_0) = 0 $ and $\frac{\partial G_{R, 2r}(\omega, t_0)}{\partial \omega} \neq 0 $ at $\omega=\omega_z(t_0)$, for some value of $r \in W$, where $(2r+1)$ is the highest order of the zero of $G_{R}(\omega, t_0)$ at $\omega=\omega_z(t_0)$. (\href{https://www.ocf.berkeley.edu/~araman/files/math_z/z_plot_1.jpg}{example plot}) \\

We use \textbf{Implicit Function Theorem} for the two dimensional case $\Re^2$ (\href{https://web.archive.org/web/20240218090425/https://www.math.ualberta.ca/~xinweiyu/217.1.13f/217-20131009.pdf}{link} and \href{https://archive.is/4zHBF}{link}). Given that $G_{R,2r}(\omega,  t_0) $ is partially differentiable with respect to $\omega$ and $t_0$ , with continuous partial derivatives, for $r \in W$ (Details in Section~\ref{sec:Section_A_1_8_1_2_0}, Section~\ref{sec:Section_A_1_8_1_2b}) and given that $G_{R,2r}(\omega, t_0) = 0$ at $\omega=\omega_z(t_0)$  and $\frac{\partial G_{R,2r}(\omega, t_0)}{\partial \omega} \neq 0 $ at $\omega=\omega_z(t_0)$, for some value of $r \in W$, where $(2r+1)$ is the highest order of the zero of $G_{R}(\omega, t_0)$ at $\omega=\omega_z( t_0)$ (using Result 4.5.1 in this section and using Section~\ref{sec:Section_order_1}), we see that $\omega_z(t_0)$ is a differentiable function of $t_0$  in the neighbourhood around $t_0$ .\\

We use Algorithm A in Section~\ref{sec:algo_A} and show that  $\omega_z(t_0)$ is a \textbf{continuous}  function of $t_0$ , for $t_0 \geq 0$. \\

We use the fact that $\omega_z(-t_0)=\omega_z(t_0)$ using Result 2.4.b in Section~\ref{sec:Section_A_1_2} and hence $\omega_z(t_0)$ is a \textbf{continuous}  function of $t_0$, for $t_0 \in \Re$. 

.\\ \tocless\subsection{\label{sec:Section_order_1} \textbf{Order of the zero in  $G_R(\omega, t_0) $ is finite.} \protect\\  \lowercase{} }

It is shown in this section that, \textbf{if} $G_{R}(\omega, t_0)=0$ at $\omega=\pm \omega_{z}(t_0)$ to satisfy Statement 1, for each fixed choice of $t_0 \in \Re$, \textbf{then} $ G_{R, 2r}(\omega, t_0) = \frac{\partial^{2r} G_{R}(\omega, t_0)}{\partial \omega^{2r}}= 0$ at $\omega= \pm \omega_{z}(t_0)$ and $\frac{\partial G_{R, 2r}(\omega, t_0)}{\partial \omega}=\frac{ \partial^{2r+1} G_{R}(\omega, t_0)}{\partial \omega^{2r+1}} \neq 0$ at $\omega= \pm \omega_{z}(t_0)$ for some value of $r \in W$ (element of set of whole numbers including zero) and $(2r+1)$ is the highest order of the zero of $G_{R}(\omega, t_0)$ at $\omega=\pm \omega_{z}(t_0)$ which is finite.\\

This is shown using Proof by Contradiction by assuming the \textbf{opposite} case that $\frac{\partial^{2r} G_{R}(\omega, t_0)}{\partial \omega^{2r}} = 0$ and $\frac{\partial^{2r+1} G_{R}(\omega, t_0)}{\partial \omega^{2r+1}} = 0$ at $\omega=\omega_{z}(t_0)$, for $r =0,1,...$, as $r \to \infty$ (\textbf{Statement D}) and show that it leads to a \textbf{contradiction}.\\

$G_{R, 2r}(\omega, t_0) $ in Eq.~\ref{sec:sec_a_1_8_1_2z0} is copied below. It is shown in \textbf{Lemma 1} in Section~\ref{sec:Section_2_1} that  $G_{R}(\omega, t_0)=0$ at $\omega=\omega_{z}(t_0)$ where it crosses the zero line to the opposite sign, \textbf{if} Statement 1 is true.
\begin{align}\label{sec:sec_order_eq_1}  
 G_{R, 2r}(\omega, t_0) = \frac{\partial^{2r} G_{R}(\omega, t_0)}{\partial \omega^{2r}} = (-1)^{r} \int_{-\infty}^{0} \tau^{2r} A(\tau, t_0) \cos{(\omega \tau)} d\tau \notag\\
A(\tau, t_0) = e^{- \sigma t_0} A_1(\tau, t_0) + e^{\sigma t_0}  A_1(\tau, -t_0) +  2 \cosh{( \sigma t_0)}  E_0(\tau) ( e^{-2 \sigma \tau} + 1) \notag\\
A_1(\tau, t_0) =E_0(\tau + t_0  ) e^{-2 \sigma \tau} +   E_0(\tau - t_0 )
\end{align} 
We compute the $(2r+1)^{th}$ and  $(2r+2)^{th}$ derivative of $G_{R}(\omega, t_0)$ and copy Eq.~\ref{sec:sec_a_1_8_1_1ab_1} and Eq.~\ref{sec:sec_a_1_8_1_1ab_1z} below.
\begin{align}\label{sec:sec_order_eq_2}    
 \frac{\partial^{2r+1} G_{R}(\omega, t_0)}{\partial \omega^{2r+1}} = (-1)^{r+1}  \int_{-\infty}^{0} \tau^{2r+1} A(\tau, t_0) \sin{(\omega \tau)} d\tau \notag\\
  \frac{\partial^{2r+2} G_{R}(\omega, t_0)}{\partial \omega^{2r+2}} = (-1)^{r+1} \int_{-\infty}^{0} \tau^{2r+2} A(\tau, t_0) \cos{(\omega \tau)} d\tau
\end{align}
We compute $C(\omega, t_0, \delta \omega) = \displaystyle\sum_{r=0}^{\infty} \frac{(\delta \omega)^{2r}}{!(2r)}   \frac{\partial^{2r} G_{R}(\omega, t_0)}{\partial \omega^{2r}} $ and $S(\omega, t_0, \delta \omega) = \displaystyle\sum_{r=0}^{\infty} \frac{(\delta \omega)^{2r+1} }{!(2r+1)}  \frac{\partial^{2r+1} G_{R}(\omega, t_0)}{\partial \omega^{2r+1}} $ below, using Eq.~\ref{sec:sec_order_eq_1} and Eq.~\ref{sec:sec_order_eq_2}, where $\delta \omega$ is real. It is noted that the integrals below converge, as $r \to \infty$, by taking the terms $\frac{1}{!(2r)}$ and $\frac{1}{!(2r+1)}$ inside the integrals.
\begin{align}\label{sec:sec_order_eq_3}    
C(\omega, t_0, \delta \omega) = \displaystyle\sum_{r=0}^{\infty} (-1)^{r} \frac{(\delta \omega)^{2r}}{!(2r)} [\int_{-\infty}^{0} \tau^{2r} A(\tau, t_0) \cos{(\omega \tau)} d\tau \notag\\
S(\omega, t_0, \delta \omega) = \displaystyle\sum_{r=0}^{\infty} (-1)^{r+1}\frac{(\delta \omega)^{2r+1}}{!(2r+1)} \int_{-\infty}^{0} \tau^{2r+1} A(\tau, t_0) \sin{(\omega \tau)} d\tau
\end{align}
We consider the integrals $C^{'}(\omega, t_0, \delta \omega) $ and $S^{'}(\omega, t_0, \delta \omega)$ obtained by interchanging the order of integration and summation  in Eq.~\ref{sec:sec_order_eq_3} as follows.
\begin{align}\label{sec:sec_order_eq_4}    
C^{'}(\omega, t_0, \delta \omega) =     \int_{-\infty}^{0} [ \displaystyle\sum_{r=0}^{\infty} (-1)^{r} \frac{(\delta \omega)^{2r}}{!(2r)} \tau^{2r}]  A(\tau, t_0)  \cos{(\omega \tau)} d\tau \notag\\
S^{'}(\omega, t_0, \delta \omega) =      \int_{-\infty}^{0} [\displaystyle\sum_{r=0}^{\infty} (-1)^{r+1} \frac{(\delta \omega)^{2r+1}}{!(2r+1)}  \tau^{2r+1}] A(\tau, t_0) \sin{(\omega \tau)} d\tau  
\end{align}
We use $\displaystyle\sum_{r=0}^{\infty} (-1)^{r} \frac{(\delta \omega)^{2r}}{!(2r)} \tau^{2r} = \cos{(  (\delta \omega) \tau)} $ and $\displaystyle\sum_{r=0}^{\infty} (-1)^{r+1} \frac{(\delta \omega)^{2r+1}}{!(2r+1)} \tau^{2r+1} =  -\sin{( (\delta \omega) \tau)} $ and write Eq.~\ref{sec:sec_order_eq_4} as follows. 
\begin{align}\label{sec:sec_order_eq_5}    
C^{'}(\omega, t_0, \delta \omega) =      \int_{-\infty}^{0}  A(\tau, t_0)  \cos{(  (\delta \omega) \tau)}   \cos{(\omega \tau)}  d\tau \notag\\
S^{'}(\omega, t_0, \delta \omega) = -     \int_{-\infty}^{0}   A(\tau, t_0) \sin{( (\delta \omega) \tau)}  \sin{(\omega \tau)}   d\tau 
\end{align}
The integrands in Eq.~\ref{sec:sec_order_eq_5} are absolutely integrable because the terms $ E_{0}(\tau +M t_0) e^{-2 \sigma \tau}$ and $ E_{0}(\tau +M t_0)$ in $A(\tau, t_0)$ in Eq.~\ref{sec:sec_order_eq_1}, for integer $M$,  have \textbf{exponential} asymptotic fall-off rate as $|\tau| \to \infty$ (Details in Section~\ref{sec:appendix_C_5b}) and the integrals  $C^{'}(\omega, t_0, \delta \omega) $ and $S^{'}(\omega, t_0, \delta \omega)$ converge. Hence we can interchange the order of integration and summation  in the integrals for $C^{'}(\omega, t_0, \delta \omega) $ and $S^{'}(\omega, t_0, \delta \omega)$ in Eq.~\ref{sec:sec_order_eq_4} using \textbf{Fubini's theorem} and we get  $C(\omega, t_0, \delta \omega) $ and $S(\omega, t_0, \delta \omega)$ in Eq.~\ref{sec:sec_order_eq_3} respectively. \href{https://archive.is/4h0fG}{(link)}  \\

Hence we see that $C(\omega, t_0, \delta \omega)=C^{'}(\omega, t_0, \delta \omega)$ and $S(\omega, t_0, \delta \omega)=S^{'}(\omega, t_0, \delta \omega)$ and interchanging the order of integration and summation in Eq.~\ref{sec:sec_order_eq_3} is justified.\\

We compute $C_S(\omega, t_0, \delta \omega)=C(\omega, t_0, \delta \omega)  + S(\omega, t_0, \delta \omega) =C^{'}(\omega, t_0, \delta \omega)+S^{'}(\omega, t_0, \delta \omega)$ using Eq.~\ref{sec:sec_order_eq_5} as follows, using the identity $ \cos{( (\omega+\delta \omega) \tau)}= \cos{(\omega \tau)} \cos{(  (\delta \omega) \tau)}   - \sin{(\omega \tau)}  \sin{( (\delta \omega) \tau)}$. 
\begin{align}\label{sec:sec_order_eq_6}    
C_S(\omega, t_0, \delta \omega) =     \int_{-\infty}^{0} A(\tau, t_0) \cos{( (\omega+\delta \omega) \tau)}  d\tau 
\end{align}
\textbf{If} Statement D is true, \textbf{then} $\frac{\partial^{2r} G_{R}(\omega, t_0)}{\partial \omega^{2r}} = 0$ and $\frac{\partial^{2r+1} G_{R}(\omega, t_0)}{\partial \omega^{2r+1}} = 0$ at $\omega=\omega_{z}(t_0)$ in Eq.~\ref{sec:sec_order_eq_1} and Eq.~\ref{sec:sec_order_eq_2}, for $r =0,1,...$, as $r \to \infty$, and $C(\omega, t_0, \delta \omega) = \displaystyle\sum_{r=0}^{\infty} \frac{(\delta \omega)^{2r}}{!(2r)}   \frac{\partial^{2r} G_{R}(\omega, t_0)}{\partial \omega^{2r}} =0$ and \\ $S(\omega, t_0, \delta \omega) = \displaystyle\sum_{r=0}^{\infty} \frac{(\delta \omega)^{2r+1} }{!(2r+1)}  \frac{\partial^{2r+1} G_{R}(\omega, t_0)}{\partial \omega^{2r+1}} =0$ at $\omega=\omega_{z}(t_0)$ in Eq.~\ref{sec:sec_order_eq_3} and hence \\ $C_S(\omega, t_0, \delta \omega)=C(\omega, t_0, \delta \omega)  + S(\omega, t_0, \delta \omega) =0$ at $\omega=\omega_{z}(t_0)$ in Eq.~\ref{sec:sec_order_eq_6}, shown as follows.
\begin{align}\label{sec:sec_order_eq_7a}    
C_S(\omega_{z}(t_0), t_0, \delta \omega)  =  \int_{-\infty}^{0} A(\tau, t_0) \cos{( (\omega_{z}(t_0)+\delta \omega) \tau)}  d\tau = 0 \notag\\
A(\tau, t_0) = e^{- \sigma t_0} A_1(\tau, t_0) + e^{\sigma t_0}  A_1(\tau, -t_0) +  2 \cosh{( \sigma t_0)}  E_0(\tau) ( e^{-2 \sigma \tau} + 1) \notag\\
A_1(\tau, t_0) =E_0(\tau + t_0  ) e^{-2 \sigma \tau} +   E_0(\tau - t_0 )
\end{align}
We expand Eq.~\ref{sec:sec_order_eq_7a} as follows, substituting $A(\tau, t_0)$ in the equation for $C_S(\omega_{z}(t_0), t_0, \delta \omega) $.
\begin{align}\label{sec:sec_order_eq_7}   
C_S(\omega_{z}(t_0), t_0, \delta \omega)  = e^{- \sigma t_0}  [\int_{-\infty}^{0}   E_0(\tau + t_0  ) e^{-2 \sigma \tau}  \cos{((\omega_{z}(t_0)+\delta \omega) \tau)} d\tau \notag\\
+ \int_{-\infty}^{0}   E_0(\tau - t_0 )   \cos{( (\omega_{z}(t_0)+\delta \omega)) \tau)} d\tau ] \notag\\
+ e^{\sigma t_0}  [\int_{-\infty}^{0}   E_0(\tau - t_0  )  e^{-2 \sigma \tau}  \cos{( (\omega_{z}(t_0)+\delta \omega) \tau)} d\tau 
+ \int_{-\infty}^{0}    E_0(\tau + t_0 )  \cos{( (\omega_{z}(t_0)+\delta \omega) \tau)} d\tau ] \notag\\ + 2 \cosh{( \sigma t_0)} \int_{-\infty}^{0}  E_0(\tau) ( e^{-2 \sigma \tau} + 1) \cos{( (\omega_{z}(t_0)+\delta \omega) \tau )} d\tau = 0
\end{align}
Eq.~\ref{sec:sec_order_eq_7} is similar to the equation for $G_{R}(\omega, t_0)$ in Eq.~\ref{sec:sec_2_1_eq_9a} in Section~\ref{sec:Section_A_1_2}, evaluated at $\omega=\omega_{z}(t_0)$ and equated to zero, to satisfy Statement 1, with $ \omega_{z}(t_0) $ replaced by $ \omega_{z}(t_0)+\delta \omega$.\\

Eq.~\ref{sec:sec_order_eq_7} holds for all real $\delta \omega$, including the case, as $\delta \omega \to 0$, \textbf{if} Statement D is true. This \textbf{contradicts} Result 2.1.h in Lemma 1 in Section~\ref{sec:Section_2_1} which requires $G_{R}(\omega, t_0)=0$ at $\omega=\omega_{z}(t_0)$ where it \textbf{crosses} the zero line to the \textbf{opposite sign}, to satisfy Statement 1. \\

Hence we see that,  \textbf{if} Statement 1 is true, \textbf{then} Statement D is false and hence there exists \textbf{at least one finite} $s \in N$ (element of set of natural numbers excluding zero) for which the $(s)^{th}$ derivative of $G_{R}(\omega, t_0)$ given by $G_{R, s}(\omega, t_0) = \frac{\partial^{s} G_{R}(\omega, t_0)}{\partial \omega^{s}} \neq 0$ at $\omega=\omega_{z}(t_0)$, where $s=2r$ is even or $s=2r+1$ is odd, for $r \in W$, for each fixed $t_0 \in \Re$.\\

We choose the \textbf{minimum} value of $s \in N$, for which $G_{R, s}(\omega, t_0) = \frac{\partial^{s} G_{R}(\omega, t_0)}{\partial \omega^{s}} \neq 0$ at $\omega=\omega_{z}(t_0)$ and hence $G_{R, s-1}(\omega, t_0) = \frac{\partial^{s-1} G_{R}(\omega, t_0)}{\partial \omega^{s-1}} = 0$ at $\omega=\omega_{z}(t_0)$ (\textbf{Result 4.6}). The first $(s-1)$ derivatives of $ G_{R}(\omega, t_0)$ equal zero at  $\omega=\omega_{z}(t_0)$ and hence  $s$ is the order of the zero of $G_{R}(\omega, t_0)$ at $\omega=\omega_{z}(t_0)$, and the order of this zero is \textbf{finite}. (\textbf{Result 4.6.1}) \\

 We see that $G_{R}(\omega, t_0)$ is an \textbf{even} function of $\omega$ because $g(t, t_0)$ is a real function of variable $t$ (Details in ~\ref{sec:appendix_I_2}) and hence \textbf{if} $G_{R}(\omega, t_0)$ has a zero at $\omega= + \omega_{z}(t_0)$, \textbf{then} it also has a zero at  $\omega= - \omega_{z}(t_0)$. (\textbf{Result 4.6.2}) \\ 
 
Using Result 4.6.1 and Result 4.6.2, we can write $G_{R}(\omega, t_0)=(\omega_{z}(t_0)^{2}- \omega^{2})^{s} N'(\omega, t_0)$ (\textbf{Eq.4.6}), for $s \in N$, where  $N'(\omega, t_0) \neq 0$ at $\omega= \pm \omega_{z}(t_0)$, for each fixed $t_0 \in \Re$ and $s$ is the highest order of the zero at $\omega=\omega_{z}(t_0)$ which is finite. ( It is noted that $\omega_{z}(t_0)$ represents the \textbf{zero crossing} in $G_{R}(\omega, t_0)$, for each fixed $t_0 \in \Re$. It is noted that $N'(\omega, t_0)$ may or may not be zero at $\omega \neq \pm \omega_{z}(t_0)$ and we \textbf{do not} claim otherwise.)\\

 The case of $s=2r$ in $G_{R}(\omega, t_0)=(\omega_{z}(t_0)^{2}- \omega^{2})^{s} N'(\omega, t_0)$ in Eq.4.6, is \textbf{ruled out} because $G_{R}(\omega, t_0)$ \textbf{changes sign} at $\omega= \pm \omega_{z}(t_0)$ (using Result 2.1.h in Lemma 1 in Section~\ref{sec:Section_2_1}), while $N'(\omega, t_0) \neq 0$ at $\omega= \pm \omega_{z}(t_0)$, \textbf{does not} change sign at $\omega= \pm \omega_{z}(t_0)$ and $(\omega_{z}(t_0)^{2}- \omega^{2})^{2r} \geq 0$ for real $\omega$, also \textbf{does not} change sign at $\omega= \pm \omega_{z}(t_0)$. (\textbf{Result 4.6.a}) \\

Hence  $G_{R}(\omega, t_0)=(\omega_{z}(t_0)^{2}- \omega^{2})^{2r+1} N'(\omega, t_0)$ in Eq.4.6 and $s=2r+1$ is the order of the zero of $G_{R}(\omega, t_0)$ at $\omega=\omega_{z}(t_0)$ and hence $G_{R, 2r+1}(\omega, t_0) = \frac{\partial^{2r+1} G_{R}(\omega, t_0)}{\partial \omega^{2r+1}} \neq 0$ at $\omega=\omega_{z}(t_0)$ and $G_{R, 2r}(\omega, t_0) = \frac{\partial^{2r} G_{R}(\omega, t_0)}{\partial \omega^{2r}} = 0$ at $\omega=\omega_{z}(t_0)$, using Result 4.6. \\

As an example, for the case $r=0$, we get $G_{R}(\omega, t_0)=(\omega_{z}(t_0)^{2}- \omega^{2}) N'(\omega, t_0)$ in Eq.4.6 and the order of the zero is $2r+1=1$ and we get $\frac{\partial G_{R}(\omega, t_0)}{\partial \omega} =(\omega_{z}(t_0)^{2}- \omega^{2}) \frac{ \partial N'(\omega, t_0)}{\partial \omega}  + N'(\omega, t_0) (-2 \omega)  \neq 0$  at $\omega=\omega_{z}(t_0)$, given that the first term is zero and the second term is not zero at $\omega=\omega_{z}(t_0)$ because $N'(\omega, t_0) \neq 0$ at $\omega= \pm \omega_{z}(t_0)$ and $\omega_{z}(t_0) \neq 0$. The above results hold for $r \in W$.\\

We have shown that, \textbf{if} $G_{R}(\omega, t_0)=0$ at $\omega=\pm \omega_{z}(t_0)$ to satisfy Statement 1, for each fixed choice of $t_0 \in \Re$, \textbf{then} $ G_{R, 2r}(\omega, t_0) = \frac{\partial^{2r} G_{R}(\omega, t_0)}{\partial \omega^{2r}}= 0$ at $\omega= \pm \omega_{z}(t_0)$ and $\frac{\partial G_{R, 2r}(\omega, t_0)}{\partial \omega}=\frac{ \partial^{2r+1} G_{R}(\omega, t_0)}{\partial \omega^{2r+1}} \neq 0$ at $\omega= \pm \omega_{z}(t_0)$ for some value of $r \in W$ and $(2r+1)$ is the highest order of the zero of $G_{R}(\omega, t_0)$ at $\omega=\pm \omega_{z}(t_0)$ which is finite.

.\\ \tocless\subsection{\label{sec:algo_A} \textbf{Algorithm A to find and track unique zero crossing function $\omega_z(t_0)$, for  $t_0 \geq 0$} \protect \\ \lowercase{} }

It is shown in \textbf{Lemma 1} in Section~\ref{sec:Section_2_1} that  $G_{R}(\omega, t_0)$ has \textbf{at least one} zero crossing at finite $\omega=\omega_z(t_0) >0$, where it crosses the zero line to the opposite sign, \textbf{if} Statement 1 is true. We see that $G_{R}(\omega, t_0)$ is an \textbf{even} function of $\omega$, because $g(t, t_0)$ is a real function of variable $t$ (\href{https://archive.is/CLVqX}{link}). If $G_{R}(\omega, t_0)$ has only one zero crossing to the right of origin, we get \textbf{at least two} zero crossings for $-\infty < \omega < \infty$. The algorithm below applies to $G_{R}(\omega, t_0)$ with one or more zero crossing to the right of origin.\\

$\bullet$  \textbf{Step A:} We consider the function $G_{R}(\omega, t_0)$ in the interval $-\omega_{max} \leq \omega \leq \omega_{max}$ and $0 \leq t_0 \leq t_{0_{max}}$, where $ t_{0_{max}}, \omega_{max}$ are real and arbitrarily large and finite. \\

We choose $\omega_{max}$ larger than the first zero crossing to the right of origin. We see that the zeros of $G_{R}(\omega, t_0)$  are isolated. (Details in Note 1 in Section~\ref{sec:Section_C_6_0z}) We compute the \textbf{minimum} separation between adjacent zero crossings in $G_{R}(\omega, t_0)$ given by $\delta_{\omega_{{min}}}$, in the interval $-\omega_{max} \leq \omega \leq \omega_{max}$, for \textbf{all possible} combinations of $t_0$ in the intervals $0 \leq t_0 \leq t_{0_{max}}$. (\textbf{Result 1}) \\

We see that $\delta_{\omega_{{min}}} >0$, given that the zeros of $G_{R}(\omega, t_0)$  are isolated (Details in Note 1 in Section~\ref{sec:Section_C_6_0z}). It is shown in Lemma 1 in Section~\ref{sec:Section_2_1} that $G_{R}(\omega, t_0)=0$ at $\omega=\omega_z(t_0)$ where it \textbf{crosses} the zero line to the opposite sign, if Statement 1 is true. So, adjacent zero crossings in $G_{R}(\omega, t_0)$ are separated by a \textbf{positive} value, for $ t_0 \in \Re $. (\textbf{Result 2})\\

$\bullet$  \textbf{Step 0:} We choose  \textbf{Segment $S_0$} in  $G_{R}(\omega, t_0)$ around the \textbf{first zero crossing} to the right of origin at $\omega_z(t_0)$, in the interval
$\omega_z(t_0) - \delta_{\omega} \leq \omega \leq \omega_z(t_0) + \delta_{\omega}  $, \textbf{evaluated at} $t_0= 0$, where $0 < \delta_{\omega} < \delta_{\omega_{{min}}}$. \\

We choose an arbitrarily small positive values of $ \delta_{t_0} << 1$ such that Implicit Function Theorem (IFT) can be used with this choice, to track this zero crossing at $\omega_z(t_0)$, in the intervals $0 \leq t_0 \leq t_{0_{max}}$, as detailed below. (See Note 2 in in Section~\ref{sec:Section_C_6_0z} for choice of $ \delta_{t_0}$)\\

We note that $ \omega_z(0) = \omega_{z_0} >0$ using Lemma 1.  We define $\omega_z(  \delta_{t_0} ) = \omega_{z_1} > 0$. We see that  $G_{R}(\omega, t_0)$ has \textbf{only one} zero crossing in Segment $S_0$, using Result 1 and 2 and  $0 < \delta_{\omega} < \delta_{\omega_{{min}}}$ in previous paras and we \textbf{back-off} and set $\delta_{\omega}=\frac{\delta_{\omega_{{min}}}}{4}$. (\textbf{Result $A_0$})\\

In this Segment $S_0$, we can apply IFT because there is \textbf{only one} zero crossing in this segment and $G_{R}(\omega, t_0)$ satisfies all criteria for IFT (\href{https://web.archive.org/web/20240218090425/https://www.math.ualberta.ca/~xinweiyu/217.1.13f/217-20131009.pdf}{link}), as shown in Section~\ref{sec:Section_A_1_8_1_2_0}, Section~\ref{sec:Section_A_1_8_1_2b} and Section~\ref{sec:Section_order_1} and hence $\omega_z(t_0)$ is a continuous function of  $t_0$ in the neighbourhood around $t_0=0$. (\textbf{Result $B_0$}) \\

We note that $| \omega_{z_0}  - \omega_{z_1}|  \leq \delta_{\omega}$ and $G_{R}(\omega, t_0)$  has \textbf{only one }zero crossing in Segment $S_0$ for all values of $0 \leq t_0 \leq  \delta_{t_0}$, using Result 1 and Result $A_0$. As we increase $t_0$, from $0$, over above intervals, at every infinitesimal step, we use Result $B_0$ and show that $\omega_z(t_0)$ is a continuous function of $t_0$ in the neighbourhood around that infinitesimal step. (\textbf{Result $C_0$})\\

We use IFT and see that $\omega_z(t_0)$ is a \textbf{continuous} function of $t_0$,  in the interval $0 \leq t_0 \leq  \delta_{t_0}$, using Result $C_0$. Hence we are able to \textbf{track} the unique first zero crossing, as it moves from $\omega_{z_0}$ to $\omega_{z_1}$. (\textbf{Result $D_0$})\\
 
$\bullet$  \textbf{Step n:} In Step 1 to N, we set $n=1,2,3,,,N$ and we \textbf{track} the zero crossing $\omega_{z_1}$ specified in Step 0, which is \textbf{unique}. We choose  \textbf{Segment $S_n$} in  $G_{R}(\omega, t_0)$ around the \textbf{unique zero crossing}  at $\omega_{z_{n}}$, in the interval $\omega_{z_{n}} - \delta_{\omega} \leq \omega \leq \omega_{z_{n}} + \delta_{\omega}  $.\\

We note that $ \omega_z( n \delta_{t_0} ) = \omega_{z_{n}} >0$ and $ \omega_z( (n+1) \delta_{t_0} ) = \omega_{z_{n+1}} >0$ using Lemma 1. We see that  $G_{R}(\omega, t_0)$ has \textbf{only one} zero crossing in Segment $S_n$, using Result 1 and 2 and $\delta_{\omega}=\frac{\delta_{\omega_{{min}}}}{4}$. (\textbf{Result $A_n$})\\ 

In Step 0, we chose $ \delta_{t_0}$ such that we can use IFT to track this unique zero crossing at $\omega_z(t_0)$ in the intervals $0 \leq t_0 \leq t_{0_{max}}$, as detailed below. (Details in Note 2 in Section~\ref{sec:Section_C_6_0z})\\

In this Segment $S_n$, we can apply Implicit Function Theorem (IFT) because there is \textbf{only one} zero crossing in this segment and $G_{R}(\omega, t_0)$ satisfies all criteria for IFT (\href{https://web.archive.org/web/20240218090425/https://www.math.ualberta.ca/~xinweiyu/217.1.13f/217-20131009.pdf}{link}), as shown in 
Section~\ref{sec:Section_A_1_8_1_2_0}, Section~\ref{sec:Section_A_1_8_1_2b} and Section~\ref{sec:Section_order_1}  and hence $\omega_z(t_0)$ is a continuous function of $t_0$ in the neighbourhood around $t_0= n \delta_{t_0}$. (\textbf{Result $B_n$}) \\
 
We note that $| \omega_{z_n}  - \omega_{z_{n+1}}|  \leq \delta_{\omega}$ and $G_{R}(\omega, t_0)$  has \textbf{only one }zero crossing in Segment $S_n$ for all values of $n \delta_{t_0} \leq t_0 \leq (n+1) \delta_{t_0}$, using Result 1 and Result $A_n$. As we increase $t_0$, from $ n \delta_{t_0}$, over above intervals, at every infinitesimal step, we use Result $B_n$ and show that $\omega_z(t_0)$ is a continuous function of  $t_0$ in the neighbourhood around that infinitesimal step. (\textbf{Result $C_n$})\\

We use IFT and see that $\omega_z(t_0)$ is a \textbf{continuous} function of  $t_0$,  in the interval $n \delta_{t_0} \leq t_0 \leq (n+1) \delta_{t_0}$, using Result $C_n$. Hence we are able to \textbf{track} the unique first zero crossing, as it moves from $\omega_{z_n}$ to $\omega_{z_{n+1}}$. (\textbf{Result $D_n$})\\

$\bullet$ We do not need to use induction hypothesis and induction step for Step $n$, given that
in this algorithm, Step $n$ from n=1,2,..N  has statements which use Result 1 and 2 derived before Step 1  and uses results derived within the \textbf{same} step $n$ and \textbf{not} previous steps $1,2,..n-1$.\\

$\bullet$  We use Step A and Step n for $n=1,2,...N$ where $N= \lceil \frac{t_{0_{max}}}{\delta_{t_0}} \rceil - 1$. \\

As an example, we choose $t_{0_{max}}= 10^{100}$ and $\omega_{max}=10^{100}$ and $\delta_{t_0}=10^{-100}$, where $\delta_{t_0}$ satisfies the conditions in Note 2. We use Step A and Step $n$ for $n=1,2,..N$, and using IFT,  we see that $\omega_z(t_0)$ is a \textbf{continuous} function of $t_0$, for $0 \leq t_0 \leq t_{0_{max}}$.\\

We can increase $t_{0_{max}}= 10^{1000}$ and $\omega_{max}=10^{1000}$ and $\delta_{t_0}=10^{-1000}$. Result 1 and 2 continue to hold and $\delta_{\omega_{{min}}} >0$. We use Step A and Step $n$ for $n=1,2,..N$ and using IFT, we see that $\omega_z(t_0)$ is a \textbf{continuous} function of $t_0$, for $0 \leq t_0 \leq t_{0_{max}}$.\\

As $t_{0_{max}}$ and $\omega_{max}$ increase to a larger and larger finite value without bounds, and $\delta_{t_0}$  decreases to a smaller and smaller positive value without bounds tending towards zero, Result 1 and 2 continue to hold and $\delta_{\omega_{{min}}} >0$ and hence $t_0$ spans the real interval $[0, \infty)$. Hence $\omega_z(t_0)$ is a continuous function of $t_0$, for  $ t_0 \geq 0 $.\\

In case the unique zero crossing reaches $\omega=0$, where $\omega_z(t_0)=0$, this case is detailed in Note 3 in Section~\ref{sec:Section_C_6_0z}.

.\\ \tocless\subsubsection{\label{sec:Section_C_6_0z} \textbf{ Note 1, 2 and 3.} \protect\\  \lowercase{} }

$\bullet$ \textbf{Note 1:}  It is shown in Section~\ref{sec:Section_order_1} that $\frac{\partial^{2r+1} G_{R}(\omega, t_0)}{\partial \omega^{2r+1}} \neq 0 $ at $\omega=\omega_z(t_0)$, for some value of $r \in W$ where $(2r+1)$ is the highest order of the zero of $G_{R}(\omega, t_0)$ at $\omega=\omega_z(t_0)$, if Statement 1 is true. This holds for every zero in   $G_{R}(\omega, t_0)$ and $ t_0 \in \Re$ (\textbf{Result 4.7.1.1}). \textbf{If} $G_{R}(\omega, t_0)$ had non-isolated zeros for some $\omega=\omega_z(t_0^{'})$, \textbf{then} we require $\frac{\partial^{2r+1} G_{R}(\omega, t_0)}{\partial \omega^{2r+1}}= 0 $ at that $\omega=\omega_z(t_0^{'})$, for all $r \in W$, which produces a \textbf{contradiction} of Result 4.7.1.1. Hence the zeros of $G_{R}(\omega, t_0)$  are \textbf{isolated}.\\
 
$\bullet$ \textbf{Note 2: Choice of $\delta_{t_{0}}$:} $G_{R}(\omega, t_0)$ satisfies all criteria for IFT (\href{https://web.archive.org/web/20240218090425/https://www.math.ualberta.ca/~xinweiyu/217.1.13f/217-20131009.pdf}{link}), as shown in 
Section~\ref{sec:Section_A_1_8_1_2_0}, Section~\ref{sec:Section_A_1_8_1_2b} Section~\ref{sec:Section_order_1} and hence we can apply IFT at each zero crossing $\omega_z(t_0)$, for each choice of $t_0$ in the intervals  $0 \leq t_0 \leq t_{0_{max}}$. For each choice of $t_0$ in the above intervals, we consider each of the zero crossings at $\omega_z(t_0)$, and compute the \textbf{maximum} value of the neighbourhood intervals around that choice of $t_0$, given by $\delta_{t_{0}}(k)$, for which IFT is valid. ($\delta_1$ in IFT proof for 2-D case corresponds to the  interval around $t_0$ given by $\delta_{t_0}(k)$ in our case and is small and positive. $\delta_1 \leq \delta_2$ in Point 2  in Line 5 of Proof section in \href{https://web.archive.org/web/20240218090425/https://www.math.ualberta.ca/~xinweiyu/217.1.13f/217-20131009.pdf}{Link to IFT proof},  where $\delta_2$ corresponds to the interval around $\omega$ in our case.) We clarify this point in next 2 paras, as follows.\\

In the 2-D case for IFT, we consider a function $f(x,y)$ which is partially differentiable with continuous partial derivatives and $f(x_0,y_0)=0$ and $a := \frac{\partial f(x_0, y_0)}{\partial y} \neq 0 $ at the zero crossing point at $y=y_0$ for $x=x_0$. We choose an open interval $I \times J  := (x_0 -  \delta_{1}, x_0 +  \delta_{1}) \times  (y_0 -  \delta_{2}, y_0 +  \delta_{2})$, where $\delta_{2}$ and $\delta_{1}$ are positive and small, satisfying the following conditions: \href{https://web.archive.org/web/20240218090425/https://www.math.ualberta.ca/~xinweiyu/217.1.13f/217-20131009.pdf}{Link to IFT proof}\\
1) We choose $\delta_2$ such that $\frac{\partial f(x, y)}{\partial y} \in (\frac{a}{2}, \frac{3a}{2})$ for all $\parallel (x,y) - (x_0, y_0 ) \parallel < 2 \delta_2$. \\
2) For the above $\delta_2$, we choose $\delta_{1} \leq \delta_2$ such that, for all $x \in (x_0 -  \delta_{1}, x_0 +  \delta_{1})  $, we satisfy $ |\frac{2 f(x,y_0)}{a}| < \delta_2$. \\

For any given $k$, we replace $\delta_{t_0}(k)$ by  $\delta_{t_0}$ to illustrate our 2-D case, similar to the 2-D case in above para. Let $a := \frac{\partial G_{R}(t_{0f},\omega_z(t_{0f}))}{\partial \omega} \neq 0 $ at the zero crossing point at $\omega=\omega_z(t_{0f})$ for $t_0=t_{0f}$. We choose an open interval $I \times J  :=   (t_{0f} -  \delta_{t_0}, t_{0f} +  \delta_{t_0}) \times  (\omega_z(t_{0f}) -  \delta_{2}, \omega_z(t_{0f}) +  \delta_{2}) $, where $  \delta_{t_0}$ is small and positive and using conditions 1 and 2 similar to 2-D case in above para, we get $0 < \delta_{t_0} \leq \delta_{t_{0_{max}}}$. We choose the \textbf{maximum} value of $\delta_{t_0}$ given by $ \delta_{t_{0_{max}}}$.\\

Using above para, we compute $\delta_{t_{0_{max}}}(k)$  for $k=1,2,...$, over all possible zero crossings and for all choices of $t_0$, in the intervals  $0 \leq t_0 \leq t_{0_{max}}$ and compute the \textbf{minimum} value of $\delta_{t_{0_{max}}}(k)$, over all $k$, given by $\delta_{t_{0_{min}}}$. We choose $0 < \delta_{t_0} \leq \delta_{t_{0_{min}}}$, for Algorithm A.\\

We can show that the minimum values given by $\delta_{t_{0_{min}}}$ is \textbf{non-zero} as follows, using proof by contradiction. We assume that $\delta_{t_{0_{min}}}=0$. This means that there is at least one value of $t_0= t_{0_{min}}$ in the interval  $0 \leq t_0 \leq t_{0_{max}}$, for one of the zero crossings, for which $\delta_{t_0}(k)=0$, for a specific $k=k'$, using IFT.($\delta_{t_0}(k)$ is the interval over which IFT results apply regarding differentiability of $\omega_z(t_0)$.) This \textbf{contradicts} Point 2  in Line 5 of Proof section in Page 1 of IFT proof,  which says $\delta_{t_0}(k)$ is small and positive and hence non-zero. ( $\delta_1$ in 2-D IFT proof corresponds to $\delta_{t_0}(k)$ in this section. \href{https://web.archive.org/web/20240218090425/https://www.math.ualberta.ca/~xinweiyu/217.1.13f/217-20131009.pdf}{Link to IFT proof}) Therefore, $\delta_{t_{0_{min}}}$ is \textbf{not }zero.\\

$\bullet$  \textbf{Note 3:} Algorithm A to find and track unique first zero crossing to the right of origin given by $\omega_z(t_0)$, for real $t_0$ is shown in this section in above paras. \textbf{In case} this unique zero crossing reaches $\omega=0$ for $t_0=t_{00}$, where $\omega_z(t_0)=0$, the algorithm \textbf{works till that point} and $\omega_z(t_0)$ is a continuous function of $t_0$, for $0 \leq t_0 < t_{00}$. (\textbf{Result 4.7.1.2}) \\

It is shown in Section~\ref{sec:Section_order_1} that $\frac{\partial^{s} G_{R}(\omega, t_0)}{\partial \omega^{s}} \neq 0 $ at $\omega=\omega_z(t_0)$, for some $s \in N$, for each choice of $ t_0$, \textbf{if} $G_{R}(\omega, t_0)=0$ at $\omega=\omega_z(t_0)$. This result holds even if Lemma 1 is not used and the result holds for $\omega_z(t_0)=0$ also. ( Eq.~\ref{sec:sec_order_eq_7a} in Section~\ref{sec:Section_order_1} holds for $\omega_z(t_0)=0$ and all real $\delta \omega$, including the case, as $\delta \omega \to 0$, \textbf{if} Statement D is true. This \textbf{contradicts} the fact that the zeros of $G_{R}(\omega, t_0)$  are \textbf{isolated}, shown in Note 1, given that $G_{R}(\omega, t_0)$ is \textbf{not} an all zero function of variable $\omega$ shown in Section~\ref{sec:Section_2_1_a0}. Hence Statement D is false.) More details in Note 3.2.\\

It is shown in Section~\ref{sec:Section_A_1_8_1_2_0}, Section~\ref{sec:Section_A_1_8_1_2b} that $G_{R, 2r}(\omega, t_0)$ is partially differentiable at least twice with respect to $\omega$, $t_0$ respectively. \\

Hence \textbf{if} the unique zero crossing $\omega_z(t_0)$ in Algorithm A reaches $\omega=0$ for $t_0=t_{00}$, where $\omega_z(t_0)=0$,  we \textbf{can} use IFT in the interval $t_{00} \leq t_0 \leq t_{00}+ \delta_{t_0}$, using the results in above 2 paras and show that $\omega_z(t_0)$ is a continuous function of $t_0$, in that interval. (\textbf{Result 4.7.1.3})\\

Now we proceed with the Algorithm A detailed in this section for $ t_0 > t_{00}+ \delta_{t_0}$. In case the unique zero crossing $\omega_z(t_0)$ goes to negative side, we choose the positive value, given that $G_{R}(\omega, t_0)$ is an even function of $\omega$. In case the unique zero crossing $\omega_z(t_0)$ reaches $\omega=0$ several more times, we use the argument in this Note 3, and show that $\omega_z(t_0)$ is a continuous function of  $t_0$, for $ t_0 \geq 0 $ using Result 4.7.1.2, Result 4.7.1.3 and this para.\\

\textbf{Note 3.1:} We note that it is shown in Lemma 1 in Section~\ref{sec:Section_2_1} that $G_{R}(\omega, t_0)$ must have\textbf{ at least one zero} at finite $\omega =  \omega_{z}(t_0) \neq 0$ where it \textbf{crosses} the zero line to the opposite sign, for each $t_0 \in \Re$, to satisfy \textbf{Statement 1}. In Algorithm A, we start with $t_0= 0$ and take the \textbf{first} zero crossing to the right of origin given by $\omega_{z}(t_0)$ and \textbf{track} it, as it moves along $0 \leq \omega < \infty$, for $t_0 \geq 0$ and show that this $\omega_z(t_0)$ is a continuous function of  $t_0$, which is required to set $\omega_z(t_{0c}) t_{0c}=  \pi$ in Section~\ref{sec:Section_A_1_3}.\\

\textbf{In case} this unique zero crossing reaches $\omega=0$ for $t_0=t_{00}$, where $\omega_z(t_0)=0$, it \textbf{does not} contradict the result $\omega_{z}(t_0) \neq 0$ in Lemma 1, given that we \textbf{do not} use Lemma 1 to get this result and instead $\omega_{z}(t_0) = 0$ is obtained during tracking of the first zero crossing to the right of origin, using Algorithm A.\\

It is shown in Section~\ref{sec:Section_C_6_0} that $\omega_z(t_0) t_0$ increases, as $t_0$ increases to a larger and larger finite value without bounds and $\omega_z(t_0) t_0 =  \pi$  can be reached, for specific  $t_0= t_{0c}$, as finite $t_0$ increases without bounds. Hence $\omega_{z}(t_0)$ and $\omega_z(t_0) t_0$ \textbf{do not} get stuck at zero, as $t_0$ increases to a larger and larger finite value without bounds.\\

\textbf{Note 3.2:} In the last 4 paras in Section~\ref{sec:Section_order_1}, it is stated that "The case of $s=2r$ in $G_{R}(\omega,t_0)=(\omega_z(t_0)^{2}- \omega^{2})^{s} N'(\omega, t_0)$  is \textbf{ruled out} because $G_{R}(\omega, t_0)$ \textbf{changes sign} at $\omega= \pm \omega_z(t_0)$ (using Result 2.1.h in Lemma 1 in Section~\ref{sec:Section_2_1})". This statement holds for the choice of $\omega_z(t_0) \neq 0$ only. For the case this unique first zero crossing reaches $\omega=0$, for some $t_0$, where $\omega_z(t_0)=0$, during tracking in Algorithm A, Lemma 1 is not used to derive $\omega_z(t_0) = 0$ and the zero crossing requirement vanishes and hence $s=2r$ is possible.\\

For the case this unique first zero crossing reaches $\omega=0$, for some $t_0=t_{00}$, $G_{R}(\omega, t_{00})$ and its even derivatives $G_{R, 2r}(\omega, t_{00}) = \frac{\partial^{2r} G_{R}(\omega, t_{00})}{\partial \omega^{2r}}$ are \textbf{even} function of $\omega$ and its odd derivatives $G_{R, 2r+1}(\omega, t_{00}) =\frac{\partial^{2r+1} G_{R}(\omega,  t_{00})}{\partial \omega^{2r+1}}$ are \textbf{odd} function of $\omega$, for all $r \in W$, given that $G_{R}(\omega,  t_0)$ is real and even function of $\omega$. Hence $G_{R, 2r+1}(\omega,  t_{00})=0$ at $\omega=0$, for all $r \in W$. \\

In general, $G_{R, 2r}(\omega,  t_{00})= \frac{\partial^{2r} G_{R}(\omega,  t_{00})}{\partial \omega^{2r}}$ may be non-zero, at $\omega=0$, for some $r \in W$, given that $G_{R}(\omega,  t_0)$ is \textbf{not} an all zero function of variable $\omega$, as shown in Section~\ref{sec:Section_2_1_a0}. This property enables us to use IFT at $\omega=0$, for the case this unique first zero crossing reaches $\omega=0$.

.\\ \tocless\section{\label{sec:Section_C_6_0} \textbf{ $\omega_z(t_0) t_0 = \pi $ can be reached for specific $t_0$} \protect\\  \lowercase{} }

In Lemma 1 in Section~\ref{sec:Section_2_1}, it is shown that $0 < \omega_z(t_0) < \infty$, for all  $ t_0 \in \Re $. For $t_0 >0$, we see that $\omega_z(t_0) t_0 > 0$. 
In Section~\ref{sec:Section_A_1_8_1}, it is shown that $\omega_z(t_0)$ is a \textbf{continuous}  function of variable $t_0$, for $ t_0 \in \Re $. Hence $\omega_z(t_0) t_0$ is a positive continuous function for $t_0 >0$.\\

We \textbf{require} $\omega_z(t_0) t_0 = \pi$, in Section 3, for a specific $t_0= t_{0c}$. 
In the section below, we consider $t_0 > 0$ and it is shown that $\omega_z(t_0) t_0$ increases to a larger and larger finite value, as $t_0$ increases to a larger and larger finite value without bounds. Hence, as we increase $t_0$, from zero, we can find a specific $t_0= t_{0c}$, for which $\omega_z(t_0) t_0 =  \pi$  is reached, using the method in Section~\ref{sec:Section_C_6_0a}.\\

Hence $\omega_z(t_0) t_0 = \pi$ can be reached, for specific  $t_0= t_{0c}$, as finite $t_0$ increases from zero, without bounds, given  that $\omega_z(t_0) t_0$ is a \textbf{continuous}  function of variable $t_0$, for all $ t_0 \in \Re $, using Intermediate Value Theorem.\\

In the section below, it is shown that, as $t_0$ increases to a larger and larger finite value without bounds, $\omega_z(t_0) t_0$ does not go to zero (Case 1) and does not approach a constant (Case 2) and instead, $\omega_z(t_0) t_0$ increases to a larger and larger finite value.

.\\ \tocless\subsection{\label{sec:Section_C_6_0a} \textbf{Method 1} \protect\\  \lowercase{} }

We consider $G_{R}(\omega, t_0) $ in Eq.~\ref{sec:sec_2_1_eq_9a} and set it to zero at $\omega=\omega_z(t_0) $ below.
\begin{align}\label{sec:sec_a_1_8_1_1z0}  
G_{R}(\omega_z( t_0), t_0) = e^{- \sigma t_0}  [\int_{-\infty}^{0}   E_0(\tau + t_0  ) e^{-2 \sigma \tau}  \cos{(\omega_z( t_0) \tau)} d\tau 
+ \int_{-\infty}^{0}   E_0(\tau - t_0 )   \cos{(\omega_z( t_0) \tau)} d\tau ] \notag\\
+ e^{\sigma t_0}  [\int_{-\infty}^{0}   E_0(\tau - t_0  )  e^{-2 \sigma \tau}  \cos{(\omega_z( t_0) \tau)} d\tau 
+ \int_{-\infty}^{0}    E_0(\tau + t_0 )  \cos{(\omega_z( t_0) \tau)} d\tau ] \notag\\ + 2 \cosh{( \sigma t_0)} \int_{-\infty}^{0}  E_0(\tau) ( e^{-2 \sigma \tau} + 1) \cos{(\omega_z( t_0) \tau )} d\tau = 0
\end{align}
We combine first and fourth integrals in Eq.~\ref{sec:sec_a_1_8_1_1z0}  and write in the first line in Eq.~\ref{sec:sec_a_1_8_1_1z}. We combine third and second integrals in Eq.~\ref{sec:sec_a_1_8_1_1z0}  and write in the second line in Eq.~\ref{sec:sec_a_1_8_1_1z} as follows. We use $D(t_0) = D(-t_0) = \cosh{( \sigma t_0)}  \int_{-\infty}^{0}  E_0(\tau) ( e^{-2 \sigma \tau} + 1) \cos{(\omega_z( t_0) \tau )} d\tau$ using Result 2.4.c in Section~\ref{sec:Section_A_1_2}. 
\begin{align}\label{sec:sec_a_1_8_1_1z}  
G_{R}(\omega_z( t_0), t_0) = e^{- \sigma t_0}  \int_{-\infty}^{0}   E_0(\tau + t_0  ) e^{-2 \sigma \tau}  \cos{(\omega_z( t_0) \tau)} d\tau + e^{\sigma t_0} \int_{-\infty}^{0}    E_0(\tau + t_0 )  \cos{(\omega_z( t_0) \tau)} d\tau  \notag\\
+ e^{\sigma t_0}  \int_{-\infty}^{0}   E_0(\tau - t_0  )  e^{-2 \sigma \tau}  \cos{(\omega_z( t_0) \tau)} d\tau + e^{- \sigma t_0} \int_{-\infty}^{0}   E_0(\tau - t_0 )   \cos{(\omega_z( t_0) \tau)} d\tau  
 + 2 D(t_0) = 0
\end{align}
We substitute  $\tau + t_0= t$ in the first two integrals in Eq.~\ref{sec:sec_a_1_8_1_1z} and get $\tau = t -  t_0$ and $d\tau = dt$ and we substitute back $t=\tau$ and write in the first line in Eq.~\ref{sec:sec_x_9az}. We substitute  $\tau - t_0= t$ in the third and fourth integrals in Eq.~\ref{sec:sec_a_1_8_1_1z} and get $\tau = t + t_0$ and $d\tau = dt$ and we substitute back $t=\tau$ and write in the second line in Eq.~\ref{sec:sec_x_9az} as follows. We use $e^{- \sigma t_0} e^{2 \sigma t_0} = e^{ \sigma t_0}$ in the first integral and $e^{ \sigma t_0} e^{-2 \sigma t_0} = e^{ -\sigma t_0}$ in the third integral, similar to the substitution in Eq.~\ref{sec:sec_az3_1_3_eq_1}.
\begin{align}\label{sec:sec_x_9az}  
G_{R}( \omega_z( t_0), t_0)  =    e^{\sigma t_0}  \int_{-\infty}^{ t_0}    E_0(\tau) e^{-2 \sigma \tau} \cos{( \omega_{z}(t_0) (\tau- t_0)} d\tau 
+  e^{ \sigma t_0}  \int_{-\infty}^{ t_0}   E_0(\tau) \cos{( \omega_{z}(t_0) (\tau - t_0)}      d\tau \notag\\
+ e^{-\sigma t_0}  \int_{-\infty}^{- t_0}    E_0(\tau) e^{-2 \sigma \tau} \cos{( \omega_{z}(t_0) (\tau + t_0)} d\tau 
+  e^{ -\sigma t_0}  \int_{-\infty}^{- t_0}   E_0(\tau) \cos{( \omega_{z}(t_0) (\tau + t_0)}      d\tau  + 2 D(t_0) = 0
\end{align}
We split the range of integration $\int_{-\infty}^{t_0} = \int_{-\infty}^{-t_0}  + \int_{-t_0}^{t_0}$. We compute the maximum value of $E_0(\tau) e^{-2 \sigma \tau}$, given by $E_{max}$, for real $\tau$ and $0 < \sigma < \frac{1}{2}$ and compute the minimum value $Min(E_{max},1)$ given by $E_{min}$. We choose sufficiently large $t_0 > 0$ such that $E_0(\tau) e^{-2 \sigma \tau} << E_{min} $ and $E_0(\tau) << E_{min}$ in the interval $|\tau| \geq t_0$. \\

 We take the $ \int_{-t_0}^{t_0} $ part of first and second integrals in Eq.~\ref{sec:sec_x_9az} and combine them to get the \textbf{dominant} term $F(t_0)$ in Eq.~\ref{sec:sec_x_9a1}. We use $2 D(t_0) = 2 \cosh{( \sigma t_0)}  \int_{-\infty}^{0}  E_0(\tau) ( e^{-2 \sigma \tau} + 1) \cos{(\omega_z( t_0) \tau )} d\tau$ in Result 2.4.c in Section~\ref{sec:Section_A_1_2} and use $2 \cosh{( \sigma t_0)} = e^{\sigma t_0} + e^{-\sigma t_0}$ and get $2 D(t_0) = ( e^{\sigma t_0}  + e^{-\sigma t_0}) \int_{-\infty}^{0}  E_0(\tau) ( e^{-2 \sigma \tau} + 1) \cos{(\omega_z( t_0) \tau )} d\tau$ and then split it into two integrals $ \int_{-\infty}^{0}  = \int_{-\infty}^{-t_0}+  \int_{-t_0}^{0} $ and take the integral $e^{\sigma t_0} \int_{-t_0}^{0}$ in $2 D(t_0)$ and combine it with the \textbf{dominant} term $F(t_0)$ in Eq.~\ref{sec:sec_x_9a1}. \\

Then we club the rest of the integrals in Eq.~\ref{sec:sec_x_9az} in the \textbf{negligible} term $X(t_0)$ and write as follows. It is shown in Section~\ref{sec:Section_C_6_2s} that $|X(t_0)|<< |F(t_0)|$ and $|X(t_0)|<< 1$, as $t_0$ increases.
\begin{align}\label{sec:sec_x_9a1}  
G_{R}( \omega_z( t_0), t_0)  = F_X(t_0)=  F(t_0)  + X(t_0) = 0    \notag\\ 
F(t_0) = e^{\sigma t_0}   [ \int_{-t_0}^{t_0}    E_0(\tau) [ e^{-2 \sigma \tau} + 1 ] \cos{( \omega_{z}(t_0) (\tau- t_0)} d\tau +  \int_{-t_0}^{0}  E_0(\tau) ( e^{-2 \sigma \tau} + 1) \cos{(\omega_z( t_0) \tau )}  d\tau ] \notag\\
X(t_0)=   e^{\sigma t_0}    \int_{-\infty}^{-t_0}    E_0(\tau) [ e^{-2 \sigma \tau} + 1 ] \cos{( \omega_{z}(t_0) (\tau- t_0)} d\tau  \notag\\
+ e^{-\sigma t_0}  \int_{-\infty}^{- t_0}    E_0(\tau) e^{-2 \sigma \tau} \cos{( \omega_{z}(t_0) (\tau + t_0)} d\tau 
+  e^{ -\sigma t_0}  \int_{-\infty}^{- t_0}   E_0(\tau) \cos{( \omega_{z}(t_0) (\tau + t_0)}      d\tau \notag\\
+ e^{\sigma t_0}  \int_{-\infty}^{-t_0}  E_0(\tau) ( e^{-2 \sigma \tau} + 1) \cos{(\omega_z( t_0) \tau )} d\tau 
+ e^{-\sigma t_0}  \int_{-\infty}^{0}  E_0(\tau) ( e^{-2 \sigma \tau} + 1) \cos{(\omega_z( t_0) \tau )} d\tau
\end{align}
It is noted that Eq.~\ref{sec:sec_x_9a1} and the equations in Section~\ref{sec:Section_C_6_2s1} and Section~\ref{sec:Section_C_6_2s2} below, apply to any real $t_0$.
In Section~\ref{sec:Section_C_6_2s1} and Section~\ref{sec:Section_C_6_2s2}, it is shown that, as $t_0$ increases to a larger and larger finite value without bounds, $\omega_z(t_0) t_0$ does not approach zero (Case 1) and does not approach a constant (Case 2). We \textbf{rule out} Case 1 and Case 2, by showing that $|F(t_0)| > 0$ and $|X(t_0)| << |F(t_0)|$ and that $F_X(t_0)= F(t_0)  + X(t_0)$ in Eq.~\ref{sec:sec_x_9a1} increases \textbf{exponentially}, as $t_0$ increases to a larger and larger value, while the right hand side in $F_X(t_0)$ in Eq.~\ref{sec:sec_x_9a1} remains zero, producing a \textbf{contradiction}. Hence we conclude that $\omega_z(t_0) t_0$ increases to a larger and larger finite value, as $t_0$ increases to a larger and larger finite value without bounds. 

.\\ \tocless\subsubsection{\label{sec:Section_C_6_2s1} \textbf{ $\omega_z(t_0) t_0$ \textbf{does not} approach zero, as  $t_0$ increases  } \protect \\ \lowercase{} }

\textbf{Case 1:} We assume that $\omega_{z}(t_0)= o[\frac{1}{t_0}]$ \href{https://www.stat.cmu.edu/~cshalizi/uADA/13/lectures/app-b.pdf}{((Little o Notation) )} (\textbf{Statement E})  and hence becomes smaller and smaller approaching zero and  $\omega_{z}(t_0) t_0$ becomes smaller and smaller approaching zero, as $t_0$ increases to a larger and larger finite value without bounds(\textbf{Result 5.1.1.a}). It is shown in this section that Statement E is false, if Statement 1 is true.  We copy Eq.~\ref{sec:sec_x_9a1} below.
\begin{align}\label{sec:sec_C_eq_1}  
F_X(t_0)= F(t_0)  + X(t_0) = 0    \notag\\ 
F(t_0) = e^{\sigma t_0}   [ \int_{-t_0}^{t_0}    E_0(\tau) ( e^{-2 \sigma \tau} + 1 ) \cos{( \omega_{z}(t_0) (\tau- t_0))} d\tau +  \int_{-t_0}^{0}  E_0(\tau) ( e^{-2 \sigma \tau} + 1) \cos{(\omega_z( t_0) \tau )}  d\tau ]
\end{align}
We find the \textbf{minimum} value of $F(t_0)$ as follows. We use  $\omega_{z}(t_0)= o[\frac{1}{t_0}]$ in Eq.~\ref{sec:sec_C_eq_1}. For sufficiently large $t_0$, we  get $\omega_z( t_0) \tau << 1 $  in the interval $|\tau| \leq t_0$ and hence $\cos{(\omega_z( t_0) \tau )}$ very close to unity. We use $\cos(x) \geq 1 - \frac{x^{2}}{2}$ for $x \in \Re$ \href{https://archive.is/mQSG0}{(link)} and use a loose bound. Hence $\cos{(\omega_z( t_0) \tau )} > \frac{1}{2}$ (\textbf{Result 5.1.1.b}).  \\

We use $ |\omega_{z}(t_0) (\tau- t_0)| = |\omega_{z}(t_0) t_0 ( \frac{\tau}{t_0} - 1) |= | \omega_{z}(t_0) t_0 (1 - \frac{\tau}{t_0})|  << 1 $, for sufficiently large $t_0$, in the interval $|\tau| \leq t_0$ and hence $\cos{( \omega_{z}(t_0) (\tau- t_0))}$ very close to unity and hence $\cos{(\omega_z( t_0) (\tau- t_0) )} > \frac{1}{2}$. (\textbf{Result 5.1.1.c})\\

For sufficiently large $t_0$, we write Eq.~\ref{sec:sec_C_eq_1} as follows, using Result 5.1.1.b and Result 5.1.1.c. We use $\int_{-t_0}^{t_0} = \int_{-t_0}^{0}   + \int_{0}^{t_0}$  below.
\begin{align}\label{sec:sec_C_eq_2}  
F(t_0) >  \frac{1}{2} e^{\sigma t_0}  [  \int_{-t_0}^{t_0}    E_0(\tau) ( e^{-2 \sigma \tau} + 1 )  d\tau + \int_{-t_0}^{0}  E_0(\tau) ( e^{-2 \sigma \tau} + 1)  d\tau ] \notag\\
F(t_0) >  \frac{1}{2} e^{\sigma t_0}   [ \int_{-t_0}^{0}    E_0(\tau) ( e^{-2 \sigma \tau} + 1 )  d\tau + \int_{0}^{t_0}    E_0(\tau) ( e^{-2 \sigma \tau} + 1 )  d\tau +  \int_{-t_0}^{0}  E_0(\tau) ( e^{-2 \sigma \tau} + 1)  d\tau ]
\end{align}
We substitute $\tau=-\tau$ in the integral $\int_{-t_0}^{0}$ and use $E_0(-\tau) =E_0(\tau) $ and $-\int_{t_0}^{0} = \int_{0}^{t_0} $ as follows.
\begin{align}\label{sec:sec_C_eq_3}  
 F(t_0) >  \frac{1}{2} e^{\sigma t_0}   [- \int_{t_0}^{0}    E_0(\tau) ( e^{2 \sigma \tau} + 1 )  d\tau + \int_{0}^{t_0}    E_0(\tau) ( e^{-2 \sigma \tau} + 1 )  d\tau 
 -  \int_{t_0}^{0}  E_0(\tau) ( e^{2 \sigma \tau} + 1)  d\tau ] \notag\\
F(t_0) >  \frac{1}{2} e^{\sigma t_0}   [ \int_{0}^{t_0}    E_0(\tau) ( e^{2 \sigma \tau} + 1 )  d\tau + \int_{0}^{t_0}    E_0(\tau) ( e^{-2 \sigma \tau} + 1 )  d\tau 
 +   \int_{0}^{t_0}  E_0(\tau) ( e^{2 \sigma \tau} + 1)  d\tau ] \notag\\
F(t_0) >  \frac{1}{2}  e^{\sigma t_0}   \int_{0}^{t_0}    E_0(\tau) ( 3 + 2 \cosh{(2 \sigma \tau)} +   e^{2 \sigma \tau} ) d\tau  = F_{min}(t_0)
\end{align}
We see that $F(t_0) > F_{min}(t_0)$ converges to a \textbf{finite positive} value, which increases exponentially, given that the integrand is positive, for positive $t_0$ and $E_0(\tau) e^{ \pm 2 \sigma \tau} >0$ (Details in ~\ref{sec:appendix_C_1}) and the integral converges as $t_0$ increases to $t_0' >> t_0$, as shown in Section~\ref{sec:Section_C_6_2s0a}. (\textbf{Result 5.1.1.d})\\

We choose sufficiently large $t_0$ such that $|X(t_0) | << F_{min}(t_0)$ is arbitrarily small. As we increase $t_0$ to $t_0' >> t_0$, we see that $F_{min}(t_0') > F_{min}(t_0)$ and  $|X(t_0') | << |X(t_0) | << F_{min}(t_0')$ (\textbf{Result 5.1.1.e}) \\

Hence $F_X(t_0)= F(t_0)  + X(t_0) = 0 $ is \textbf{not} possible in Eq.~\ref{sec:sec_C_eq_1}, which is required by Statement E. Hence Statement E is \textbf{false}, if Statement 1 is true and $\omega_z(t_0) t_0$ \textbf{does not} approach zero, as  $t_0$ increases to a larger and larger finite value without bounds.

.\\ \tocless\subsubsection{\label{sec:Section_C_6_2s2} \textbf{ $\omega_z(t_0) t_0$ \textbf{does not} approach a non-zero constant, as  $t_0$ increases  } \protect \\ \lowercase{} }

 \textbf{Case 2:} We assume that $\omega_z(t_0)= O[\frac{1}{t_0}]$ (Big O Notation), and hence approaches zero and  $\omega_z(t_0) t_0$ approaches a constant $K_w \neq 0$, as $t_0$ increases to a larger and larger finite value without bounds (\textbf{Statement F}). It is shown in this section that Statement F is false, if Statement 1 is true. (\textbf{Result 5.1.2.a})\\
 

We expand $F(t_0)$ in Eq.~\ref{sec:sec_x_9a1} using the identity $\cos{( \omega_{z}(t_0) (\tau- t_0))}= \cos{( \omega_{z}(t_0) t_0)} \cos{( \omega_{z}(t_0) \tau)} + \sin{( \omega_{z}(t_0) t_0)} \sin{( \omega_{z}(t_0) \tau)}$ and show that it converges to a finite positive value. We use $2 \cosh{( \sigma t_0)} = e^{\sigma t_0} + e^{-\sigma t_0}$ and rearrange terms.
\begin{align}\label{sec:sec_C_eq_5a}  
F_X(t_0)= F(t_0)  + X(t_0) = 0    \notag\\ 
F(t_0) = e^{\sigma t_0}  [ \cos{( \omega_{z}(t_0) t_0)} \int_{-t_0}^{t_0}    E_0(\tau) [ e^{-2 \sigma \tau} + 1 ] \cos{( \omega_{z}(t_0) \tau)} d\tau \notag\\
+ \sin{( \omega_{z}(t_0) t_0)} \int_{-t_0}^{t_0}    E_0(\tau) [ e^{-2 \sigma \tau} + 1 ] \sin{( \omega_{z}(t_0) \tau)} d\tau +  \int_{-t_0}^{0}  E_0(\tau) ( e^{-2 \sigma \tau} + 1) \cos{(\omega_z( t_0) \tau )} d\tau  ]
\end{align}
We write succinctly as follows, using  $\theta_{t_0} = \omega_{z}(t_0) t_0$.
\begin{align}\label{sec:sec_C_eq_5}  
F(t_0) = e^{\sigma t_0} [ \cos{(\theta_{t_0})} F_c(t_0) +  \sin{( \theta_{t_0})}  F_s(t_0) + D_1(t_0) ]  \notag\\
D_1(t_0) =    \int_{-t_0}^{0}  E_0(\tau) ( e^{-2 \sigma \tau} + 1) \cos{(\omega_z( t_0) \tau )} d\tau   \notag\\
F_c(t_0) =  \int_{-t_0}^{t_0}    E_0(\tau) [ e^{-2 \sigma \tau} + 1 ] \cos{( \omega_{z}(t_0) \tau)} d\tau  , \quad
F_s(t_0) =   \int_{-t_0}^{t_0}    E_0(\tau) [ e^{-2 \sigma \tau} + 1 ] \sin{( \omega_{z}(t_0) \tau)} d\tau
\end{align}
We consider the general case $\theta_{t_0} = \omega_{z}(t_0) t_0 \neq m \pi$ and $\theta_{t_0} \neq (2 m +1 ) \frac{\pi}{2}$, for integer $m$, where both $\cos{( \omega_{z}(t_0) t_0)}$ and $\sin{( \omega_{z}(t_0) t_0)}$ are \textbf{non-zero}.\\

We take the term $\cos{( \omega_z(t_0) \tau)}$ in the interval $|\tau| \leq t_0$ and hence $\omega_z(t_0) \tau \leq 1$ using Statement F. We write $\cos{( \omega_z(t_0) \tau)} = 1 + d_c(\tau, t_0)$ where $d_c(\tau, t_0)= - \frac{(\omega_z(t_0) \tau)^{2}}{!2} + \frac{(\omega_z(t_0) \tau)^{4}}{!4}+... $.  We take the term $\sin{( \omega_z(t_0) \tau)}$ in the interval $|\tau| \leq t_0$ and write $\sin{( \omega_z(t_0) \tau)} = \omega_z(t_0) \tau + d_s(\tau, t_0)$ where $d_s(\tau, t_0)=  - \frac{(\omega_z(t_0) \tau)^{3}}{!3} + \frac{(\omega_z(t_0) \tau)^{5}}{!5}+...  $ and write Eq.~\ref{sec:sec_C_eq_5}  as follows. (\textbf{Result 5.1.2.b})
\begin{align}\label{sec:sec_x_9a3_a}  
F_c(t_0) = F_0(t_0) + F_{d_c}(t_0), \quad F_s(t_0) = F_{s_1}(t_0) + F_{d_s}(t_0), \quad D_1(t_0) =  D_0(t_0)  + D_{d_c}(t_0) ] \notag\\
F_0(t_0) =  \int_{-t_0}^{t_0}    E_0(\tau) [ e^{-2 \sigma \tau} + 1 ]  d\tau, \quad
F_{d_c}(t_0) = \int_{-t_0}^{t_0}    E_0(\tau) [ e^{-2 \sigma \tau} + 1 ] d_c(\tau, t_0)  d\tau  \notag\\
F_{s_1}(t_0) =  \omega_z(t_0)  \int_{-t_0}^{t_0}   \tau  E_0(\tau) [ e^{-2 \sigma \tau} + 1 ]  d\tau, \quad F_{d_s}(t_0) =   \int_{-t_0}^{t_0}    E_0(\tau) [ e^{-2 \sigma \tau} + 1 ] d_s(\tau, t_0)  d\tau \notag\\
D_0(t_0)  =  \int_{-t_0}^{0}  E_0(\tau) ( e^{-2 \sigma \tau} + 1)  d\tau, \quad 
D_{d_c}(t_0)  =  \int_{-t_0}^{0}  E_0(\tau) ( e^{-2 \sigma \tau} + 1) d_c(\tau, t_0)  d\tau \notag\\
\end{align}

We see that $ F_{0}(t_0) =\int_{-t_0}^{t_0}    E_0(\tau) ( e^{-2 \sigma \tau} + 1)   d\tau = 2 \int_{0}^{t_0}    E_0(\tau) ( \cosh{(2 \sigma \tau)} + 1)   d\tau$ and $D_0(t_0)  =   \int_{-t_0}^{0}  E_0(\tau) ( e^{-2 \sigma \tau} + 1)  d\tau = \int_{0}^{t_0}  E_0(\tau) ( e^{2 \sigma \tau} + 1)  d\tau   $ converge to a \textbf{non-zero} positive finite value, as $t_0$ increases, as shown in Section~\ref{sec:Section_C_6_2s0a}. (\textbf{Result 5.1.2.c}) We write Eq.~\ref{sec:sec_C_eq_5} as follows using Eq.~\ref{sec:sec_x_9a3_a}.
\begin{align}\label{sec:sec_C_eq_6}  
F(t_0) = e^{\sigma t_0} [ F_A(t_0) + F_B(t_0)] , \quad 
F_A(t_0) =  \cos{(\theta_{t_0})} F_0(t_0)  + D_0(t_0) \notag\\
 F_B(t_0) = \cos{(\theta_{t_0})} F_{d_c}(t_0) +  \sin{( \theta_{t_0})}  ( F_{s_1}(t_0) + F_{d_s}(t_0) ) + D_{d_c}(t_0) 
\end{align}
We will show that, \textbf{compared} to $F_0(t_0)$ and $D_0(t_0)$ in $F_A(t_0)$, rest of the terms in $F_B(t_0)$ becomes smaller and smaller approaching zero, as $t_0$ increases to a larger and larger finite value without bounds, using Statement F and Result 5.1.2.b. We consider 4 cases below.\\

\textbf{Case 2A:} We consider the case $F_A(t_0) =  \cos{(\theta_{t_0})} F_0(t_0)  + D_0(t_0)$ in Eq.~\ref{sec:sec_C_eq_6}, approaches a \textbf{non-zero} constant $F_{A_{min}}$, as $t_0$ increases to a larger and larger finite value without bounds, using Statement F and Result 5.1.2.c and $\theta_{t_0}  = \omega_{z}(t_0) t_0 \neq m \pi$ and $\theta_{t_0} \neq (2 m +1 ) \frac{\pi}{2}$, for integer $m$, where both $\cos{( \theta_{t_0})}$ and $\sin{( \theta_{t_0})}$ are \textbf{non-zero}. (\textbf{Result 5.1.2.d})\\

As we increase $t_0$ to $t_0'$, we see that $F_{0}(t_0') > F_{0}(t_0) > 0$ and $D_{0}(t_0') > D_{0}(t_0) > 0$ using Result 5.1.2.c (Details in Section~\ref{sec:Section_C_6_2s0a}). \\ 

We see that the first term in $F_B(t_0)$  in Eq.~\ref{sec:sec_C_eq_6}, given by $F_{d_c}(t_0)= \int_{-t_0}^{t_0}    E_0(\tau)  ( e^{-2 \sigma \tau} + 1 ) \\ d_c(\tau, t_0)  d\tau $ in Eq.~\ref{sec:sec_x_9a3_a}, can be written as $F_{d_c}(t_0) = - \frac{(\omega_z(t_0)^{2}}{!2} \int_{-t_0}^{t_0}   (\tau)^{2}  E_0(\tau) ( e^{-2 \sigma \tau} + 1 )  d_c^{'}(\tau, t_0) d\tau $, using $d_c(\tau, t_0) = - \frac{(\omega_z(t_0) \tau)^{2}}{!2} + \frac{(\omega_z(t_0) \tau)^{4}}{!4}+...  = -\frac{(\omega_z(t_0) \tau)^{2}}{!2} d_c^{'}(\tau, t_0) d\tau $, where $d_c^{'}(\tau, t_0) = [ 1 - (!2) \frac{(\omega_z(t_0) \tau)^{2}}{!4}+... ]  $, in the interval $|\tau| \leq t_0$, using Result 5.1.2.b and the terms in $d_c^{'}(\tau, t_0)$ converge, given $\omega_z(t_0)= O[\frac{1}{t_0}]$ and $\omega_z(t_0) \tau \leq 1$.\\

 We see that the integral in $F_{d_c}(t_0) = - \frac{(\omega_z(t_0)^{2}}{!2}  I(t_0)$ given by $I(t_0)=\int_{-t_0}^{t_0}   (\tau)^{2}  E_0(\tau) ( e^{-2 \sigma \tau} + 1 )  d_c^{'}(\tau, t_0) d\tau $ converges to a finite value, as $t_0$ increases to a larger and larger finite value without bounds (Details in Section~\ref{sec:Section_C_6_2s0a}) and hence  $F_{d_c}(t_0) = - \frac{(\omega_z(t_0)^{2}}{!2}  I(t_0)$ becomes smaller and smaller approaching zero, due to the term $(\omega_z(t_0) )^{2}$ where $\omega_z(t_0)= O[\frac{1}{t_0}]$, as $t_0$ increases. Hence $|F_{d_c}(t_0) \cos{(\theta_{t_0})} | << |F_{A_{min}}|$ (\textbf{Result 5.1.2.e})\\

We see that the third term in $F_B(t_0)$  in Eq.~\ref{sec:sec_C_eq_6}, is given by $D_{d_c}(t_0)= \int_{-t_0}^{0}    E_0(\tau)  ( e^{-2 \sigma \tau} + 1 ) \\ d_c(\tau, t_0)  d\tau $ in Eq.~\ref{sec:sec_x_9a3_a}, and $|D_{d_c}(t_0)| << |F_{A_{min}}|$ using arguments in Result 5.1.2.e. (\textbf{Result 5.1.2.f})\\

We see that the second term in in $F_B(t_0)$  in Eq.~\ref{sec:sec_C_eq_6}, given by  $F_s(t_0)=F_{s_1}(t_0)+F_{d_s}(t_0) \\ =   \int_{-t_0}^{t_0}    E_0(\tau) ( e^{-2 \sigma \tau} + 1 ) \sin{( \omega_{z}(t_0) \tau)} d\tau  = \omega_z(t_0)  \int_{-t_0}^{t_0} \tau   E_0(\tau) ( e^{-2 \sigma \tau} + 1 )  d_s^{'}(\tau, t_0)  d\tau $, where $\sin{( \omega_z(t_0) \tau)}  =  \omega_z(t_0) \tau -\frac{(\omega_z(t_0) \tau)^{3}}{!3} + \frac{(\omega_z(t_0) \tau)^{5}}{!5}+... = (\omega_z(t_0) \tau) d_s^{'}(\tau, t_0) $ where $d_s^{'}(\tau, t_0) =[ 1 - \frac{(\omega_z(t_0) \tau)^{2}}{!3}+... ] $, in the interval $|\tau| \leq t_0$, using Result 5.1.2.b and the terms in $d_s^{'}(\tau, t_0)$ converge, given $\omega_z(t_0)= O[\frac{1}{t_0}]$  and $\omega_z(t_0) \tau \leq 1$. \\

We see that the integral in $F_{s}(t_0) = \omega_z(t_0)  J(t_0)$ given by $J(t_0)=\int_{-t_0}^{t_0}   \tau E_0(\tau) (e^{-2 \sigma \tau}+1)  d_s^{'}(\tau, t_0) d\tau $ converges to a finite value, as $t_0$ increases to a larger and larger finite value without bounds(Details in Section~\ref{sec:Section_C_6_2s0a}) and hence  $F_{s}(t_0) = \omega_z(t_0)  J(t_0)$ becomes smaller and smaller approaching zero, due to the term $\omega_z(t_0)$ where $\omega_z(t_0)= O[\frac{1}{t_0}]$, as $t_0$ increases. Hence $|F_{s}(t_0)  \sin{(\theta_{t_0})}|$ is $<< |F_{A_{min}}|$ (\textbf{Result 5.1.2.g})\\

We use Result 5.1.2.d, Result 5.1.2.e, Result 5.1.2.f and Result 5.1.2.g and see that $ | F_B(t_0)| \leq |\cos{(\theta_{t_0})} F_{d_c}(t_0)| +  |\sin{( \theta_{t_0})}  ( F_{s_1}(t_0) + F_{d_s}(t_0) )| + |D_{d_c}(t_0)| <<  |F_{A_{min}}|  $ and hence $|F(t_0)|= e^{\sigma t_0} | F_A(t_0) + F_B(t_0)|$ in Eq.~\ref{sec:sec_C_eq_6} approaches $ e^{\sigma t_0} |F_{A_{min}}| \neq 0$, as $t_0$ increases to a larger and larger finite value without bounds.\\

We see that $  |X(t_0)| << 1$ in Section~\ref{sec:Section_C_6_2s}  and we can choose $ |X(t_0)| << |F(t_0)|$ by increasing $t_0$ to $t_0' >> t_0$ and hence $F(t_0)$ does not cancel $ X(t_0)$ in Eq.~\ref{sec:sec_C_eq_5a} (\textbf{Result 5.1.2.h}). Hence we see that $F_X(t_0)= F(t_0)  + X(t_0)$ \textbf{cannot} equal zero in Eq.~\ref{sec:sec_C_eq_5a}, using Result 5.1.2.h   and hence Statement F is \textbf{false} for Case 2A.\\

Hence Statement F is \textbf{false} for Case 2A and  $\omega_z(t_0) t_0$ \textbf{does not} approach a constant, as  $t_0$ increases to a larger and larger finite value without bounds and hence  $\omega_z(t_0) t_0$ increases with $t_0$, as $t_0$ increases. \\

\textbf{Case 2B:} We consider the case $F_A(t_0) =  \cos{(\theta_{t_0})} F_0(t_0)  + D_0(t_0)$ in Eq.~\ref{sec:sec_C_eq_6}, becomes smaller and smaller approaching \textbf{zero}, as $t_0$ increases to a larger and larger finite value without bounds, using Statement F and Result 5.1.2.c and $\theta_{t_0} = \omega_z(t_0) t_0 \neq m \pi$ and $\theta_{t_0} \neq (2 m +1 ) \frac{\pi}{2}$, for integer $m$, where both $\cos{( \theta_{t_0})}$ and $\sin{( \theta_{t_0})}$ are \textbf{non-zero}. (\textbf{Result 5.1.2.i})\\

We take $ F_B(t_0) = \cos{(\theta_{t_0})} F_{d_c}(t_0) +  \sin{( \theta_{t_0})}  ( F_{s_1}(t_0) + F_{d_s}(t_0) ) + D_{d_c}(t_0) $ in Eq.~\ref{sec:sec_C_eq_6} and  we can show that, compared to the second term $F_{s_1}(t_0) = \omega_z(t_0)  \int_{-t_0}^{t_0}   \tau  E_0(\tau) [ e^{-2 \sigma \tau} + 1 ]  d\tau $, the other terms have \textbf{faster }fall-off rate as $t_0$ increases, due to terms $(\omega_z(t_0)^{r}$ where $r > 1$ and $\omega_z(t_0)= O[\frac{1}{t_0}]$ and hence \textbf{negligible}, as shown below.\\

In Eq.~\ref{sec:sec_C_eq_5},  Eq.~\ref{sec:sec_x_9a3_a} and Eq.~\ref{sec:sec_C_eq_6}, we consider the integral $F_{s}(t_0) = \int_{-t_0}^{t_0}    E_0(\tau) ( e^{-2 \sigma \tau} + 1)  \sin{( \omega_z(t_0) \tau)} d\tau = F_{s_1}(t_0) +  F_{d_s}(t_0)= \omega_z(t_0)  \int_{-t_0}^{t_0}   \tau  E_0(\tau) ( e^{-2 \sigma \tau} + 1)   d\tau +  \int_{-t_0}^{t_0}    E_0(\tau) ( e^{-2 \sigma \tau} + 1) d_s(\tau, t_0) d\tau $, using $d_s(\tau, t_0)= - \frac{(\omega_z(t_0) \tau)^{3}}{!3} + \frac{(\omega_z(t_0) \tau)^{5}}{!5}+...$ where $ d_s(\tau, t_0) =  - \frac{(\omega_z(t_0) \tau)^{3}}{!3}   d_s^{'}(\tau, t_0)$ and $d_s{'}(\tau, t_0)= [ 1 - (!3)\frac{(\omega_z(t_0) \tau)^{2}}{!5}+... ]  $, in the interval $|\tau| \leq t_0$, using Result 5.1.2.b and the terms in $d_s^{'}(\tau, t_0)$ converge, given $\omega_z(t_0)= O[\frac{1}{t_0}]$ and $\omega_z(t_0) \tau \leq 1$. We see that $F_{s_1}(t_0)= \omega_z(t_0)  K(t_0)$ with $ K(t_0) = \int_{-t_0}^{t_0}   \tau  E_0(\tau) ( e^{-2 \sigma \tau} + 1)   d\tau$ and $K(t_0)$ converges to a \textbf{non-zero} and finite value, as $t_0$ increases to a larger and larger finite value without bounds, as shown in Section~\ref{sec:Section_C_6_2s0a}. (\textbf{Result 5.1.2.j}) \\

The second integral in $F_{s}(t_0)$ in above para, given by $F_{d_s}(t_0)=\int_{-t_0}^{t_0}    E_0(\tau) ( e^{-2 \sigma \tau} + 1) d_s(\tau, t_0) d\tau \\ = - \frac{(\omega_z(t_0)^{3} }{!3}   \int_{-t_0}^{t_0} \tau^{3}   E_0(\tau) ( e^{-2 \sigma \tau} + 1)  d_s^{'}(\tau, t_0) d\tau $ with the \textbf{term} $(\omega_z(t_0)^{3}$, becomes smaller and smaller approaching zero \textbf{faster than} the first integral $F_{s_1}(t_0)=(\omega_z(t_0) \int_{-t_0}^{t_0}   \tau  E_0(\tau) ( e^{-2 \sigma \tau} + 1)  d\tau$ with the \textbf{term} $(\omega_z(t_0)$, using $\omega_z(t_0)= O[\frac{1}{t_0}]$, as $t_0$ increases to a larger and larger finite value without bounds.\\

Hence $|F_{d_s}(t_0)| << |F_{s_1}(t_0)|$ and $ e^{\sigma t_0}|F_{d_s}(t_0)| << e^{\sigma t_0} |F_{s_1}(t_0)|$ and $ e^{\sigma t_0} F_{s_1}(t_0)$ is the \textbf{dominant} term in Case 2B, which increases exponentially as $t_0$ increases, due to the factor $\frac{e^{\sigma t_0}}{t_0}$. It is noted that the multiplication factor $e^{\sigma t_0}$ is taken from the equation for $F(t_0)$ in Eq.~\ref{sec:sec_C_eq_6} (\textbf{Result 5.1.2.k}) \\

We see that the terms $F_{d_c}(t_0)$ and $D_{d_c}(t_0)$ in Eq.~\ref{sec:sec_C_eq_6} have terms of the form $(\omega_z(t_0)^{r}$ where $r \geq 2$ and hence we can use arguments in above two paras and Result 5.1.2.j and Result 5.1.2.k, to show that these terms have \textbf{faster }fall-off rates of $O[\frac{1}{t_0^{r}}]$ and are \textbf{negligible}, compared to the term $F_{s_1}(t_0)$ which has a fall-off rate of $O[\frac{1}{t_0}]$, as $t_0$ increases to a larger and larger finite value without bounds. (\textbf{Result 5.1.2.m}) \\

We use Result 5.1.2.j, Result 5.1.2.k and Result 5.1.2.m and see that, as $t_0$ increases to a larger and larger finite value without bounds, $ | F_B(t_0)| $ is \textbf{nearly equal} to $  |\sin{( \theta_{t_0})}   F_{s_1}(t_0) | > 0 $ and hence $|F(t_0)|= e^{\sigma t_0} | F_A(t_0) + F_B(t_0)| $ is nearly equal to $ e^{\sigma t_0} |F_B(t_0)| $ in Eq.~\ref{sec:sec_C_eq_6}, which approaches $ e^{\sigma t_0} | \sin{( \theta_{t_0})}  F_{s_1}(t_0)| $ which is of the order $O[\frac{e^{\sigma t_0}}{t_0}]$, using Result 5.1.2.k.\\

We see that $  |X(t_0)| << 1$ in Section~\ref{sec:Section_C_6_2s}  and we can choose $ |X(t_0)| << |F(t_0)|$ by increasing $t_0$ to $t_0' >> t_0$ and hence $F(t_0)$ does not cancel $ X(t_0)$ in Eq.~\ref{sec:sec_C_eq_5a}. Hence we see that $F_X(t_0)= F(t_0)  + X(t_0)$ \textbf{cannot} equal zero in Eq.~\ref{sec:sec_C_eq_5a}   and hence Statement F is \textbf{false} for Case 2B. (\textbf{Result 5.1.2.n})\\

Hence Statement F is \textbf{false} for Case 2B and  $\omega_z(t_0) t_0$ \textbf{does not} approach a constant, as  $t_0$ increases to a larger and larger finite value without bounds and hence  $\omega_z(t_0) t_0$ increases with $t_0$, as $t_0$ increases, using Result 5.1.2.n. \\

\textbf{Case 2C:} We consider the specific case $\theta_{t_0} = \omega_z(t_0) t_0 =  m \pi$ in Eq.~\ref{sec:sec_C_eq_6} and we get $\sin{( \theta_{t_0})}=0$ and $\cos{( \theta_{t_0})} = (-1)^{m}$ and hence $F(t_0) = e^{\sigma t_0}  [ (-1)^{m} F_0(t_0)  + D_0(t_0) + (-1)^{m} F_{d_c}(t_0)  + D_{d_c}(t_0)]$ in Eq.~\ref{sec:sec_C_eq_6}. ( \textbf{Result 5.1.2.o}) We will show that $F(t_0)$ converges to a non-zero finite value, as $t_0$ increases. \\

Using Result 5.1.2.c, we see that the term $A_1(t_0)=(-1)^{m} F_0(t_0)  + D_0(t_0) =  \int_{0}^{t_0}    E_0(\tau) ( (-1)^{m} 2 ( \cosh{(2 \sigma \tau)} + 1) + e^{2 \sigma \tau} + 1)   d\tau$ and $A_1(t_0)=  \int_{0}^{t_0}    E_0(\tau)   ( 2 \cosh{(2 \sigma \tau)} + 2  + e^{2 \sigma \tau} + 1)   d\tau  > 0$ for even $m$ and $A_1(t_0)=  \int_{0}^{t_0}    E_0(\tau)   ( -2 \cosh{(2 \sigma \tau)} - 2  + e^{2 \sigma \tau} + 1)   d\tau  = \int_{0}^{t_0}    E_0(\tau)   ( -e^{2 \sigma \tau} - e^{-2 \sigma \tau}  - 2  + e^{2 \sigma \tau} + 1)   d\tau = - \int_{0}^{t_0}    E_0(\tau)   (  e^{-2 \sigma \tau}   + 1)   d\tau  < 0$ for odd $m$ and hence $A_1(t_0)$ converges to a \textbf{non-zero} finite value for integer $m$, as $t_0$ increases to a larger and larger finite value without bounds, as shown in Section~\ref{sec:Section_C_6_2s0a}. ( \textbf{Result 5.1.2.p})\\

The term $A_2(t_0)=(-1)^{m} F_{d_c}(t_0)  + D_{d_c}(t_0)$ in Result 5.1.2.o has fall-off rates of $O[\frac{1}{t_0^{2}}]$, as $t_0$ increases to a larger and larger finite value without bounds, and \textbf{negligible} compared to $A_1(t_0)$, using arguments used in Result 5.1.2.e and Result 5.1.2.f. ( \textbf{Result 5.1.2.q})\\

Hence $F(t_0)$ nearly equals  $e^{\sigma t_0} A_1(t_0)$, which converges to a \textbf{non-zero} finite value, which grows exponentially, using Result 5.1.2.o, Result 5.1.2.p and 5.1.2.q.\\

We see that $  |X(t_0)| << 1$ in Section~\ref{sec:Section_C_6_2s}  and we can choose $ |X(t_0)| << |F(t_0)|$ by increasing $t_0$ to $t_0' >> t_0$ and hence $F(t_0)$ does not cancel $ X(t_0)$ in Eq.~\ref{sec:sec_C_eq_5a}. Hence we see that $F_X(t_0)= F(t_0)  + X(t_0)$ \textbf{cannot} equal zero in Eq.~\ref{sec:sec_C_eq_5a} and hence Statement F is \textbf{false}.\\

\textbf{Case 2D:} We consider the specific case $\theta_{t_0} = \omega_z(t_0) t_0 =  (2 m + 1 ) \frac{\pi}{2}$ in Eq.~\ref{sec:sec_C_eq_6} and we get $\sin{( \theta_{t_0})}= (-1)^{m} $ and $\cos{( \theta_{t_0})} = 0$ and hence $F(t_0) = e^{\sigma t_0}  [  D_0(t_0) + (-1)^{m} ( F_{s_1}(t_0) + F_{d_s}(t_0) )   + D_{d_c}(t_0)]$ in Eq.~\ref{sec:sec_C_eq_6}. ( \textbf{Result 5.1.2.r}) We will show that $F(t_0)$ converges to a non-zero finite value, as $t_0$ increases. \\

Using Result 5.1.2.c, we see that the term $ D_0(t_0) =  \int_{0}^{t_0}    E_0(\tau)  ( e^{2 \sigma \tau} + 1)   d\tau > 0$ and converges to a \textbf{non-zero} finite value for integer $m$, as $t_0$ increases to a larger and larger finite value without bounds, as shown in Section~\ref{sec:Section_C_6_2s0a}. ( \textbf{Result 5.1.2.s})\\

The term $A_3(t_0)=(-1)^{m} ( F_{s_1}(t_0) + F_{d_s}(t_0) )   + D_{d_c}(t_0)$ in Result 5.1.2.r, has fall-off rates of $O[\frac{1}{t_0}]$, as $t_0$ increases to a larger and larger finite value without bounds, and \textbf{negligible} compared to $D_0(t_0)$, using arguments used in Result 5.1.2.f and Result 5.1.2.g. ( \textbf{Result 5.1.2.t})\\

Hence $F(t_0)$ nearly equals  $e^{\sigma t_0} D_0(t_0)$, which converges to a \textbf{non-zero} finite value, which grows exponentially, using Result 5.1.2.r, Result 5.1.2.s and 5.1.2.t.\\

We see that $  |X(t_0)| << 1$ in Section~\ref{sec:Section_C_6_2s}  and we can choose $ |X(t_0)| << |F(t_0)|$ by increasing $t_0$ to $t_0' >> t_0$ and hence $F(t_0)$ does not cancel $ X(t_0)$ in Eq.~\ref{sec:sec_C_eq_5a}. Hence we see that $F_X(t_0)= F(t_0)  + X(t_0)$ \textbf{cannot} equal zero in Eq.~\ref{sec:sec_C_eq_5a} and hence Statement F is \textbf{false}.

.\\ \tocless\subsection{\label{sec:Section_C_6_2s} \textbf{ $|X(t_0| << 1$  } \protect \\ \lowercase{} }

We compute the \textbf{maximum} value of the absolute value of  $X(t_0)$ in Eq.~\ref{sec:sec_x_9a1}, for sufficiently large $t_0$ and $0 < \sigma < \frac{1}{2}$. 
\begin{align}  \label{sec:sec_az3_1_3_eq_1_5_cs}   
X(t_0)=   e^{\sigma t_0}    \int_{-\infty}^{-t_0}    E_0(\tau) [ e^{-2 \sigma \tau} + 1 ] \cos{( \omega_{z}(t_0) (\tau- t_0)} d\tau   \notag\\
+ e^{-\sigma t_0}  \int_{-\infty}^{- t_0}    E_0(\tau) e^{-2 \sigma \tau} \cos{( \omega_{z}(t_0) (\tau + t_0)} d\tau 
+  e^{ -\sigma t_0}  \int_{-\infty}^{- t_0}   E_0(\tau) \cos{( \omega_{z}(t_0) (\tau + t_0)}      d\tau \notag\\
+ e^{\sigma t_0}  \int_{-\infty}^{-t_0}  E_0(\tau) ( e^{-2 \sigma \tau} + 1) \cos{(\omega_z( t_0) \tau )} d\tau 
+ e^{-\sigma t_0}  \int_{-\infty}^{0}  E_0(\tau) ( e^{-2 \sigma \tau} + 1) \cos{(\omega_z( t_0) \tau )} d\tau
 \end{align}
 We use $|f(t_0)| \leq \int | f_i(\tau, t_0) | d\tau$ where $f(t_0) =  \int  f_i(\tau, t_0)  d\tau$ (\href{https://web.archive.org/web/20240218091537/http://www.math.ualberta.ca/~isaac/math311/s14/abs_value.pdf}{link} ) and use $|\cos{( \omega_{z}(t_0) (\tau- t_0)}| \leq 1$, $|\cos{( \omega_{z}(t_0) (\tau + t_0)}| \leq 1$ and $|\cos{( \omega_{z}(t_0) \tau )}| \leq 1$ and $E_0(\tau) > 0$ for real $\tau$. We combine first and fourth integrals and combine second and third integrals in Eq.~\ref{sec:sec_az3_1_3_eq_1_5_cs}  and rearrange as follows.
 \begin{align}  \label{sec:sec_az3_1_3_eq_1_5_cs0}   
| X(t_0) | \leq   2 e^{\sigma t_0}    \int_{-\infty}^{-t_0}    E_0(\tau) [ e^{-2 \sigma \tau} + 1 ] d\tau + e^{-\sigma t_0}  \int_{-\infty}^{- t_0}    E_0(\tau) [ e^{-2 \sigma \tau} + 1 ] d\tau    \notag\\
+ e^{-\sigma t_0}  \int_{-\infty}^{0}  E_0(\tau) ( e^{-2 \sigma \tau} + 1)  d\tau
 \end{align}
The second and third integrals in Eq.~\ref{sec:sec_az3_1_3_eq_1_5_cs0} converge given that the integrands have asymptotic exponential fall-off rate using ~\ref{sec:appendix_C_5} and multiplied by the term $ e^{-\sigma t_0} $, these integrals goes to zero, as $t_0$ increases to a larger and larger finite value without bounds (\textbf{Result 5.2.a}).\\

We consider the first integral $I(t_0)= 2 e^{\sigma t_0}  \int_{-\infty}^{-t_0}  E_0(\tau) ( e^{-2 \sigma \tau} + 1 ) d\tau $ and substitute $\tau=-\tau$ and get $ I(t_0) = 2 e^{\sigma t_0}   \int_{t_0}^{\infty}  E_0(\tau) ( e^{2 \sigma \tau} + 1)  d\tau  $. (\textbf{Result 5.2.a.1}) \\

 We see that $E_0(\tau)  e^{ 2 \sigma \tau}$ has an asymptotic minimum fall-off rate of $o[e^{-0.5 |\tau| } ] $ where $0 < \sigma < \frac{1}{2}$ (Details in ~\ref{sec:appendix_C_5}). Hence we write $|I(t_0)| \leq 2 e^{\sigma t_0} K_I \int_{t_0}^{\infty}  e^{-0.5 \tau}  d\tau = 2  e^{\sigma t_0}  K_I  [ \frac{e^{- 0.5 \tau}}{-0.5}  ]_{t_0}^{\infty} =  4 e^{\sigma t_0} K_I e^{- 0.5 t_0} = K(t_0)$.  (\textbf{Result 5.2.b}) The constant $K_I$ is chosen as below.\\
 
 It is noted that $K_I =\frac{E_0(t_0) ( e^{2 \sigma t_0} + 1)}{ e^{-0.5 t_0} }$ is chosen such that the integrand $ e^{-0.5 \tau}  $ in above para, when evaluated at the lower limit $t_0$ and multiplied by the term $K_I$ in the equation for $K(t_0)$, gives the \textbf{same} value $ K_I e^{- 0.5 t_0}$, as the integrand $  E_0(\tau) ( e^{2 \sigma \tau} + 1)  $ in $I(t_0)$ in Result 5.2.a.1, when evaluated at the lower limit $t_0$ (\textbf{Note 5.2}). \\
 
 For $0 < \sigma < \frac{1}{2}$, $|I(t_0)|$ has an asymptotic \textbf{minimum} fall-off rate of $o[e^{- (0.5-\sigma) | t_0| } ] $ using Result 5.2.b and hence $| X(t_0) | $ goes to zero, using Result 5.2.a, Result 5.2.a.1 and Result 5.2.b, as $t_0$ increases to a larger and larger finite value without bounds. 
 
.\\ \tocless\subsection{\label{sec:Section_C_6_2s0a} \textbf{ Derivation of Result A  } \protect \\ \lowercase{} }

In this section, it is shown that $I(t_0, r) = \int_{-t_0}^{t_0} ( \tau)^{r}   E_0(\tau) ( e^{-2 \sigma \tau} + 1)  d\tau$ converges to a non-zero and finite value, as $t_0$ increases to a larger and larger finite value without bounds, for non-negative integer $r$ and $-\frac{1}{2} < \sigma < \frac{1}{2}$. Hence $(\omega_z(t_0))^{r}  |I(t_0, r)|$ goes to zero, as $t_0$ increases to a larger and larger finite value without bounds, for $\omega_z(t_0)= o[\frac{1}{t_0}]$ and $\omega_z(t_0)= O[\frac{1}{t_0}]$ and $r \geq 1$.\\

We see that $E_0(\tau) (e^{-2 \sigma \tau}+1)$ is an absolutely integrable function, for $0 \leq |\sigma| < \frac{1}{2}$ given that it has exponential fall-off rates as $|\tau| \to \infty$. (Details in ~\ref{sec:appendix_C_5} and ~\ref{sec:appendix_C_5z}).\\

Hence $( \tau)^{r}   E_0(\tau) (e^{-2 \sigma \tau} +1)$ also has exponential fall-off rates as $|\tau| \to \infty$, and is an absolutely integrable function. Hence $I(t_0, r) = \int_{-t_0}^{t_0} ( \tau)^{r}   E_0(\tau) ( e^{-2 \sigma \tau}   + 1 ) d\tau < \infty$.\\

We see that $I(t_0, r) = \int_{-t_0}^{t_0} ( \tau)^{r}   E_0(\tau) ( e^{-2 \sigma \tau} + 1)  d\tau = \int_{0}^{t_0} ( \tau)^{r}   E_0(\tau)( e^{-2 \sigma \tau}+1)  d\tau + \int_{-t_0}^{0} ( \tau)^{r}   E_0(\tau) (e^{-2 \sigma \tau}+1)  d\tau $. We substitute $\tau=-\tau$ in second integral and use $E_0(-\tau)=E_0(\tau)$ and get $\int_{-t_0}^{0} ( \tau)^{r}   E_0(\tau)( e^{-2 \sigma \tau} +1) d\tau = (-1)^{r} \int_{0}^{t_0} ( \tau)^{r}   E_0(\tau) (e^{2 \sigma \tau} +1) d\tau $. Hence $I(t_0, r) =\int_{0}^{t_0} ( \tau)^{r}   E_0(\tau) (e^{-2 \sigma \tau}+1)  d\tau + (-1)^{r} \int_{0}^{t_0} ( \tau)^{r}   E_0(\tau) (e^{2 \sigma \tau} +1) d\tau$. \\

Hence $I(t_0, r) = 2 \int_{0}^{t_0} ( \tau)^{r}   E_0(\tau) (\cosh{(2 \sigma \tau)} +1) d\tau  $ for \textbf{even} $r$ and $I(t_0, r) = -2 \int_{0}^{t_0} ( \tau)^{r}   E_0(\tau) \sinh{(2 \sigma \tau)}  d\tau  $ for \textbf{odd} $r$ . Hence $I(t_0, r)$ is\textbf{ non-zero} and finite, given that $E_0(\tau)>0$ and the integrands in $I(t_0, r)$ are positive in the interval $0 \leq \tau \leq t_0$, for non-negative integer $r$. (\textbf{Result 5.3.1}).\\

As we increase $t_0$ to $t_0' >> t_0$, we see that $|I(t_0', r)| > |I(t_0, r)| > 0$, using Result 5.3.1, given that $ ( \tau)^{r}  E_0(\tau)  (\cosh{(2 \sigma \tau)} +1 ) >0$ and  $ ( \tau)^{r}  E_0(\tau)  \sinh{(2 \sigma \tau)} >0$ for $0 \leq \tau \leq t_0$ (\textbf{Result 5.3.2}).\\

We can show that $ (\omega_z(t_0))^{r} |I(t_0, r)|$ goes to zero, as $t_0$ increases to $t_0' >> t_0$, for $\omega_z(t_0)= o[\frac{1}{t_0}]$ and $\omega_z(t_0)= O[\frac{1}{t_0}]$ and $r \geq 1$. For \textbf{even} $r$, using Result 5.3.1, we see that $I(t_0', r) =  2 \int_{0}^{t_0'} ( \tau)^{r}   E_0(\tau) (\cosh{(2 \sigma \tau)} +1)  d\tau =   \int_{0}^{t_0} ( \tau)^{r}   E_0(\tau) ( e^{2 \sigma \tau} + e^{-2 \sigma \tau} +2)  d\tau +   \int_{t_0}^{t_0'} (\tau)^{r}   E_0(\tau) ( e^{2 \sigma \tau} + e^{-2 \sigma \tau} +2)  d\tau  = I(t_0, r) + J(t_0, t_0', r) $. (\textbf{Result 5.3.3}) \\

Given that the integrand in $J(t_0, t_0', r)$ has minimum asymptotic exponential fall-off rate of $o[e^{-0.5 \tau}]$ using ~\ref{sec:appendix_C_5}, we get $J(t_0, t_0', r) \leq  K_I \int_{t_0}^{t_0'} e^{-0.5 \tau} d\tau $, for sufficiently large $t_0' >> t_0$, where $K_I = \frac{(t_0)^{r}   E_0(t_0) ( e^{2 \sigma t_0} + e^{-2 \sigma t_0} +2) }{e^{-0.5 t_0}}$ is a constant chosen as in Note 5.2 in previous subsection. (\textbf{Result 5.3.4})\\

We get $J(t_0, t_0', r) \leq  K_I  [ \frac{e^{-0.5 \tau}}{-0.5} ]_{t_0}^{t_0'} = -2 K_I ( e^{-0.5 t_0'} - e^{-0.5 t_0} ) = 2 K_I e^{-0.5 t_0} (1 - e^{-0.5 (t_0' - t_0) }  ) $. As  we increase $t_0' >> t_0$, we see that $J(t_0, t_0', r) \leq 2 K_I e^{-0.5 t_0}$ becomes smaller and smaller approaching zero and $I(t_0', r) =  I(t_0, r) + J(t_0, t_0', r) $ \textbf{approaches} $I(t_0, r)$ in Result 5.3.1, which converges. Hence $I(t_0', r)$ does not continue to increase with increasing $t_0$. We get similar results for $I(t_0', r)$ for \textbf{odd} $r$. Hence $(\omega_z(t_0))^{r}  |I(t_0, r)|$ goes to zero, as $t_0$ increases to a larger and larger finite value without bounds, for $\omega_z(t_0)= o[\frac{1}{t_0}]$ and $\omega_z(t_0)= O[\frac{1}{t_0}]$ and $r \geq 1$. (\textbf{Result 5.3.5})


.\\ \tocless\section{\label{sec:Section_A_1_6} \textbf{Strictly decreasing $E_0(t)$ for $t > 0$ } \protect\\  \lowercase{} }

Let us consider $E_{0}(t)  =  \Phi(t) =     \sum_{n=1}^{\infty}  [ 4 \pi^{2} n^{4} e^{4t}    - 6 \pi n^{2}   e^{2t} ]  e^{- \pi n^{2} e^{2t}} e^{\frac{t}{2}} $ in Eq.~\ref{sec_intro_eq_1}, where $t$ is real, whose Fourier Transform is given by the entire function $E_{0\omega}(\omega) = \xi(\frac{1}{2}+ i \omega)$. It is  known that $\Phi(t)$ is positive for $|t| < \infty$ and its first derivative is negative for $t>0$ and hence $\Phi(t)$  is a \textbf{strictly decreasing} function for $t>0$. (\href{https://web.archive.org/web/20240218091913/https://www.ams.org/notices/200303/fea-conrey-web.pdf#page=5}{Conrey}, \href{https://www.ocf.berkeley.edu/~araman/files/math_z/decr/E0_t_decreasing_April_16_2022_Akhila_Raman.pdf}{link 2},  \href{https://www.ocf.berkeley.edu/~araman/files/math_z/decr/zeta_paper_Raman_Sep_25_2022.pdf#page=12}{Pages 12-17}, {\citep{TB}}). This is shown  below. We take the term $2 \pi n^{2} $ out of the brackets.
\begin{align}\label{sec:sec_a_1_eq_1}   
E_{0}(t) = \Phi(t) =     \sum_{n=1}^{\infty}  [ 4 \pi^{2} n^{4} e^{4t}    - 6 \pi n^{2}   e^{2t} ]  e^{- \pi n^{2} e^{2t}} e^{\frac{t}{2}}  = \sum_{n=1}^{\infty} 2 \pi n^{2} e^{- \pi n^{2} e^{2t}} e^{\frac{t}{2}} [ 2 \pi n^{2} e^{4t}    - 3 e^{2t}  ] 
\end{align}
We show that $X(t)= \frac{E_0(t)}{2}$ is a \textbf{strictly decreasing} function for $t > 0$ as follows.\\

$\bullet$ In Section~\ref{sec:Section_A_1_6_0}, it is shown that the first derivative of $X(t)$, given by $\frac{d X(t)}{dt} <0$  for $t > t_z $ where $t_z = \frac{1}{2} \log{\frac{y_z}{\pi}}$ and $y_z= 3.16$.\\

$\bullet$  In Section~\ref{sec:Section_A_1_6_1}, it is shown that, $\frac{d X(t)}{dt} < 0$ for $0 < t \leq t_z$.\\

Hence $\frac{d X(t)}{dt} < 0$ for $t >0$. Given that $2 X(t)=E_0(t)=E_0(-t) > 0$ for $|t| < \infty$ (Details in ~\ref{sec:appendix_C_5a} and ~\ref{sec:appendix_C_7}), we see that $X(t)$  is strictly decreasing for  $t>0$ and $E_0(t)=2 X(t)$  is \textbf{strictly decreasing} for $t>0$.

.\\ \tocless\subsection{\label{sec:Section_A_1_6_0} \textbf{ $\frac{d X(t)}{dt} < 0$ for $t > t_z $} \protect\\  \lowercase{} }

We consider $X(t)= \frac{E_0(t)}{2} =\sum_{n=1}^{\infty} \pi n^{2} e^{- \pi n^{2} e^{2t}} e^{\frac{t}{2}} [ 2 \pi n^{2} e^{4t}    - 3 e^{2t}  ] $ in Eq.~\ref{sec:sec_a_1_eq_1} and take the first derivative of $X(t)$. We note that $E_0(t)$ and $X(t)$ are analytic functions for real $t$ and infinitely differentiable in that interval. We compute $\frac{d X(t)}{dt}$ below and take the term $e^{2t}$ out, in the last line below.
\begin{align}\label{sec:sec_a_1_eq_1_a}   
 \frac{d X(t)}{dt} = \sum_{n=1}^{\infty} \pi n^{2} e^{- \pi n^{2} e^{2t}} e^{\frac{t}{2}} [ 8 \pi n^{2} e^{4t}    - 6 e^{2t} + (2 \pi n^{2} e^{4t}    - 3 e^{2t}) (\frac{1}{2} - 2 \pi n^{2} e^{2t} ) ] \notag\\
 \frac{d X(t)}{dt} = \sum_{n=1}^{\infty} \pi n^{2} e^{- \pi n^{2} e^{2t}} e^{\frac{t}{2}} [ 8 \pi n^{2} e^{4t}    - 6 e^{2t} + ( \pi n^{2} e^{4t}    - \frac{3}{2} e^{2t} - 4 \pi^2 n^{4} e^{6t} +  6 \pi n^{2} e^{4t} ) ] \notag\\
  \frac{d X(t)}{dt} = \sum_{n=1}^{\infty} \pi n^{2} e^{- \pi n^{2} e^{2t}} e^{\frac{t}{2}} [  - 4 \pi^2 n^{4} e^{6t} + 15 \pi n^{2} e^{4t}    - \frac{15}{2} e^{2t}  ] \notag\\
\frac{d X(t)}{dt} =  \sum_{n=1}^{\infty} \pi n^{2} e^{- \pi n^{2} e^{2t}} e^{\frac{t}{2}}  e^{2t}  [  - 4 \pi^2 n^{4} e^{4t} + 15 \pi n^{2} e^{2t}    - \frac{15}{2} ]
\end{align}
We substitute $y = \pi e^{2t} $ in Eq.~\ref{sec:sec_a_1_eq_1_a} and define $A(y)$ such that $\frac{d X(t)}{dt} = \pi e^{\frac{5 t}{2}} A(y)$.  {\citep{TB}}

\begin{align}  \label{sec:sec_a_1_eq_1_b}   
A(y) =  \sum_{n=1}^{\infty}  n^{2} e^{-  n^{2} y}  [  - 4  n^{4} y^2 + 15  n^{2} y    - \frac{15}{2} ]
\end{align}

We see that $A(y)=0$ at $y=\pi$ which corresponds to $t=0$ given $y = \pi e^{2t} $ and $\frac{d X(t)}{dt} = \pi e^{\frac{5 t}{2}} A(y)$, given that $\frac{d X(t)}{dt} = 0$ at $t=0$. (\textbf{Result 6.1.a}) Because $X(t)= \frac{E_0(t)}{2} $ is an even function of variable $t$(Details in ~\ref{sec:appendix_C_7}) and hence $\frac{d X(t)}{dt}$ is  an \textbf{odd} function of variable $t$.\\

The quadratic expression $ B(y, n) = (- 4  n^{4} y^2 + 15  n^{2} y    - \frac{15}{2})$ in Eq.~\ref{sec:sec_a_1_eq_1_b} equals zero at $y = y_r =\frac{-15 n^2 \pm \sqrt{225 n^4 - 120 n^4}}{-8 n^4} = \frac{(15 \pm \sqrt{105})}{8 n^2}$ {\citep{TB}}. We see that the first derivative of $B(y, n)$ is given by $\frac{dB(y, n)}{dy} = - 8  n^{4} y + 15  n^{2}$ equals zero at $y=y_m=\frac{15}{8 n^2}$. The second derivative of $B(y, n)$ given by $\frac{d^2 B(y, n)}{dy^2} = - 8  n^{4}$, is negative for all $y$ and $n \geq 1$ and hence $B(y, n)$ is a \textbf{concave down} function for each $n$, which reaches a maximum at $y=y_m=\frac{15}{8 n^2}$ and has zeros at $y = y_r = \frac{(15 \pm \sqrt{105})}{8 n^2}$ and given the dominant term $ - 4  n^{4} y^2$ in Eq.~\ref{sec:sec_a_1_eq_1_b} , we see that $B(y, n) < 0$, for $ y  > \frac{(15 + \sqrt{105})}{8 } = 3.1559$, for $n \geq 1$ and hence $B(y, n) < 0$ for $ y > 3.16 = y_z$ and  hence $A(y) < 0$ for finite $y > y_z$. Using $y = \pi e^{2t}$ and  $\frac{d X(t)}{dt} = \pi e^{\frac{5 t}{2}} A(y)$, we see that  $\frac{d X(t)}{dt} < 0$ for $t > \frac{1}{2} \log{\frac{y_z}{\pi}} = t_z $(\textbf{Result 6.1}). (\href{https://archive.is/dpXEm}{concave down function})\\ 

We show in the next section that  $\frac{d X(t)}{dt} < 0$ for $0 < t \leq t_z$. It suffices to show that $\frac{d A(y)}{dy} < 0$  for $\pi \leq y \leq y_z=3.16$ and hence $A(y) < 0$ for $\pi < y \leq y_z=3.16$, given that  $A(y)=0$ at $y=\pi$, using Result 6.1.a. [ We use $y = \pi e^{2t} $ and $\frac{d X(t)}{dt} = \pi e^{\frac{5 t}{2}} A(y)$ and $\frac{d X(t)}{dt} =0$ at $t=0$.]

.\\ \tocless\subsection{\label{sec:Section_A_1_6_1} \textbf{$\frac{d X(t)}{dt} < 0$ for $0 < t \leq t_z$} \protect\\  \lowercase{} }

It is shown in this section that $\frac{d A(y)}{dy} < 0$  for $\pi \leq y \leq 3.16$ and hence $A(y) < 0$ for $\pi < y \leq 3.16$ {\citep{TB}} , given that  $A(y)=0$ at $y=\pi$. We take the  derivative of $A(y)$ in Eq.~\ref{sec:sec_a_1_eq_1_b} and take the factor $n^2$ out of the brackets in the last line below.
\begin{align}\label{sec:sec_A_1_eq_2}   
\frac{d A(y)}{dy} =  \sum_{n=1}^{\infty}  n^{2} e^{-  n^{2} y}  [  - 8  n^{4} y + 15  n^{2}  +(- 4  n^{4} y^2 + 15  n^{2} y    - \frac{15}{2} ) (-n^2)   ] \notag\\
\frac{d A(y)}{dy} =   \sum_{n=1}^{\infty}  n^{4} e^{-  n^{2} y}  [   -8  n^{2} y + 15    + 4  n^{4} y^2 - 15  n^{2} y    + \frac{15}{2}    ] = \sum_{n=1}^{\infty}  n^{4} e^{-  n^{2} y}  [  4  n^{4} y^2 - 23  n^{2} y  + \frac{45}{2}    ]
\end{align}
We examine the term $ C(y, n) = n^{4} e^{-  n^{2} y}  ( 4  n^{4} y^2 - 23  n^{2} y  + \frac{45}{2})$ in Eq.~\ref{sec:sec_A_1_eq_2} in the interval $\pi \leq y \leq 3.16$ and show that $\frac{d A(y)}{dy} = C(y, 1) + \sum_{n=2}^{\infty} C(y, n) < 0$, as follows.{\citep{TB}} We want the maximum value of $C(y, n)$ for each $n$ and we consider the maximum value of positive terms and minimum value of absolute value of negative terms in the paragraphs below. (\textbf{Result 6.2.a})\\

For $n=1$, we see that $C(y, 1) = e^{- y} ( 4  y^2 - 23 y  + \frac{45}{2}) =   4  y^2  e^{- y} - 23 y e^{- y}   + \frac{45}{2} e^{- y}  < 0$ in the interval $\pi \leq y \leq 3.16$ as follows. Given that $3.16^2 < 10$ and $\pi > 3.14$,  in the interval $\pi \leq y \leq 3.16$, we see that $ C(y, 1) <   4  * 10  e^{-3.14} - 23 * 3.14 e^{- 3.16}  + \frac{45}{2} e^{-3.14}  = -0.3588 <  -6 e^{- 3} = C_{max}(1) $ where $C_{max}(1)$ is the maximum value of $C(y,1)$ in the interval $\pi \leq y \leq 3.16$, using Result 6.2.a. 
\begin{align}  \label{sec:sec_A_1_eq_2_a}   
C(y, 1) = e^{- y} ( 4  y^2 - 23 y  + \frac{45}{2}) < -6 e^{- 3} , \quad  \pi \leq y \leq 3.16
\end{align}
For $n>1$, in the interval $\pi \leq y \leq 3.16$, we can write $C(y, n)$ as follows, given that $\pi > 3.14$ and $3.16^2 < 10$ and  the term $- 23  n^{2} y  < 0$ is omitted below, given that we want the maximum value of $C(y, n)$. We write the term $ \frac{45}{2} < 4  n^{4} * 0.5 $ and $e^{- 0.14 n^{2} } * 10.5 <10 $ for $n \geq 2$, using Result 6.2.a .
\begin{align}\label{sec:sec_A_1_eq_2_b}   
C(y, n) = n^{4} e^{-  n^{2} y}  ( 4  n^{4} y^2 - 23  n^{2} y  + \frac{45}{2}) < n^{4} e^{- \pi n^{2} }  ( 4  n^{4} ((3.16)^2+0.5) )  < 4 n^8  e^{- 3 n^{2} } e^{- 0.14 n^{2} } *10.5 < 40 n^8  e^{- 3 n^{2} }
\end{align}
We want to show that $\frac{d A(y)}{dy} = C(y, 1) + \sum_{n=2}^{\infty} C(y, n) < 0$ in the interval $\pi \leq y \leq 3.16$. Using Eq.~\ref{sec:sec_A_1_eq_2_a}  and 
Eq.~\ref{sec:sec_A_1_eq_2_b}, we write as follows. We multiply both sides by $e^{3}$ in the second line below.
\begin{align}\label{sec:sec_A_1_eq_2_c}   
\frac{d A(y)}{dy} = C(y, 1) + \sum_{n=2}^{\infty} C(y, n) < -6 e^{- 3}  +  \sum_{n=2}^{\infty}   40 n^8  e^{- 3 n^{2} } \notag\\
e^{3} \frac{d A(y)}{dy} < -6  + \sum_{n=2}^{\infty}   40 n^8  e^{3 - 3 n^{2} }
\end{align}
We want to show that $e^{3} \frac{d A(y)}{dy} < 0$ in the interval $\pi \leq y \leq 3.16$.  We compute $\log{(n^8  e^{3 - 3 n^{2}})}$ as follows. We note that $f(x)= \log{x} $ is a \textbf{concave down} function whose second derivative given by $-\frac{1}{x^2} < 0$ for $|x| < \infty$ and we can write $f(x)=\log{x} \leq f(x_0) + [\frac{df}{dx}]_{x=x_0} ( x - x_0)$ using its \textbf{tangent line} equation. We see that $\frac{df}{dx}=\frac{1}{x}$. We set $x=n$ and $x_0=2$ and get $\log{n} \leq \log{2} + \frac{1}{2} (n-2)$  below.
\begin{align}\label{sec:sec_A_1_eq_2_d}   
\log{(n^8  e^{3 - 3 n^{2}})} =8 \log{n} + (3 - 3 n^{2}) \leq  8 (\log{2} + \frac{1}{2} (n-2) ) + (3 - 3 n^{2}) \notag\\
\log{(n^8  e^{3 - 3 n^{2}})} \leq 8 \log{2} + 4 n - 5 - 3 n^{2} 
\end{align}
We note that $g(x)=  4 x - 5 - 3 x^{2} $ in Eq.~\ref{sec:sec_A_1_eq_2_d} is a \textbf{concave down} function (\href{https://archive.is/dpXEm}{concave down function}), whose second derivative given by $-6 < 0$ for all $x$ and we can write $g(x) \leq g(x_0) + [\frac{dg}{dx}]_{x=x_0} ( x - x_0)$ using its \textbf{tangent line} equation. We see that $\frac{dg}{dx}= 4 - 6 x$. We set $x=n$ and $x_0=2$ and get $g(n) \leq g(2) + [ 4 - 6 x]_{x=2} ( n - 2)= - 9 -8(n-2)  $ and write Eq.~\ref{sec:sec_A_1_eq_2_d} as follows. We take the exponent $e$ of both sides in the second line below.
\begin{align}\label{sec:sec_A_1_eq_2_e}   
\log{(n^8  e^{3 - 3 n^{2}})} \leq 8 \log{2} - 9 -8(n-2) \leq  8 \log{2} -1 + 8 (1-n) \notag\\
n^8  e^{3 - 3 n^{2}} \leq e^{8 \log{2}  -1 + 8 (1-n)} = 2^8 e^{-1} e^{8 (1-n)}
\end{align}
We substitute the result in Eq.~\ref{sec:sec_A_1_eq_2_e} in Eq.~\ref{sec:sec_A_1_eq_2_c} and simplify as follows. 
\begin{align}\label{sec:sec_A_1_eq_2_f}   
e^{3} \frac{d A(y)}{dy} < -6  +  40 * 2^8 e^{-1} \sum_{n=2}^{\infty}   e^{8 (1-n)} \notag\\
e^{3} \frac{d A(y)}{dy} < -6  +  40 * 2^8 e^{-1} * e^{8} \sum_{n=2}^{\infty}   e^{-8n} \notag\\
e^{3} \frac{d A(y)}{dy} < -6  +  40 * 2^8 e^{-1} * e^{8} \frac{e^{-8*2}}{1-e^{-8}} \notag\\
e^{3} \frac{d A(y)}{dy} < -6  +  40 * 2^8 e^{-1} *  \frac{e^{-8}}{1-e^{-8}} \notag\\
e^{3} \frac{d A(y)}{dy} < -6  +  40 * 2^8 e^{-1} *  \frac{1}{e^{8}-1} 
\end{align}
We multiply Eq.~\ref{sec:sec_A_1_eq_2_f} by $\frac{(e^{8}-1)}{6}$ and write as follows. The symbol $\approx$ means "approximately equals".
\begin{align} \label{sec:sec_A_1_eq_2_g}   
e^{3} \frac{d A(y)}{dy} \frac{(e^{8}-1)}{6} < - e^{8} + 1  +  40 e^{-1} * \frac{256}{6}  \approx -2352
\end{align}
We see that $e^{3} \frac{d A(y)}{dy} \frac{(e^{8}-1)}{6} <0$ in Eq.~\ref{sec:sec_A_1_eq_2_g}. Hence  $\frac{d A(y)}{dy}< 0$ in the interval $\pi \leq y \leq 3.16$, given that $e^{3}  \frac{(e^{8}-1)}{6} > 0$ and $e > 2$. Given that $A(y)=0$ at $y=\pi$ using Result 6.1.a, we see that $A(y) < 0$ in Eq.~\ref{sec:sec_a_1_eq_1_b} , for $\pi < y \leq 3.16$ and $\frac{d X(t)}{dt} = \pi e^{\frac{5 t}{2}} A(y) < 0$ in the interval $0 < t \leq t_z$.(\textbf{Result 6.2})\\

In Section~\ref{sec:Section_A_1_6_0}, it is shown that  $\frac{d X(t)}{dt} <0$  for $t > t_z $ (using Result 6.1). In this section, we have shown that  $\frac{d X(t)}{dt} <0$  for $0 < t \leq t_z$. Hence $\frac{d X(t)}{dt} <0$  for $t>0$.\\

Hence $E_0(t)= 2 X(t)$ is a \textbf{strictly decreasing function} for  $t > 0$ and strictly increasing for $t < 0 $, given that $E_0(t)=E_0(-t) > 0$ for $t > 0$ (Details in ~\ref{sec:appendix_C_5a} and ~\ref{sec:appendix_C_7}).

.\\ \tocless\section{\label{sec:Section_3} \textbf{Hurwitz Zeta Function and related functions }\protect  \lowercase{} }

We can show that the new method is \textbf{not} applicable to Hurwitz zeta function and related zeta functions and \textbf{does not} contradict the existence of their non-trivial zeros away from the critical line given by $Re[s] = \frac{1}{2}$. The new method requires the \textbf{symmetry} relation $\xi(s) = \xi(1-s)$ and hence $\xi(\frac{1}{2} + i \omega)=\xi(\frac{1}{2} - i \omega)$ when evaluated at the critical line $s = \frac{1}{2} + i \omega$. This means  $\xi(\frac{1}{2} + i \omega) = E_{0\omega}(\omega) = E_{0\omega}(-\omega)$ and $E_{0}(t)=E_{0}(-t)$ (Details in ~\ref{sec:appendix_C_7}) where $E_0(t) =  \sum_{n=1}^{\infty}  [ 4 \pi^{2} n^{4} e^{4t}    - 6 \pi n^{2}   e^{2t} ]  e^{- \pi n^{2} e^{2t}} e^{\frac{t}{2}} $ and this condition is satisfied for Riemann's Zeta function.\\

It is \textbf{not} known that Hurwitz Zeta Function  given by $\zeta(s, a) = \sum\limits_{m=0}^{\infty} \frac{1}{(m+a)^{s}}$ satisfies a symmetry relation similar to $\xi(s) = \xi(1-s)$ where $\xi(s)$ is an entire function, for $a \neq 1$ and hence the condition $E_{0}(t)=E_0(-t)$ is \textbf{not} known to be satisfied {\citep{WK1}}. Hence the new method is \textbf{not} applicable to Hurwitz zeta function and  \textbf{does not} contradict the existence of their non-trivial zeros away from the critical line.\\

Dirichlet L-functions satisfy a symmetry relation $\xi(s, \chi) = \epsilon(\chi) \xi(1-s, \bar{\chi)}$ {\citep{WK2}} which does \textbf{not} translate to $E_{0}(t)=E_0(-t)$ required by the new method and hence this proof is \textbf{not} applicable to them. This proof does not need or use Euler product.\\

We know that $\zeta(s) = \sum\limits_{m=1}^{\infty} \frac{1}{m^{s}}$ diverges for $Re[s] \leq 1$. Hence we derive a convergent and entire function $\xi(s)$ using the well known theorem $F(x) = 1 + 2 \displaystyle\sum_{n=1}^{\infty} e^{-\pi n^2 x} = \frac{1}{\sqrt{x}} (1 + 2 \displaystyle\sum_{n=1}^{\infty} e^{-\pi \frac{n^2}{x} })$, where $x > 0$ is real {\citep{FWE}}\href{https://www.ocf.berkeley.edu/~araman/files/math_z/Ellison_p147-152.pdf#page=5}{(link)} and then derive $E_0(t) =  \sum_{n=1}^{\infty}  [ 4 \pi^{2} n^{4} e^{4t}    - 6 \pi n^{2}   e^{2t} ]  e^{- \pi n^{2} e^{2t}} e^{\frac{t}{2}}   $. In the case of \textbf{Hurwitz zeta} function and \textbf{other zeta functions} with non-trivial zeros away from the critical line, it is \textbf{not} known if a corresponding relation similar to $F(x)$ exists, which enables derivation of a convergent and entire function $\xi(s)$ and results in $E_0(t)$ as a Fourier transformable, real, even and analytic function. Hence the new method presented in this paper is \textbf{not} applicable to Hurwitz zeta function and related zeta functions.\\

The proof of Riemann Hypothesis presented in this paper is \textbf{only} for the specific case of Riemann's Zeta function and \textbf{only} for the \textbf{critical strip} $0 \leq |\sigma| < \frac{1}{2}$. This proof requires both $E_{p}(t)$ and $E_{p\omega}(\omega)$ to be Fourier transformable where $E_p(t) = E_0(t) e^{-\sigma t}$ is a real analytic function and uses the fact that $E_0(t)$ is an \textbf{even} function of variable $t$ and  $ E_0(t)  >0 $ for $|t| < \infty$ (Details in ~\ref{sec:appendix_C_5a}) and $E_0(t)$ is \textbf{strictly decreasing} function for $t > 0$ (Details in Section~\ref{sec:Section_A_1_6}).  These conditions may \textbf{not} be satisfied for many other functions including those which have non-trivial zeros away from the critical line and hence the new method may \textbf{not} be applicable to such functions.\\

The case of analytic continuation of Riemann's zeta function derived from Dirichlet Eta function is discussed in ~\ref{sec:Section_2a} to ~\ref{sec:appendix_C_5b}.

.\\ \tocless\subsection{\textbf{Acknowledgments}}

I am very grateful to Srinivas M. Aji and Thomas Browning for their detailed feedback on my manuscript drafts and suggestions. My thanks to Brad Rodgers, Ken Ono and Anurag Sahay for helpful suggestions. I am grateful to Bhaskar Ramamurthi for introducing Fourier Transforms and inspiring my interest in mathematics. My thanks to Rukmini Dey whose love for mathematics inspired me. I would like to thank my parents for their encouragement of my mathematical pursuits and my father who introduced numbers to me and inspired my interest in mathematics.

\appendix.\\ \tocless\section{\label{sec:appendix_A} \textbf{Derivation of $E_p(t)$} \protect\\  \lowercase{} }

Let us start with Riemann's Xi Function $\xi(s)$ evaluated at $s = \frac{1}{2} + i \omega$ given by $\xi(\frac{1}{2} + i \omega)= E_{0\omega}(\omega)$. Its inverse Fourier Transform is given by $ E_0(t)=  \frac{1}{2 \pi} \int_{-\infty}^{\infty} E_{0\omega}(\omega) e^{i\omega t} d\omega =  \sum_{n=1}^{\infty}  [ 4 \pi^{2} n^{4} e^{4t}    - 6 \pi n^{2}   e^{2t} ]  e^{- \pi n^{2} e^{2t}} e^{\frac{t}{2}}$ using Eq.~\ref{sec_intro_eq_1}. \\ 

We will show in this section that the inverse Fourier Transform of the function  $\xi(\frac{1}{2} + \sigma + i \omega) =  E_{p\omega}(\omega)$,  is given by $E_{p}(t) =   E_{0}(t) e^{-\sigma t}$  where $0 < |\sigma| < \frac{1}{2}$ is real. We use $E_{p\omega}(\omega) = E_{0\omega}(\omega - i \sigma)$ below. 
\begin{align}\label{sec_app_A_1_eq_1}
\xi(\frac{1}{2} + \sigma + i \omega) = \xi(\frac{1}{2} + i (\omega - i \sigma) ) = E_{p\omega}(\omega) = E_{0\omega}(\omega - i \sigma) \notag\\
E_{p}(t) =   \frac{1}{2 \pi} \int_{-\infty}^{\infty} E_{p\omega}(\omega) e^{i\omega t} d\omega =   \frac{1}{2 \pi} \int_{-\infty}^{\infty} E_{0\omega}(\omega - i \sigma) e^{i\omega t} d\omega 
\end{align}
 We substitute $\omega' = \omega - i \sigma$ in Eq.~\ref{sec_app_A_1_eq_1} as follows. We get $\omega=\omega'+i \sigma$ and $d\omega=d\omega'$.
\begin{align}  \label{sec_app_A_1_eq_1_1}
E_{p}(t) =    e^{-\sigma t}  \frac{1}{2 \pi} \int_{-\infty  - i \sigma}^{\infty  - i \sigma} E_{0\omega}(\omega') e^{i\omega' t} d\omega'
 \end{align}
We can evaluate the above integral in the complex plane using contour integration, substituting $\omega' = z = x + i y$ and we use a rectangular contour comprised of $C_1$ along the line $ z = (-\infty, \infty)$, $C_2$ along the line $z = (\infty, \infty-i \sigma)$, $C_3$ along the line $z = ( \infty-i\sigma, -\infty-i\sigma)$ and then $C_4$ along the line $z = (-\infty-i \sigma, -\infty)$. We can see that $E_{0\omega}(z)=\xi(\frac{1}{2}+ i z)$ has no singularities in the region bounded by the contour because $\xi(\frac{1}{2}+ i z)$ is an entire function in the complex plane.\\

We use the fact that $E_{0\omega}(z) = \xi(\frac{1}{2}+ i z) = \xi(\frac{1}{2} -y + i x) =  \int_{-\infty}^{\infty} E_0(t) e^{-i z t} dt  =  \int_{-\infty}^{\infty} E_0(t) e^{y t} e^{-i x t} dt $, \textbf{goes to zero} as $x \to \pm \infty$ when $-\sigma \leq y \leq 0$, as per Riemann-Lebesgue Lemma \href{https://archive.is/HV6zJ}{(link)}, because $E_0(t) e^{y t}$ is a absolutely integrable function   for real $t$(Details in ~\ref{sec:appendix_C_6}). Hence the integral in Eq.~\ref{sec_app_A_1_eq_1_1} \textbf{vanishes }along  the contours $C_2$ and $C_4$. Using Cauchy's Integral theroem, we can write Eq.~\ref{sec_app_A_1_eq_1_1} as follows.
\begin{align}\label{sec_app_A_1_eq_2}
E_{p}(t) =    e^{-\sigma t}  \frac{1}{2 \pi} \int_{-\infty}^{\infty} E_{0\omega}(\omega') e^{i\omega' t} d\omega'  \notag\\
E_{p}(t) =    E_{0}(t) e^{-\sigma t} =  \sum_{n=1}^{\infty}  [ 4 \pi^{2} n^{4} e^{4t}    - 6 \pi n^{2}   e^{2t} ]  e^{- \pi n^{2} e^{2t}} e^{\frac{t}{2}} e^{-\sigma t} 
\end{align}
Thus we have arrived at the desired result $E_{p}(t) =   E_{0}(t) e^{-\sigma t}$. 
\textbf{Alternate} derivation of $E_0(t)$ and $E_p(t)$ are in ~\ref{sec:appendix_H_1}.

.\\ \tocless\subsection{\label{sec:appendix_C_6} \textbf{$E_y(t) = E_0(t)  e^{y t}$ is an  absolutely integrable function} \protect\\  \lowercase{} }

We see that $E_0(t) > 0$ and finite for $-\infty < t < \infty$ (Details in ~\ref{sec:appendix_C_5a}). Hence $E_y(t) = E_0(t)  e^{y t} > 0 $ and finite for all $-\infty < t < \infty$, for $- \sigma \leq y \leq 0$ and $0 \leq |\sigma| < \frac{1}{2}$ (\textbf{Result A.1.1}).\\

$E_0(t)$ has an asymptotic \textbf{exponential} fall-off rate of $o[e^{-1.5 |t| } ]   $ (using ~\ref{sec:appendix_C_5}) and hence $E_y(t) = E_0(t)  e^{y t}$ has an asymptotic \textbf{exponential} fall-off rate of $o[e^{- |t|} ] $ (using $o[e^{-(1.5 + y) |t| } ]$ ), for $- \sigma \leq y \leq 0$ and $0 \leq |\sigma| < \frac{1}{2}$. Hence $E_y(t)= E_0(t)  e^{y t}$ decays exponentially, at $t \to \pm \infty$.(\textbf{Result A.1.2}) \\

Using Result A.1.1 and A.1.2, we can write $\int_{-\infty}^{\infty} |E_y(t)| dt$ is finite and $E_y(t)$ is an  absolutely \textbf{integrable function} (Details in ~\ref{sec:appendix_C_5z}) and its Fourier transform $ E_{y\omega}(\omega)$ goes to zero as $\omega \to \pm \infty$, as per Riemann Lebesgue Lemma \href{https://archive.is/HV6zJ}{(link)}.

.\\ \tocless\section{\label{sec:appendix_H} \textbf{Details of entire function $\xi(s)$ } \protect\\  \lowercase{} }

In this section, we will start with Riemann's Xi function $\xi(s)$ and take the inverse Fourier Transform of  $\xi(\frac{1}{2} + i \omega)= E_{0\omega}(\omega)$ and show the result $E_{0}(t) =    \sum_{n=1}^{\infty}  [ 4 \pi^{2} n^{4} e^{4t}    - 6 \pi n^{2}   e^{2t} ]  e^{- \pi n^{2} e^{2t}} e^{\frac{t}{2}} $ and $E_p(t)=E_0(t) e^{-\sigma t}$. \\
 
We will use the equation for $\xi(s)$ derived in Ellison's book "Prime Numbers" pages 151-152  which uses \textbf{the well known theorem} $ 1 + 2 w(x) = \frac{1}{\sqrt{x}} (1 + 2 w(\frac{1}{x}) )$, where  $w(x)= \displaystyle\sum_{n=1}^{\infty} e^{-\pi n^2 x}$ and $x > 0$ is real.{\citep{FWE}} \href{https://www.ocf.berkeley.edu/~araman/files/math_z/Ellison_p147-152.pdf#page=5}{(link)}. 
\begin{align} \label{sec:App_H_eq_3_1} 
\xi(s) = \frac{1}{2} s (s-1) \Gamma(\frac{s}{2}) \pi^{-\frac{s}{2}} \zeta(s) =  \frac{1}{2} [ 1 +  s (s-1)  \int_{1}^{\infty} ( x^{\frac{s}{2}} + x^{\frac{1-s}{2}}) w(x)  \frac{dx}{x}  ]
\end{align}
We see that $\xi(s)$ is an entire function, for all values of  $s$ in the complex plane and hence we get an analytic continuation of $\xi(s)$ over the entire complex plane. We see that $\xi(s) = \xi(1-s)$  {\citep{FWE}}.

.\\ \tocless\subsection{\label{sec:appendix_H_1} \textbf{Derivation of $E_p(t)$ and $E_0(t)$ } \protect\\  \lowercase{} }

Given that $w(x) = \sum\limits_{n=1}^{\infty} e^{- \pi n^{2}x }$, we substitute $x= e^{2t}, \frac{dx}{x}= 2 dt$ in Eq.~\ref{sec:App_H_eq_3_1} and evaluate at  $s= \frac{1}{2} + \sigma + i \omega$ as follows.
\begin{align}  \label{sec:App_H_eq_4} 
\xi(\frac{1}{2}  + \sigma + i \omega) =   \frac{1}{2} [ 1 +  2 (\frac{1}{2} + \sigma + i \omega) (-\frac{1}{2} + \sigma + i \omega)  \int_{0}^{\infty}  \sum\limits_{n=1}^{\infty} e^{- \pi n^{2} e^{2t} }  ( e^{\frac{t}{2}} e^{\sigma t} e^{ i \omega t} +  e^{\frac{t}{2}} e^{-\sigma t} e^{ -i \omega t} )  dt   ]
\end{align}
We can substitute $t=-t$ in the first term in above integral and simplify above equation as follows.
\begin{align}\label{sec:App_H_eq_5} 
\xi(\frac{1}{2}  + \sigma + i \omega) =   \frac{1}{2}  +   (- \frac{1}{4} +  \sigma^{2} -  \omega^{2} + i \omega (2 \sigma)  )  [ \int_{-\infty}^{0}  \sum\limits_{n=1}^{\infty} e^{- \pi n^{2} e^{-2t} }  e^{\frac{-t}{2}} e^{-\sigma t} e^{ -i \omega t} dt \notag\\
+ \int_{0}^{\infty}  \sum\limits_{n=1}^{\infty} e^{- \pi n^{2} e^{2t} }   e^{\frac{t}{2}} e^{-\sigma t} e^{ -i \omega t}   dt  ] 
\end{align}
We can write this as follows.
\begin{align}  \label{sec:App_H_eq_6} 
\xi(\frac{1}{2}  + \sigma + i \omega) =   \frac{1}{2}  +   (- \frac{1}{4} +  \sigma^{2} -  \omega^{2} + i \omega (2 \sigma)  )   \int_{-\infty}^{\infty} [ \sum\limits_{n=1}^{\infty} e^{- \pi n^{2} e^{-2t} }  e^{\frac{-t}{2}} u(-t) + \sum\limits_{n=1}^{\infty} e^{- \pi n^{2} e^{2t} }   e^{\frac{t}{2}}   u(t) ]  e^{-\sigma t} e^{ -i \omega t} dt 
\end{align}
We define $A(t) =  [ \sum\limits_{n=1}^{\infty} e^{- \pi n^{2} e^{-2t} }  e^{\frac{-t}{2}} u(-t) + \sum\limits_{n=1}^{\infty} e^{- \pi n^{2} e^{2t} }   e^{\frac{t}{2}}   u(t) ]  e^{-\sigma t}  $ and get the \textbf{inverse Fourier transform} of $\xi(\frac{1}{2} + \sigma + i \omega) = E_{p\omega}(\omega)$ in above equation given by $E_p(t)$ as follows. We use dirac delta function  $ \delta(t) $.
\begin{align}\label{sec:App_H_eq_10a} 
E_p(t) =  \frac{1}{2} \delta(t) + (- \frac{1}{4} +  \sigma^{2}) A(t) + 2 \sigma \frac{dA(t)}{dt}  + \frac{d^{2}A(t)}{dt^{2}}  \notag\\
A(t) = [ \sum\limits_{n=1}^{\infty} e^{- \pi n^{2} e^{-2t} }  e^{\frac{-t}{2}} u(-t) + \sum\limits_{n=1}^{\infty} e^{- \pi n^{2} e^{2t} }   e^{\frac{t}{2}}   u(t) ]  e^{-\sigma t} 
\end{align}
We compute the derivatives of $A(t)$ as follows. We note that $[A(t)]_{t=0+} = [A(t)]_{t=0-} =\sum\limits_{n=1}^{\infty} e^{- \pi n^{2}} $.
\begin{align}\label{sec:App_H_eq_10} 
\frac{dA(t)}{dt} = \sum\limits_{n=1}^{\infty} e^{- \pi n^{2} e^{-2t} } e^{\frac{-t}{2}}  e^{-\sigma t} [ -\frac{1}{2} - \sigma + 2 \pi n^{2} e^{-2t} ] u(-t) + \sum\limits_{n=1}^{\infty} e^{- \pi n^{2} e^{2t} } e^{\frac{t}{2}}  e^{-\sigma t} [ \frac{1}{2} - \sigma - 2 \pi n^{2} e^{2t} ] u(t) \notag\\
\frac{d^{2}A(t)}{dt^{2}}  = \sum\limits_{n=1}^{\infty} e^{- \pi n^{2} e^{-2t} } e^{\frac{-t}{2}} e^{-\sigma t} [ - 4 \pi n^{2} e^{-2t} + (-\frac{1}{2} - \sigma + 2 \pi n^{2} e^{-2t} )^{2} ] u(-t) \notag\\+  \sum\limits_{n=1}^{\infty} e^{- \pi n^{2} e^{2t} } e^{\frac{t}{2}} e^{-\sigma t} [ - 4 \pi n^{2} e^{2t} + (\frac{1}{2} - \sigma - 2 \pi n^{2} e^{2t} )^{2} ] u(t) + A_0  \delta(t) 
\end{align}
We use $A_0=[\frac{dA(t)}{dt}]_{t=0+} - [\frac{dA(t)}{dt}]_{t=0-} = \sum\limits_{n=1}^{\infty} e^{- \pi n^{2}} (\frac{1}{2} - \sigma - 2 \pi n^{2} -( -\frac{1}{2} - \sigma + 2 \pi n^{2} )  ) = \sum\limits_{n=1}^{\infty} e^{- \pi n^{2}} (1 - 4 \pi n^{2} ) $. We can simplify above equation as follows.
\begin{align}\label{sec:App_H_eq_11} 
\frac{d^{2}A(t)}{dt^{2}}  =   \sum\limits_{n=1}^{\infty} e^{- \pi n^{2} e^{-2t} } e^{\frac{-t}{2}}  e^{-\sigma t}  [ \frac{1}{4} + \sigma^2 +\sigma +   4 \pi^{2} n^{4} e^{-4t}  - 6 \pi n^{2} e^{-2t} - 4 \sigma \pi n^{2} e^{-2t}  ] u(-t) \notag\\
+ \sum\limits_{n=1}^{\infty} e^{- \pi n^{2} e^{2t} } e^{\frac{t}{2}}  e^{-\sigma t}  [ \frac{1}{4} + \sigma^2 -\sigma +   4 \pi^{2} n^{4} e^{4t}  - 6 \pi n^{2} e^{2t} + 4 \sigma \pi n^{2} e^{2t}  ] u(t) + \delta(t) [\sum\limits_{n=1}^{\infty} e^{- \pi n^{2}} (1 - 4 \pi n^{2} ) ] 
\end{align}
We use the fact that  $F(x) = 1 + 2 w(x) = \frac{1}{\sqrt{x}} (1 + 2 w(\frac{1}{x}) )$, where $w(x) = \sum\limits_{n=1}^{\infty} e^{- \pi n^{2}x }$ and $x > 0$ is real {\citep{FWE}}, and we take the first derivative of $F(x)$ and evaluate it at $x=1$ and get $\sum\limits_{n=1}^{\infty} e^{- \pi n^{2}} (1 - 4 \pi n^{2} )  = -\frac{1}{2}$ (Details in ~\ref{sec:appendix_H_2}) and hence \textbf{dirac delta terms cancel each other} in  Eq. ~\ref{sec:App_H_eq_10a} written as follows using Eq. ~\ref{sec:App_H_eq_10}  and Eq. ~\ref{sec:App_H_eq_11}. 
\begin{align}\label{sec:App_H_eq_12} 
E_p(t) =  \frac{1}{2} \delta(t) + (- \frac{1}{4} +  \sigma^{2}) A(t) + 2 \sigma \frac{dA(t)}{dt}  + \frac{d^{2}A(t)}{dt^{2}}  \notag\\
E_p(t) =   \sum\limits_{n=1}^{\infty} e^{- \pi n^{2} e^{-2t} } e^{\frac{-t}{2}}  e^{-\sigma t}  [ -\frac{1}{4} + \sigma^2  + 2 \sigma( -\frac{1}{2} - \sigma + 2 \pi n^{2} e^{-2t} ) \notag\\+ \frac{1}{4} + \sigma^2 +\sigma +   4 \pi^{2} n^{4} e^{-4t}  - 6 \pi n^{2} e^{-2t} - 4 \sigma \pi n^{2} e^{-2t}  ] u(-t) \notag\\
+ \sum\limits_{n=1}^{\infty} e^{- \pi n^{2} e^{2t} } e^{\frac{t}{2}}  e^{-\sigma t}  [ -\frac{1}{4} + \sigma^2  + 2 \sigma( \frac{1}{2} - \sigma - 2 \pi n^{2} e^{2t} ) + \frac{1}{4} + \sigma^2 -\sigma +   4 \pi^{2} n^{4} e^{4t}  - 6 \pi n^{2} e^{2t} + 4 \sigma \pi n^{2} e^{2t}  ] u(t) \notag\\
E_p(t) =  \sum\limits_{n=1}^{\infty} e^{- \pi n^{2} e^{-2t} } e^{\frac{-t}{2}}  e^{-\sigma t} D(t, n) u(-t) + \sum\limits_{n=1}^{\infty} e^{- \pi n^{2} e^{2t} } e^{\frac{t}{2}}  e^{-\sigma t} C(t, n) u(t)
\end{align}
We cancel the common terms in Eq.~\ref{sec:App_H_eq_12}  and simplify above equation as follows. 
\begin{align}\label{sec:App_H_eq_12_a} 
C(t, n)=-\frac{1}{4} + \sigma^2  + \sigma - 2 \sigma^2 - 4 \sigma \pi n^{2} e^{2t}  + \frac{1}{4} + \sigma^2 -\sigma +   4 \pi^{2} n^{4} e^{4t}  - 6 \pi n^{2} e^{2t} + 4 \sigma \pi n^{2} e^{2t} \notag\\
C(t, n)  = 4 \pi^{2} n^{4} e^{4t}  - 6 \pi n^{2} e^{2t} \notag\\
D(t, n)=-\frac{1}{4} + \sigma^2  - \sigma - 2 \sigma^2 + 4 \sigma \pi n^{2} e^{-2t}  + \frac{1}{4} + \sigma^2 +\sigma +   4 \pi^{2} n^{4} e^{-4t}  - 6 \pi n^{2} e^{-2t} - 4 \sigma \pi n^{2} e^{-2t} \notag\\ 
D(t, n) = 4 \pi^{2} n^{4} e^{-4t}  - 6 \pi n^{2} e^{-2t}
\end{align}
We see that $D(t, n)= C(-t, n)$. Hence we can write Eq.~\ref{sec:App_H_eq_12} as follows.
\begin{align}\label{sec:App_H_eq_13} 
E_p(t) =  [   E_0(-t) u(-t) + E_0(t) u(t) ] e^{-\sigma t} \notag\\
E_0(t) = \sum\limits_{n=1}^{\infty} C(t, n)  e^{- \pi n^{2} e^{2t} } e^{\frac{t}{2}} =  \sum\limits_{n=1}^{\infty} [  4 \pi^{2} n^{4} e^{4t}  - 6 \pi n^{2} e^{2t}  ]  e^{- \pi n^{2} e^{2t} } e^{\frac{t}{2}} 
\end{align}
We use the fact that $E_{0}(t)=E_{0}(-t)$ (Details in ~\ref{sec:appendix_C_7}) and $u(t)+u(-t)=1$ and  we arrive at the desired result for $E_p(t)$ as follows. 
\begin{align}\label{sec:App_H_eq_14} 
E_0(t) =   \sum\limits_{n=1}^{\infty} [  4 \pi^{2} n^{4} e^{4t}  - 6 \pi n^{2} e^{2t}  ]  e^{- \pi n^{2} e^{2t} } e^{\frac{t}{2}} \notag\\
E_p(t) =  E_0(t) e^{-\sigma t} =   \sum\limits_{n=1}^{\infty} [  4 \pi^{2} n^{4} e^{4t}  - 6 \pi n^{2} e^{2t}  ]  e^{- \pi n^{2} e^{2t} } e^{\frac{t}{2}} e^{-\sigma t} 
\end{align}

.\\ \tocless\subsection{\label{sec:appendix_H_2} \textbf{Derivation of  $\sum\limits_{n=1}^{\infty} e^{- \pi n^{2}} (1 - 4 \pi n^{2} )  = -\frac{1}{2}$  } \protect\\  \lowercase{} }

In this section, we derive $\sum\limits_{n=1}^{\infty} e^{- \pi n^{2}} (1 - 4 \pi n^{2} )  = -\frac{1}{2}$. We use the fact that  $F(x) = 1 + 2 w(x) = \frac{1}{\sqrt{x}} (1 + 2 w(\frac{1}{x}) )$, where $w(x) = \sum\limits_{n=1}^{\infty} e^{- \pi n^{2}x }$ and $x > 0$ is real {\citep{FWE}}, and we take the first derivative of $F(x)$ and evaluate it at $x=1$.
\begin{align}\label{sec:App_H_2_eq_1} 
F(x) = 1 + 2 w(x) = \frac{1}{\sqrt{x}} (1 + 2 w(\frac{1}{x}) ) \notag\\
F(x) = 1 + 2 \sum\limits_{n=1}^{\infty} e^{- \pi n^{2}x } = \frac{1}{\sqrt{x}} (1 + 2 \sum\limits_{n=1}^{\infty} e^{- \pi n^{2}\frac{1}{x} } ) \notag\\
\frac{dF(x)}{dx} = 2 \sum\limits_{n=1}^{\infty} (- \pi n^{2}) e^{- \pi n^{2}x } = \frac{1}{\sqrt{x}}  \sum\limits_{n=1}^{\infty} (2 \pi n^{2}) e^{- \pi n^{2}\frac{1}{x} } (\frac{1}{x^2})  + (1 + 2 \sum\limits_{n=1}^{\infty} e^{- \pi n^{2}\frac{1}{x} } )  (\frac{-1}{2}) \frac{1}{x^{\frac{3}{2}}}
\end{align}
We evaluate the above equation at $x=1$ and we simplify as follows.
\begin{align}\label{sec:App_H_2_eq_2} 
[ \frac{dF(x)}{dx} ]_{x=1} = 2 \sum\limits_{n=1}^{\infty} (- \pi n^{2}) e^{- \pi n^{2} } =  \sum\limits_{n=1}^{\infty} (2 \pi n^{2}) e^{- \pi n^{2} }   + (1 + 2 \sum\limits_{n=1}^{\infty} e^{- \pi n^{2} } )  (\frac{-1}{2})  \notag\\
\sum\limits_{n=1}^{\infty} e^{- \pi n^{2}} (1 - 4 \pi n^{2} )  = -\frac{1}{2}
\end{align}

.\\ \tocless\section{\label{sec:appendix_C} Properties of Fourier Transforms  \protect  \lowercase{} }

.\\ \tocless\subsection{\label{sec:appendix_C_1} \textbf{ $E_p(t), h(t)$ are absolutely integrable functions and their Fourier Transforms are finite. } \protect\\  \lowercase{} }

The inverse Fourier Transform of the function $ E_{p\omega}(\omega)=\xi(\frac{1}{2}+ \sigma + i \omega)$ is given by $E_p(t) = E_0(t) e^{-\sigma t}$ (~\ref{sec:appendix_A}). In Eq.~\ref{sec_intro_eq_1}, we see that $E_0(t) =  \sum_{n=1}^{\infty}  [ 4 \pi^{2} n^{4} e^{4t}    - 6 \pi n^{2}   e^{2t} ]  e^{- \pi n^{2} e^{2t}} e^{\frac{t}{2}} > 0 $  and finite for all $-\infty <  t < \infty$(Details in ~\ref{sec:appendix_C_5a}). Hence $E_p(t) = E_0(t) e^{-\sigma t} > 0 $ and finite for all $-\infty < t < \infty$ and $0 < \sigma < \frac{1}{2}$.\\

It is shown in ~\ref{sec:appendix_C_5} that $E_0(t)$ has an asymptotic \textbf{exponential} fall-off rate of $o[e^{-1.5 |t| } ] $ and hence
$E_p(t)= E_0(t) e^{-\sigma t}$ has an asymptotic \textbf{exponential} fall-off rate of $ o[e^{- |t|} ] $ (using $e^{-(1.5 -  \sigma) |t| } $), for $0 < \sigma < \frac{1}{2}$. Hence $E_p(t)$ goes to zero, at $t \to \pm \infty$ and we showed that $E_p(t) > 0 $  and finite for all $-\infty < t < \infty$ in the last paragraph.(\textbf{Result C.1.1})  Hence $E_{p\omega}(\omega) = \int_{-\infty}^{\infty} E_p(t) e^{-i \omega t} dt$, evaluated at $\omega=0$ \textbf{cannot} be zero. Hence $E_{p\omega}(\omega)$ \textbf{does not have a zero} at $\omega=0$ and hence $\omega_{0} \neq 0$.\\

Given that $\xi(\frac{1}{2} + \sigma + i \omega)= E_{p\omega}(\omega)$ is an entire function in the whole of s-plane, it is finite for real $\omega$ and also for $\omega=0$. Hence $E_{p\omega}(0) = \int_{-\infty}^{\infty} E_p(t) dt$ is finite. Using Result C.1.1, we can write $\int_{-\infty}^{\infty} |E_p(t)| dt$ is finite and $E_p(t)$ is an  absolutely \textbf{integrable function} and its Fourier transform $ E_{p\omega}(\omega)$ converges and goes to zero as $\omega \to \pm \infty$, as per Riemann Lebesgue Lemma \href{https://archive.is/HV6zJ}{(link)}.\\

Using the arguments in above paragraph, we replace $\sigma$ in $E_p(t)= E_0(t) e^{-\sigma t}$ by $0$ and $2 \sigma$ respectively and see that $E_0(t)$ and $E_0(t) e^{-2 \sigma t}$
are absolutely \textbf{integrable} functions and the integrals  $\int_{-\infty}^{\infty} |E_0(t)| dt < \infty$ and  $\int_{-\infty}^{\infty} |E_0(t) e^{-2 \sigma t}| dt < \infty$.\\

Given that $E_p(t)= E_0(t) e^{-\sigma t}$ is an absolutely integrable function, its shifted versions are absolutely integrable and we see that $E_{p}(t + t_0)$ and $E_{p}(t - t_0)$  in Eq.~\ref{sec:sec_2_1_eq_0} are absolutely integrable functions, for a  finite shift of $ t_0$. ( We substitute $t- t_0=\tau$ and $dt=d\tau$ and get $\int_{-\infty}^{\infty} |E_p(t - t_0)| dt = \int_{-\infty}^{\infty} |E_p(\tau)| d\tau < \infty$ and hence $E_p(t - t_0)$ is an absolutely integrable function, given that $E_p(t)$ is absolutely integrable. Same argument holds for $E_p(t +  t_0)$.) \\

We see that $h(t)=  e^{ \sigma t} u(-t) + e^{ - \sigma t} u(t) $ is an absolutely \textbf{integrable function} because $h(t)>0$ for real $t$ and $\int_{-\infty}^{\infty} |h(t)| dt = \int_{-\infty}^{\infty} h(t) dt =  [ \int_{-\infty}^{\infty} h(t) e^{-i \omega t} dt ]_{\omega=0} = [ \frac{1}{\sigma - i \omega} +  \frac{1}{\sigma + i \omega} ]_{\omega=0} = \frac{2}{\sigma}$, is finite for $0 < \sigma < \frac{1}{2}$ and its Fourier transform $ H(\omega)$ converges and goes to zero as $\omega \to \pm \infty$, as per Riemann Lebesgue Lemma \href{https://archive.is/HV6zJ}{(link)}.

.\\ \tocless\subsection{\label{sec:appendix_C_2} \textbf{ Convolution integral convergence  } \protect\\  \lowercase{} }

Let us consider $h(t)=  e^{ \sigma t} u(-t) + e^{ - \sigma t} u(t) $  whose first derivative given by $\frac{d h(t)}{dt} =  \sigma e^{ \sigma t} u(-t)  - \sigma e^{ - \sigma t} u(t) $ 
and $A_0=[\frac{d h(t)}{dt}]_{t=0+} - [\frac{d h(t)}{dt}]_{t=0-} = -2 \sigma \neq 0$, for $0 < \sigma < \frac{1}{2}$ and hence $\frac{d h(t)}{dt}$ is \textbf{discontinuous} at $t=0$. The second derivative of $h(t)$ given by $h_2(t)$ has a Dirac delta function $A_0 \delta(t)$ where $A_0=-2 \sigma$ and its Fourier transform $H_2(\omega)$ has a constant term $A_0$, corresponding to the Dirac delta function.\\

This means $h(t)$ is obtained by integrating $h_2(t)$ twice and its Fourier transform $H(\omega)$ has a term $\frac{A_0}{(i \omega)^2}$ \href{https://web.archive.org/web/20240218075832/https://engineering.purdue.edu/~mikedz/ee301/FourierTransformTable.pdf}{(link)}  and has a \textbf{fall off rate}  of $\frac{1}{\omega^2}$ as $|\omega| \to \infty$. (\textbf{Result C.2}) \\

We see that $G(\omega', t_0)$ converges for real $\omega'$ (Eq.~\ref{sec:sec_2_1_eq_1_2} and Eq.~\ref{sec:sec_2_1_eq_9}) and $G(\omega -\omega', t_0) $ also converges for real $\omega$ and both have a minimum fall-off rate of $\frac{1}{\omega^0}=1$ as $|\omega| \to \infty$, given that inverse Fourier transform of $G(\omega', t_0)$ given by $g(t, t_0)$ converges. Using Result C.2, the integrand  $G(\omega -\omega', t_0) H(\omega')$ in Eq.~\ref{sec_C_1_eq_2} has a \textbf{minimum} fall off rate of $\frac{1}{\omega^2}$ as $|\omega| \to \infty$. Hence the convolution integral below converges to a finite value for  real $\omega$. 
\begin{align} \label{sec_C_1_eq_2}   
F(\omega, t_0) = \frac{1}{2 \pi}  [ G(\omega, t_0) \ast H(\omega)] = \frac{1}{2 \pi}  \int_{-\infty}^{\infty} G(\omega', t_0) H(\omega - \omega') d\omega' = \frac{1}{2 \pi}  \int_{-\infty}^{\infty} G(\omega -\omega', t_0) H(\omega') d\omega'  
\end{align}

.\\ \tocless\subsection{\label{sec:appendix_C_5} \textbf{ Exponential Fall off rate of $E_0(t)$, $E_p(t)$ and  $x(t)=E_0(t) e^{-2 \sigma t}$ } \protect\\  \lowercase{} }

We can write $E_0(t) =  \sum_{n=1}^{\infty}  [ 4 \pi^{2} n^{4} e^{4t}    - 6 \pi n^{2}   e^{2t} ]  e^{- \pi n^{2} e^{2t}} e^{\frac{t}{2}} $ in Eq.~\ref{sec_intro_eq_1} as follows. We take the term $ 2 \pi n^{2}   e^{2t}$ out of the brackets below. In the term $e^{- \pi n^{2} e^{2t}} $, we use Taylor series expansion around $t=0$ for $e^{2t}= \displaystyle\sum_{r=0}^{\infty} \frac{(2t)^{r}}{!r}$, given that $e^{2t}$ is an analytic function for real $t$.
\begin{align}\label{sec_app_C_5_eq_1}
E_0(t) =  \sum_{n=1}^{\infty} 2 \pi n^{2}   e^{2t} [ 2 \pi n^{2} e^{2t}    - 3 ]  e^{- \pi n^{2} e^{2t}} e^{\frac{t}{2}} \notag\\
= \sum_{n=1}^{\infty} 2 \pi n^{2}   e^{2t} [ 2 \pi n^{2} e^{2t}    - 3 ]  e^{- \pi n^{2} (1+ 2t)} e^{- \pi n^{2} ( \frac{(2t)^2}{!2} +  \frac{(2t)^3}{!3} ....)} e^{\frac{t}{2}} 
\end{align}
We take the term $ e^{-2 \pi t}$ out of the summation, corresponding to $n=1$ and then take the term $ 2 \pi e^{4t} e^{\frac{t}{2}}  = 2 \pi e^{\frac{9 t}{2}} $ out and write Eq.~\ref{sec_app_C_5_eq_1} as follows.
\begin{align} \label{sec_app_C_5_eq_2}
E_0(t) =  2 \pi e^{-2 \pi t}  e^{\frac{9 t}{2}}  \sum_{n=1}^{\infty} n^{2}   [ 2 \pi n^{2}   - 3 e^{-2t} ]  e^{- \pi n^{2}}   e^{- 2 \pi (n^{2}-1) t} e^{- \pi n^{2} ( \frac{(2t)^2}{!2} +  \frac{(2t)^3}{!3} ....)} 
 \end{align}
For $t > 0$, we see that the term corresponding to $n=1$ in Eq.~\ref{sec_app_C_5_eq_2} has an asymptotic fall-off rate of  $o[e^{-1.5 t} ]$(using $e^{- (2 \pi - \frac{9}{2}) t} $). The terms corresponding to $n > 1$ have higher fall-off rates, due to the term $e^{- 2 \pi (n^{2}-1) t}$.\\

Hence we see that $E_0(t)$ has an asymptotic fall-off rate of $o[e^{-1.5 t} ]$, for $t>0$. Given that $E_0(t)=E_0(-t)$(Details in ~\ref{sec:appendix_C_7}), we see that $E_0(t)$ has an \textbf{exponential} asymptotic fall-off rate of  $o[e^{-1.5 |t|} ]$.\\

Similarly, $E_0(t) e^{-\sigma t}$ has an asymptotic \textbf{exponential} fall-off rate of $ o[e^{- |t|} ] $ (using $e^{-(1.5-  \sigma) |t| } $) and $E_0(t) e^{-2 \sigma t}$ has an asymptotic \textbf{exponential} fall-off rate of $ o[e^{-0.5 |t|} ] $ (using $e^{-(1.5- 2 \sigma) |t| } $), for $0 \leq |\sigma| < \frac{1}{2}$ .\\

The above results which show \textbf{exponential} fall-off rates for above mentioned functions, continue to hold, as $|t|$ increases to a larger and larger finite value, without bounds.

.\\ \tocless\subsection{\label{sec:appendix_C_5z} \textbf{ Absolutely integrable functions } \protect\\  \lowercase{} }

We see that a real function $y(t)$ which is finite for all $t$ and has an asymptotic falloff rate of $O[\frac{1}{t^{2}}]$ is an absolutely integrable function, given that $\int_{-\infty}^{\infty} |y(t)| dt = \int_{-\infty}^{-T} |y(t)| dt  + \int_{-T}^{T} |y(t)| dt + \int_{T}^{\infty} |y(t)| dt  $ is finite, for non-zero and finite $T$, because when we integrate the integrand $|y(t)|$ with order $O[\frac{1}{t^{2}}]$ , we get the result $O[\frac{1}{t}]$, which is finite at the limit $t=\pm T$ and the result $O[\frac{1}{t}]$ is zero at the limit $t \to \pm \infty$. If $y(t)$ has an exponential asymptotic falloff rate, when we integrate the integrand $|y(t)|$ with order $O[e^{-A |t|}]$ for real $A>0$, we get the result $O[\frac{1}{A} e^{-A |t|}]$, which is finite at the limit $t=\pm T$ and the result is zero at the limit $t \to \pm \infty$ and hence $y(t)$  is an absolutely integrable function. 

.\\ \tocless\subsection{\label{sec:appendix_C_5a} \textbf{ $E_0(t) > 0$ for $-\infty < t < \infty$ } \protect\\  \lowercase{} }
 
For $0 \leq t < \infty$, we can show that $E_0(t) = \sum_{n=1}^{\infty}  f(t, n) >0$ where \\ $  f(t, n) =  [ 4 \pi^{2} n^{4} e^{4t}    - 6 \pi n^{2}   e^{2t} ]  e^{- \pi n^{2} e^{2t}} e^{\frac{t}{2}} =  2 \pi n^{2}   e^{2t} [ 2 \pi n^{2} e^{2t}    - 3 ]  e^{- \pi n^{2} e^{2t}} e^{\frac{t}{2}} $ as follows.\\

The sum is positive because each summand $f(t, n)$ is positive for finite $n$, and each summand is positive because the term $ 2 \pi n^{2} e^{2t}    - 3  > 0$ for all $t \geq 0$ and $n \geq 1$, given that $\pi > 3$ and $2 \pi n^{2}   e^{2t} e^{- \pi n^{2} e^{2t}} e^{\frac{t}{2}} > 0 $ for $0 \leq t < \infty$ and finite $n \geq 1$.(\textbf{Result C.7.1}) \\

For $t=0$ and $n=1$, we see that $f(0, 1) =   2 \pi [ 2 \pi - 3 ] e^{- \pi  } > 0$.\\

For $t=0$ and for \textbf{each finite} $n \geq 1$, we see that $f(0, n)  = 2 \pi n^{2}   [ 2 \pi n^{2} - 3 ] e^{- \pi   n^{2} }  > 0$.\\

For $0 < t < \infty$ and for \textbf{each finite} $n \geq 1$, we see that $f(t, n) =  2 \pi n^{2}  e^{2t}   [ 2 \pi n^{2}  e^{2t}  - 3 ]  e^{- \pi n^{2} e^{2t}} e^{\frac{t}{2}}  > 0$, using Result C.7.1.\\

As $n \to \infty$, $f(t, n)$ goes to zero, for $0 \leq t < \infty$ due to the term $e^{- \pi n^{2} e^{2t}} $. We do summation over $n$ and see that the sum of the terms  $\sum_{n=1}^{\infty}  f(t, n) > 0$, for $0 \leq t < \infty$. \\

Hence $E_0(t) = \sum_{n=1}^{\infty}  f(t, n) > 0$ for $0 \leq t < \infty$.\\

Given that $\xi(\frac{1}{2} + i \omega)= E_{0\omega}(\omega)$ is an entire function in the whole of s-plane, it is finite for real $\omega$ and also for $\omega=0$. Hence $E_{0\omega}(0) = \int_{-\infty}^{\infty} E_0(t) dt$ is finite. We see that $E_0(t)$ is an analytic function for real $t$ (Details in Section~\ref{sec:Section_1_1}). Hence $E_0(t) = \sum_{n=1}^{\infty}  f(t, n) > 0$ is finite for $0 \leq t < \infty$.\\

Given that $E_0(t)=E_0(-t)$(Details in ~\ref{sec:appendix_C_7}), we see that $E_0(t) > 0$ and finite for all $-\infty < t < \infty$.

.\\ \tocless\subsection{\label{sec:appendix_C_7}  \textbf{$E_0(t)$ is real and even } \protect\\  \lowercase{} }

We see that  $\xi(\frac{1}{2} + i \omega) = E_{0\omega}(\omega) = E_{0\omega}(-\omega)$ (\textbf{Result C.8.1}) because $\xi(s)=\xi(1-s)$ \href{https://www.ocf.berkeley.edu/~araman/files/math_z/Ellison_p147-152.pdf#page=6}{(link)}  and hence $\xi(\frac{1}{2} + i \omega)=\xi(\frac{1}{2} - i \omega)$ when evaluated at $s = \frac{1}{2} + i \omega$.\\

We take the Inverse Fourier transform of $E_{0\omega}(\omega)$ and use $ E_{0\omega}(\omega) = E_{0\omega}(-\omega)$ from Result C.8.1 in the first line in Eq.~\ref{appendix_C_7_a} and then substitute $\omega=-\omega^{'}$ in the second line in Eq.~\ref{appendix_C_7_a}, as follows.
\begin{align}\label{appendix_C_7_a}   
E_0(t)= \frac{1}{2 \pi} \int_{-\infty}^{\infty} E_{0\omega}(\omega) e^{i \omega t} d\omega =
\frac{1}{2 \pi} \int_{-\infty}^{\infty} E_{0\omega}(-\omega) e^{i \omega t} d\omega \notag\\
=  \frac{1}{2 \pi} \int_{-\infty}^{\infty} E_{0\omega}(\omega^{'}) e^{-i \omega^{'} t} d\omega^{'} = E_0(-t)
\end{align}
We see that $E_0(t)$ in Eq.~\ref{sec_intro_eq_1} is real and $E_0(t)$ in Eq.~\ref{appendix_C_7_a} is even and hence we have derived the result that $E_0(t)$ is a \textbf{real and even} function of variable $t$.

.\\ \tocless\section{\label{sec:appendix_I} Properties of Fourier Transforms Part 1 \protect \\ \lowercase{} }
.\\ \tocless\subsection{\label{sec:appendix_I_2} \textbf{Fourier transform of Real g(t)} \protect   \lowercase{} }

In this section, we show that the Fourier transform of a \textbf{real} function $g(t)$, given by $G(\omega) =  G_R(\omega) + i G_I(\omega)$ has the properties given by $ G_R(-\omega) = G_R(\omega) $ and  $ G_I(-\omega) = -G_I(\omega)$. We use the fact that $g(t)$ is real and $ \cos{(\omega t)} $ is an \textbf{even} function of $\omega$ and  $ \sin{(\omega t)} $ is an \textbf{odd} function of $\omega$ below.
\begin{align}\label{sec_C_2_eq_1}   
G(\omega)=  \int_{-\infty}^{\infty} g(t) e^{-i \omega t} dt = G_R(\omega) + i G_I(\omega) \notag\\
G_R(\omega)=  \int_{-\infty}^{\infty} g(t) \cos{(\omega t) } dt = G_R(-\omega) \notag\\
G_I(\omega)=  - \int_{-\infty}^{\infty} g(t) \sin{(\omega t) } dt = - G_I(-\omega) 
\end{align}

.\\ \tocless\subsection{\label{sec:appendix_I_3} \textbf{Even part of g(t) corresponds to real part of Fourier transform $G(\omega)$} \protect\\  \lowercase{} }

In this section, we take the \textbf{even part} of real function $g(t)$, given by $g_{even}(t)=\frac{1}{2} [g(t)+g(-t) ] $ and show that its Fourier transform is given by the \textbf{real part} of  $G(\omega)$. 
\begin{eqnarray*}\label{sec_I_3_eq_1}   
G(\omega)=  \int_{-\infty}^{\infty} g(t) e^{-i \omega t} dt = G_R(\omega) + i G_I(\omega)\\
 \int_{-\infty}^{\infty} g_{even}(t) e^{-i \omega t} dt = \int_{-\infty}^{\infty} \frac{1}{2} [g(t)+g(-t) ] e^{-i\omega t} dt =\frac{G(\omega)}{2} + \frac{1}{2} \int_{-\infty}^{\infty}  g(-t)  e^{-i\omega t} dt
\end{eqnarray*}
\begin{align} \end{align}
We substitute $t=-t$ in the second integral in the last line of Eq.~\ref{sec_I_3_eq_1}. We use the fact that $ G_R(-\omega) = G_R(\omega) $ and  $ G_I(-\omega) = -G_I(\omega) $ for a real function $g(t)$. (Details in ~\ref{sec:appendix_I_2})
\begin{eqnarray*}\label{sec_I_3_eq_2}   
 \int_{-\infty}^{\infty} g_{even}(t) e^{-i \omega t} dt = \frac{G(\omega)}{2} + \frac{1}{2} \int_{-\infty}^{\infty}  g(t)  e^{i\omega t} dt=  \frac{G(\omega)}{2} +  \frac{G(-\omega)}{2}\\
 = \frac{1}{2}  [  G_R(\omega) + i  G_I(\omega) +  G_R(-\omega) + i  G_I(-\omega)] =  \frac{1}{2}  [  G_R(\omega) + i  G_I(\omega) +  G_R(\omega) - i  G_I(\omega)]= G_R(\omega)
\end{eqnarray*}
\begin{align} \end{align}

\clearpage
\tocless\section{\label{sec:Section_2a} \textbf{Dirichlet Eta function}  \protect\\  \lowercase{} }

We use the analytic continuation of Riemann's zeta function given by $\zeta(s)=\displaystyle\frac{\eta(s)}{1-2^{1-s}}$ where $\zeta(s)= \displaystyle\sum\limits_{n=1}^{\infty} \frac{1}{n^{s}}$ diverges for $Re[s] \leq 1$ and $\eta(s)= \displaystyle\sum\limits_{n=1}^{\infty} (-1)^{n-1} \frac{1}{n^{s}}$ is Dirichlet Eta function which converges for $Re[s] >0$. (\href{https://en.wikipedia.org/wiki/Dirichlet_eta_function}{link} and  \href{https://www.ocf.berkeley.edu/~araman/files/math_z/Titchmarsh_pp16_17.pdf}{Titchmarsh pp16-17}) \\

We see that if $\eta(s)$ has a zero in the critical strip, then $\zeta(s)$ also has a zero at the same location. We evaluate $A(s) =  \Gamma(\frac{s}{2}) \eta(s)$ at $s = \frac{1}{2} + \sigma + i \omega$ in Eq.~\ref{sec_2a_eq_1} for $0 < \sigma < \frac{1}{2}$ and compute its inverse Fourier Transform $a(t)$ in Eq.~\ref{sec_2a_eq_4}. \\

We assume that $\eta(\frac{1}{2} + \sigma + i \omega)$ has a zero at $\omega=\omega_0$ in the critical strip (\textbf{Statement A}) and show that the Fourier transform of the function $E_p(t)= E_0(t)  e^{-\sigma t}$ \textbf{also} has a zero at $\omega=\omega_0$, \textbf{if }Statement A  is true, and then prove that this leads to a \textbf{contradiction} for $0 < |\sigma| < \frac{1}{2}$ in ~\ref{sec:Section_2a_3}, where $E_0(t)= \displaystyle\sum\limits_{n=1}^{\infty} (-1)^{n-1}  ( e^{-\pi \frac{n^{2}}{4} e^{ -2 t} } - e^{-\pi n^{2} e^{ -2 t} } ) e^{ -\frac{t}{2} } $ which is derived using $a(t)$ in ~\ref{sec:Sec_A_1_7a}.

.\\  \tocless\subsection{\label{sec:Sec_A_1_7}\textbf{ Analytic continuation of  Riemann Zeta function derived from Dirichlet Eta function} \protect\\  \lowercase{} }

We consider Riemann's Xi function $ \xi(s)$, where $s=\frac{1}{2} + \sigma + i\omega$. Using the functional equation of Riemann's zeta function given by $\zeta(s)= \zeta(1-s) \Gamma(1-s) \sin{(\frac{s \pi}{2})} \pi^{(s-1)} 2^s$, we get $\xi(s)=\xi(1-s)$. \href{https://www.ocf.berkeley.edu/~araman/files/math_z/Titchmarsh_pp16_17.pdf}{Titchmarsh pp16-17})  Using $\zeta(s)=\frac{\eta(s)}{1-2^{1-s}}$, we write as follows.

\begin{eqnarray*}\label{sec:sec_a_1_7_eq_0}  
\xi(s)= \zeta(s) \Gamma(\frac{s}{2}) \pi^{\frac{-s}{2}} \frac{s (s-1)}{2} = \xi(1-s) \\
\xi(s)= \frac{\eta(s)}{1-2^{1-s}} \Gamma(\frac{s}{2}) \pi^{\frac{-s}{2}} \frac{s (s-1)}{2} 
\end{eqnarray*}
\begin{equation} \end{equation}

We define a related analytic continuation $E(s)$ as follows. Given $\xi(s)= \xi(1-s)$, we see that $E(s)=E(1-s)$ is analytic in the region $0 < Re[s] < 1$ and has simple poles at $s=0$ and $s=1$.  

\begin{eqnarray*}  \label{sec:sec_a_1_7_eq_1}   
E(s)=\frac{\xi(s) (1-2^{1-s}) (2^{s}-1)}{s (s-1)} \\
E(1-s) = \frac{\xi(1-s) (1-2^{s}) (2^{1-s}-1)}{(1-s) (-s)} =\frac{\xi(s) (2^{s}-1) (1-2^{1-s})}{(s-1) (s)}  = E(s)
\end{eqnarray*}
\begin{equation}  \end{equation}

We substitute $\xi(s)$ from Eq.~\ref{sec:sec_a_1_7_eq_0} and $\zeta(s)=\displaystyle\frac{\eta(s)}{1-2^{1-s}}$  in Eq.~\ref{sec:sec_a_1_7_eq_1} and cancel the common terms $s (s-1)$ and $(1-2^{1-s})$ as follows.

\begin{eqnarray*}\label{sec:sec_a_1_7_eq_2}   
E(s)=\frac{\eta(s)}{1-2^{1-s}} \Gamma(\frac{s}{2}) \pi^{\frac{-s}{2}} \frac{s (s-1)}{2} \frac{(1-2^{1-s}) (2^{s}-1)}{s (s-1)}  \\
E(s)=\frac{\eta(s)}{1-2^{1-s}} \Gamma(\frac{s}{2}) \pi^{\frac{-s}{2}} \frac{1}{2} (1-2^{1-s}) (2^{s}-1) \\ 
E(s) = \eta(s) \Gamma(\frac{s}{2}) \frac{\pi^{\frac{-s}{2}}}{2} (2^{s} - 1) 
\end{eqnarray*}
\begin{equation} \end{equation}

We evaluate $E(s)$ at $s= \frac{1}{2} + \sigma + i\omega$ and use $K^{i \omega} = e^{i\omega \log(K)}$ as follows. 

\begin{eqnarray*}\label{sec:sec_a_1_7_eq_3}   
E(\frac{1}{2} + \sigma + i\omega)= E_{p\omega}(\omega) = \eta(\frac{1}{2} + \sigma + i\omega) \Gamma(\frac{\frac{1}{2} + \sigma + i\omega}{2}) \frac{\pi^{\frac{-(\frac{1}{2} + \sigma)}{2}}}{2}  e^{\frac{-i\omega}{2} \log(\pi)} (2^{\frac{1}{2} + \sigma}e^{i\omega \log(2)}-1)
\end{eqnarray*}
\begin{equation} \end{equation}

We define $A_{\omega}(\omega)= \eta(\frac{1}{2} + \sigma + i\omega) \Gamma(\frac{\frac{1}{2} + \sigma + i\omega}{2}) $, and we can rearrange the terms as follows.

\begin{equation}  \label{sec:sec_a_1_7_eq_3a} 
E_{p\omega}(\omega) =  A_{\omega}(\omega)   \frac{\pi^{\frac{-(\frac{1}{2} + \sigma)}{2}}}{2} e^{\frac{-i\omega}{2} \log(\pi)} (2^{\frac{1}{2} + \sigma}e^{i\omega \log(2)} -1)   
\end{equation}

We define $a(t)$ as the Inverse Fourier Transform of  $A_{\omega}(\omega)$.
We compute the Inverse Fourier Transform of $E_{p\omega}(\omega)$ given by $E_{p}(t)$ as follows, using \href{https://lpsa.swarthmore.edu/Fourier/Xforms/FXProps.html}{time shifting property}. 
 
\begin{eqnarray*}\label{sec:sec_a_1_7_eq_4}   
E_{p}(t) =  \frac{\pi^{\frac{-(\frac{1}{2} + \sigma)}{2}}}{2}  [ 2^{\frac{1}{2} + \sigma} a(t-\frac{\log(\pi)}{2}+\log(2)) -  a(t-\frac{\log(\pi)}{2})  ]  
\end{eqnarray*}
\begin{equation} \end{equation}

 \tocless\subsection{\label{sec:Sec_A_1_7a} \textbf{Derivation of $a(t)$ and $E_p(t)$} \protect  \lowercase{} \\}

 We start with the gamma function $\Gamma(\frac{s}{2}) = \int_{0}^{\infty} y^{\frac{s}{2}-1} e^{-y} dy$.  We evaluate $A(s) = \Gamma(\frac{s}{2}) \eta(s)$ at $s=\frac{1}{2} + \sigma + i\omega$  below. We substitute $y =  \pi n^{2} x$ and $dy =   \pi n^{2} dx$ in Eq.~\ref{sec_2a_eq_1} and get $y^{\frac{s}{2}-1} dy = (\pi n^{2})^{\frac{s}{2}-1} x^{\frac{s}{2}-1} \pi n^{2} dx =  \pi^{\frac{s}{2}} n^{s} (\pi n^{2})^{-1} x^{\frac{s}{2}-1 } \pi n^{2} dx = \pi^{\frac{s}{2}} n^{s} x^{\frac{s}{2}-1 } dx $. 

\begin{equation}  \label{sec_2a_eq_1}  
A(s) = \Gamma(\frac{s}{2}) \eta(s)  = \displaystyle\sum\limits_{n=1}^{\infty} (-1)^{n-1} \frac{1}{n^{s}} \int_{0}^{\infty} y^{\frac{s}{2}-1} e^{-y} dy =  \pi^{\frac{s}{2}} \displaystyle\sum\limits_{n=1}^{\infty} (-1)^{n-1} \frac{1}{n^{s}} n^{s} \int_{0}^{\infty}  x^{\frac{s}{2}-1} e^{-\pi n^{2} x } dx 
\end{equation}

For $Re[s]>0$, the gamma function is analytic in the complex plane \href{https://proofwiki.org/wiki/Poles_of_Gamma_Function}{(link)} and $\eta(s)$ converges and hence $|A(s)| = |\Gamma(\frac{s}{2}) \eta(s)|$ converges. \\


We consider $A^{'}(s)$ obtained after interchanging the order of integration and summation  in Eq.~\ref{sec_2a_eq_1} as follows.  We use $ \frac{1}{n^{s}} n^{s} =1$.

\begin{equation} \label{sec_2a_eq_1a}  
A^{'}(s) =   \pi^{\frac{s}{2}} \int_{0}^{\infty}  \displaystyle\sum\limits_{n=1}^{\infty} (-1)^{n-1}   e^{-\pi n^{2} x }  x^{\frac{s}{2}-1} dx  
 \end{equation}

Now we substitute $x = e^{-2 t}$ and $dx = -2 e^{-2 t} dt = -2 x dt$ and write Eq.~\ref{sec_2a_eq_1a} as follows.

\begin{equation}  \label{sec_2a_eq_2}  
A^{'}(s) =  2 \pi^{\frac{s}{2}} \int_{-\infty}^{\infty}  \displaystyle\sum\limits_{n=1}^{\infty} (-1)^{n-1}   e^{-\pi n^{2} e^{-2 t} }  e^{-s t} dt 
\end{equation}

We substitute $s = \frac{1}{2} + \sigma + i \omega$ in Eq.~\ref{sec_2a_eq_2} as follows.

\begin{equation} \label{sec_2a_eq_3}  
A^{'}( \frac{1}{2} + \sigma + i\omega) = A_{\omega}^{'}(\omega) =  2 \pi^{\frac{ \frac{1}{2} + \sigma }{2}} e^{ \frac{i \omega}{2} \log{\pi} } \int_{-\infty}^{\infty}  \displaystyle\sum\limits_{n=1}^{\infty} (-1)^{n-1}   e^{-\pi n^{2} e^{ -2 t} }  e^{ -\frac{t}{2} } e^{-\sigma t} e^{-i \omega t} dt 
 \end{equation}

The integrand in Eq.~\ref{sec_2a_eq_3} is absolutely integrable given asymptotic exponential fall-off rate. (~\ref{sec:appendix_C_5c}). Hence we can interchange the order of integration and summation in
Eq.~\ref{sec_2a_eq_3} and Eq.~\ref{sec_2a_eq_1a} for $A^{'}(s)$ and get $A(s)$ in Eq.~\ref{sec_2a_eq_1}. Hence $A(s)=A^{'}(s)$ and we can interchange the order of integration and summation  in Eq.~\ref{sec_2a_eq_1} using Fubini's theorem, to obtain Eq.~\ref{sec_2a_eq_1a} and the interchange is justified.\\

We see that the inverse Fourier transform of $A_{\omega}(\omega)=A_{\omega}^{'}(\omega)$ is given by $a(t)$ as follows, using the \href{https://lpsa.swarthmore.edu/Fourier/Xforms/FXProps.html}{time shifting property}.

\begin{equation}\label{sec_2a_eq_4}
a(t) = a_0(t+ \frac{ \log{\pi}}{2}), \quad a_0(t)= 2 \pi^{\frac{1}{4} + \frac{ \sigma }{2}}  \displaystyle\sum\limits_{n=1}^{\infty} (-1)^{n-1}   e^{-\pi n^{2} e^{ -2 t} }  e^{ -\frac{t}{2} } e^{-\sigma t}
\end{equation}

We know that $\Gamma(\frac{s}{2})$ does not have zeros for any value of $s$ \href{https://www.ocf.berkeley.edu/~araman/files/math_z/Ellison_p147-152.pdf#page=4}{(link)} and the gamma function is analytic in the complex plane for $Re[s]>0$ \href{https://proofwiki.org/wiki/Poles_of_Gamma_Function}{(link)}. If $\eta(s)$ has a zero at $\omega=\omega_0$ in the critical strip to satisfy Statement A, then  $A( \frac{1}{2} + \sigma + i \omega) $ in Eq.~\ref{sec_2a_eq_1} has a zero at $\omega=\omega_0$ and the Fourier transform of $a(t) $ given by $A_{\omega}(\omega)$ in Eq.~\ref{sec_2a_eq_3} has a zero at $\omega=\omega_0$ (\textbf{Result E.0})\\

Now we substitute $a(t)$ in Eq.~\ref{sec_2a_eq_4} in Eq.~\ref{sec:sec_a_1_7_eq_4} copied below and cancel the common terms $\frac{\log(\pi)}{2}$ and $2 \pi^{\frac{1}{4} + \frac{ \sigma }{2}}$ as follows. We use $2^{\frac{1}{2} + \sigma}  2^{-(\frac{1}{2} + \sigma)} = 1$ in the first term in $E_p(t)$ below.

\begin{eqnarray*}\label{sec:sec_a_1_7_eq_5}  
E_{p}(t) =  \frac{\pi^{\frac{-(\frac{1}{2} + \sigma)}{2}}}{2}  [ 2^{\frac{1}{2} + \sigma} a(t-\frac{\log(\pi)}{2}+\log(2)) -  a(t-\frac{\log(\pi)}{2})  ]   \\
E_{p}(t) =  \frac{\pi^{-(\frac{1}{4} + \frac{ \sigma }{2})}}{2}  [ 2^{\frac{1}{2} + \sigma} a_0(t-\frac{\log(\pi)}{2}+\frac{\log(\pi)}{2}+\log(2)) -  a_0(t-\frac{\log(\pi)}{2}+\frac{\log(\pi)}{2}  ]  \\
E_{p}(t) = \frac{\pi^{-(\frac{1}{4} + \frac{ \sigma }{2})}}{2}   [ 2^{\frac{1}{2} + \sigma} a_0(t+\log(2)) -  a_0(t) ] , \quad 
a_0(t+\log(2)) = 2 * 2^{-(\frac{1}{2} + \sigma)} \pi^{\frac{1}{4} + \frac{ \sigma }{2}}  \displaystyle\sum\limits_{n=1}^{\infty} (-1)^{n-1}   e^{-\pi \frac{n^{2}}{4} e^{ -2 t} }  e^{ -\frac{t}{2} } e^{-\sigma t} \\ 
E_{p}(t) =  \displaystyle\sum\limits_{n=1}^{\infty} (-1)^{n-1}   e^{-\pi \frac{n^{2}}{4} e^{ -2 t} }  e^{ -\frac{t}{2} } e^{-\sigma t} -  \displaystyle\sum\limits_{n=1}^{\infty} (-1)^{n-1}   e^{-\pi n^{2} e^{ -2 t} }  e^{ -\frac{t}{2} } e^{-\sigma t}\\
E_{p}(t) = E_0(t) e^{-\sigma t} , \quad E_0(t) =  \displaystyle\sum\limits_{n=1}^{\infty} (-1)^{n-1}  ( e^{-\pi \frac{n^{2}}{4} e^{ -2 t} } - e^{-\pi n^{2} e^{ -2 t} } ) e^{ -\frac{t}{2} } 
\end{eqnarray*}
\begin{equation} \end{equation}


.\\

We see that $E_{0}(t)$ is the inverse Fourier transform of $E(\frac{1}{2} +  i\omega)$ (set $\sigma=0$ in Eq.~\ref{sec:sec_a_1_7_eq_3} and Eq.~\ref{sec:sec_a_1_7_eq_4}) and it obeys $E_0(t)=E_0(-t)$ given that $E(s)=E(1-s)$ using Eq.~\ref{sec:sec_a_1_7_eq_1} (We use the result in ~\ref{sec:appendix_C_7} with $\xi(s)$ replaced by $E(s)$). (\textbf{Result E.1})\\ 

Using Eq.~\ref{sec:sec_a_1_7_eq_3} and Eq.~\ref{sec:sec_a_1_7_eq_3a}, we have derived the analytic continuation of Riemann's zeta function derived from Dirichlet Eta function given by $E_{p\omega}(\omega)=  \eta(\frac{1}{2} + \sigma + i\omega) B(\omega)$  where \\$B(\omega)=\Gamma(\frac{\frac{1}{2} + \sigma + i\omega}{2})  \frac{\pi^{\frac{-(\frac{1}{2} + \sigma)}{2}}}{2} e^{\frac{-i\omega}{2} \log(\pi)} (2^{\frac{1}{2} + \sigma}e^{i\omega \log(2)} -1) $.\\

We see that, if $\eta(\frac{1}{2} + \sigma + i \omega)$ has a zero at $\omega=\omega_0$ in the critical strip, then the Fourier transform of the function $E_p(t)= E_0(t)  e^{-\sigma t}$ given by $E_{p\omega}(\omega)$  \textbf{also} has a zero at $\omega=\omega_0$, where $E_0(t)= \displaystyle\sum\limits_{n=1}^{\infty} (-1)^{n-1}  ( e^{-\pi \frac{n^{2}}{4} e^{ -2 t} } - e^{-\pi n^{2} e^{ -2 t} } ) e^{ -\frac{t}{2} } $ .

.\\  \tocless\subsection{\label{sec:Section_2a_1} \textbf{$E_0(t) > 0$ for $-\infty < t < \infty$}  \protect\\  \lowercase{} }

It is shown in this section that $E_0(t) > 0$ for $-\infty < t < \infty$. We use Eq.~\ref{sec:sec_a_1_7_eq_5}  below and take the term $e^{-\pi \frac{n^{2}}{4} e^{2 t} } e^{\frac{t}{2} }$ out of the brackets in Eq.~\ref{sec_2a_1_eq_1} for $E_0(-t)$ and use $(n+1)^{2}= n^{2}+2n+1$ and rearrange the terms in the last line below.

\begin{eqnarray*}\label{sec_2a_1_eq_1}  
E_0(-t)= \displaystyle\sum\limits_{n=1}^{\infty} (-1)^{n-1}  ( e^{-\pi \frac{n^{2}}{4} e^{ 2 t} } - e^{-\pi n^{2} e^{2 t} } ) e^{ \frac{t}{2} } \\
E_0(-t)= \displaystyle\sum\limits_{n=odd}^{\infty}  ( e^{-\pi \frac{n^{2}}{4} e^{ 2 t} } - e^{-\pi n^{2} e^{2 t} } - e^{-\pi \frac{(n+1)^{2}}{4} e^{2 t} } + e^{-\pi (n+1)^{2} e^{2 t} } ) e^{ \frac{t}{2} } \\
E_0(-t)= \displaystyle\sum\limits_{n=odd}^{\infty} e^{-\pi \frac{n^{2}}{4} e^{2 t} } e^{\frac{t}{2} } ( 1 - e^{-\pi \frac{3 n^{2}}{4} e^{2 t} } - e^{-\pi \frac{(2 n+1)}{4} e^{2 t} } + e^{-\pi \frac{3 n^{2}}{4} e^{2 t} } e^{-\pi (2 n+1) e^{2 t} } ) 
\end{eqnarray*}
\begin{equation} \end{equation}

We compute the \textbf{minimum} value of $E_0(-t)$ in Eq.~\ref{sec_2a_1_eq_1} for $0 \leq t < \infty$, by computing the minimum value of positive terms and maximum value of absolute value of negative terms. We ignore the last term $e^{-\pi \frac{3 n^{2}}{4} e^{2 t} } e^{-\pi (2 n+1) e^{2 t} } > 0$ because we want the minimum value of $E_0(-t)$.\\

The minimum value of the first term inside brackets in Eq.~\ref{sec_2a_1_eq_1} is $A_1=1$. The maximum value of the absolute value of the second term inside brackets $e^{-\pi \frac{3 n^{2}}{4} e^{2 t} }$ occurs at $n=1$ and $t=0$, given by $A_2=e^{-\pi \frac{3}{4}}$. The maximum value of the absolute value of the third term $ e^{-\pi \frac{(2 n+1)}{4} e^{2 t} } $ occurs at $n=1$ and $t=0$, given by $A_3=e^{-\pi \frac{3}{4}}$. Hence the minimum value of the terms inside the brackets is given by $A_1 - A_2 - A_3 = 1 - 2 e^{-\pi \frac{3}{4}} = 0.8104 > 0$ for all $n$ and hence $E_0(-t) > 0$ for $0 \leq t < \infty$.\\




Given that $E_0(t)=E_0(-t)$ (We use the result in ~\ref{sec:appendix_C_7} with $\xi(s)$ replaced by $E(s)$ in Eq.~\ref{sec:sec_a_1_7_eq_1}), we see that $E_0(t)>0$ for $-\infty < t < \infty$.

.\\  \tocless\subsection{\label{sec:Section_2a_2} \textbf{$E_0(t)$ is strictly decreasing for $t > 0$}  \protect\\  \lowercase{} }

We show that $E_0(t)=E_0(-t)$ is strictly decreasing for $t > 0$, by showing that $\frac{d E_0(-t)}{dt} < 0$ for $ 0 < t < \infty$. We set $y=\pi e^{2t}$ in $E_0(-t)$ in the second line in Eq.~\ref{sec_2a_1_eq_1}  and then take the first derivative of $E_0(y)$ as follows.  We see that  $ \frac{dy}{dt} = 2 \pi e^{2t}= 2y $  and $\frac{d E_0(-t)}{dt} = \frac{d E_0(-t)}{dy} \frac{dy}{dt} =  \frac{d E_0(y)}{dy} 2 y$ and hence we will show that $\frac{d E_0(y)}{dy}  < 0$ for $\pi < y < \infty$.


\begin{eqnarray*}\label{sec_2a_2_eq_1} 
E_0(-t)= \displaystyle\sum\limits_{n=odd}^{\infty} ( e^{-\pi \frac{n^{2}}{4} e^{ 2 t} } - e^{-\pi n^{2} e^{2 t} } - e^{-\pi \frac{(n+1)^{2}}{4} e^{2 t} } + e^{-\pi (n+1)^{2} e^{2 t} } ) e^{ \frac{t}{2} } \\ 
E_0(y)= (\pi)^{-\frac{1}{4}} \displaystyle\sum\limits_{n=odd}^{\infty}  e^{- \frac{n^{2}}{4} y } y^{\frac{1}{4}} - e^{- n^{2} y }  y^{\frac{1}{4}} - e^{- \frac{(n+1)^{2}}{4} y }  y^{\frac{1}{4}} + e^{- (n+1)^{2} y}  y^{\frac{1}{4}} \\
\frac{d E_0(y)}{dy} = (\pi)^{-\frac{1}{4}} \displaystyle\sum\limits_{n=odd}^{\infty}  e^{- \frac{n^{2}}{4} y } y^{\frac{1}{4}} (\frac{1}{4y} - \frac{n^{2}}{4}   )  - e^{- n^{2} y }  y^{\frac{1}{4}} (\frac{1}{4y} - n^{2}   )  \\
- e^{- \frac{(n+1)^{2}}{4} y }  y^{\frac{1}{4}} (\frac{1}{4y} - \frac{(n+1)^{2}}{4}   )  + e^{- (n+1)^{2} y}  y^{\frac{1}{4}} (\frac{1}{4y} - (n+1)^{2}   )  
\end{eqnarray*}
\begin{equation} \end{equation}


We take the common term $ e^{- \frac{n^{2}}{4} y } y^{\frac{1}{4}} $ out and use $(n+1)^{2}= n^{2}+2n+1$ and rearrange the terms in Eq.~\ref{sec_2a_2_eq_1}  as follows.

\begin{eqnarray*}\label{sec_2a_2_eq_2a}  
\frac{d E_0(y)}{dy} = (\pi)^{-\frac{1}{4}} \displaystyle\sum\limits_{n=odd}^{\infty}  e^{- \frac{n^{2}}{4} y } y^{\frac{1}{4}} [ (\frac{1}{4y} - \frac{n^{2}}{4}   )  - e^{- \frac{3 n^{2}}{4} y }  (\frac{1}{4y} - n^{2}   )  \\
- e^{- \frac{(2 n+1)}{4} y }  (\frac{1}{4y} - \frac{(n+1)^{2}}{4}   )  + e^{- \frac{3 n^{2}}{4} y }  e^{- (2 n+1) y}   (\frac{1}{4y} - (n+1)^{2}   )  ]
\end{eqnarray*}
\begin{equation} \end{equation}

We compute the \textbf{maximum} value of $\frac{d E_0(y)}{dy} $ in Eq.~\ref{sec_2a_2_eq_2a}, by computing the maximum value of positive terms and minimum value of absolute value of negative terms. We ignore the negative terms inside the brackets $- e^{- \frac{3 n^{2}}{4} y } \frac{1}{4y}$, $- e^{- \frac{(2 n+1)}{4} y } \frac{1}{4y}$  and $ - (n+1)^{2}  e^{- \frac{3 n^{2}}{4} y }  e^{- (2 n+1) y}  $ because we want the maximum value of $\frac{d E_0(y)}{dy}$ in the interval  $\pi < y < \infty$.

\begin{eqnarray*}\label{sec_2a_2_eq_2}  
\frac{d E_0(y)}{dy} < (\pi)^{-\frac{1}{4}} \displaystyle\sum\limits_{n=odd}^{\infty}  e^{- \frac{n^{2}}{4} y } y^{\frac{1}{4}} [ (\frac{1}{4y} - \frac{n^{2}}{4}   )  + e^{- \frac{3 n^{2}}{4} y }  n^{2}     \\
+ e^{- \frac{(2 n+1)}{4} y }  \frac{(n+1)^{2}}{4}     + e^{- \frac{3 n^{2}}{4} y }  e^{- (2 n+1) y}   \frac{1}{4y}   ]
\end{eqnarray*}
\begin{equation} \end{equation}

We see that $y=\pi e^{2t}$ is in the range $y=[\pi, \infty)$ for $0 \leq t < \infty$, and in the range $y=[\pi, y_a)$ for $0 \leq t < t_a=0.1$, where $y_a=\pi e^{2 t_a}=3.8371$. \\


$\bullet$ It is shown in Section~\ref{sec:Section_2a_2a} that $\frac{d E_0(y)}{dy} <0$ for $y_a \leq y < \infty$ for  $y_a= 3.8371$.\\

 
$\bullet$ It is shown in Section~\ref{sec:Section_2a_2c} that  $\frac{d^2 E_0(y)}{dy^2} <0$ for $\pi \leq y < y_a $ and hence  $\frac{d E_0(y)}{dy} <0$ for  $\pi < y < y_a $, given that $\frac{d E_0(y)}{dy} = 0$ at $y=\pi$ corresponding to $t=0$, using $y= \pi   e^{2t}$. We use $\frac{d E_0(-t)}{dt} =  \frac{d E_0(y)}{dy} 2 y$ from the line below and $\frac{d E_0(-t)}{dt} = 0$ at $t=0$, using $E_0(t)=E_0(-t)$. (\textbf{Result E.5.a})\\

$\bullet$ Hence  $\frac{d E_0(y)}{dy} <0$ for $\pi < y < \infty$. Given $y= \pi   e^{2t} $ and $ \frac{dy}{dt} = 2 \pi e^{2t}= 2y $  and $\frac{d E_0(-t)}{dt} = \frac{d E_0(-t)}{dy} \frac{dy}{dt} =  \frac{d E_0(y)}{dy} 2 y$, we see that  $\frac{d E_0(-t)}{dt}  < 0$ for $t > 0$. Hence $E_0(t)=E_0(-t)$ is strictly decreasing for $t > 0$.

.\\  \tocless\subsubsection{\label{sec:Section_2a_2a} \textbf{$\frac{d E_0(y)}{dy} <0$ for $y_a \leq y < \infty$  for $y_a=3.8371$}  \protect\\  \lowercase{} }

We see that the \textbf{maximum} value of the \textbf{first term} inside brackets $(\frac{1}{4y} - \frac{n^{2}}{4}  )$ in Eq.~\ref{sec_2a_2_eq_2} occurs at $n=1$ and $y=y_{a}=3.8371$ given by  $D_1 = \frac{1}{4 y_{a}} - \frac{1}{4} =  \frac{1}{4 * 3.8371} - \frac{1}{4} = -0.1848$.\\

We consider the \textbf{second term} inside brackets in Eq.~\ref{sec_2a_2_eq_2} given by $I(y,n)=  n^{2} e^{- \frac{3 n^{2}}{4} y }  $. It is a \textbf{strictly decreasing} function in the region $y_a \leq y < \infty$, with \textbf{maximum} value at $y=y_a$, for each $n$. \\

We set $y=y_{a}=3.8371$ and compute $\frac{d I(y_a,n)}{dn} = e^{- \frac{3 n^{2}}{4} y_a }  [ 2 n  + n^{2} ( - \frac{6 n y_a}{4} ) ]$ which has an inflection point at $ 2 n  + n^{2} ( - \frac{6 n y_a}{4} ) = 0$. We cancel common term $n$ and get $ 2 + n^{2}  (-\frac{6 y_a}{4})= 0$ which has roots at $n^{2}=\frac{4}{3 y_a}$ given by $n = \pm 0.5895$. Hence we choose $n=0.5895$ as a positive solution and the nearest positive integer is $n=1$, where $I(y_a,n)$ has a \textbf{maximum} value for all positive integer $n$. Given that $I(y_a,n)> 0$ for all finite $n$ and goes to zero as $n \to \infty $ due to the term $e^{- \frac{3 n^{2}}{4} y_a } $, and there is only one positive inflection point at $n=0.5895$, this inflection point is  a \textbf{maximum} point and $I(y_a,n)$ is \textbf{strictly decreasing} for $n > 0.5895$. (\textbf{Result E.5.1})\\ 

Hence the \textbf{maximum} value of $I(y,n)$ in the interval $y_a \leq y < \infty$, is at $y=y_a$ and $n=1$  given by $I(y_a,1)=  e^{- \frac{3}{4} y_a } = 0.0563 = D_2$.\\

We consider the \textbf{third term} inside brackets in Eq.~\ref{sec_2a_2_eq_2} given by $J(y,n)=   \frac{(n+1)^{2}}{4}   e^{- \frac{(2 n+1)}{4} y }  $ which is \textbf{strictly decreasing} function in the interval $y_a \leq y < \infty$, with \textbf{maximum} value at $y=y_a$, for each $n$. \\

We set $y=y_{a}=3.8371$ and compute $\frac{d J(y_a,n)}{dn} =   e^{- \frac{(2 n+1)}{4} y_a } [ \frac{2 (n+1)}{4} + \frac{(n+1)^{2}}{4} (- \frac{(2 y_a)}{4} ) ] $ which has an inflection point at $  \frac{2 (n+1)}{4} + \frac{(n+1)^{2}}{4} (- \frac{(2 y_a)}{4} ) = 0$. We cancel common term $\frac{2 (n+1)}{4}$ and get $1 - (n+1) \frac{y_a}{4} = 0$ which has roots at $n+1=\frac{4}{y_a}=1.0424$ given by $n = 0.0424$. Hence  $J(y_a,n)$ is \textbf{strictly decreasing} for $n > 0.0424$ and the nearest positive integer is $n=1$ where $J(y_a,n)$ has a \textbf{maximum} value for all positive integer $n$. Given that $J(y_a,n)> 0$ for all finite $n$ and goes to zero as $n \to \infty $ due to the term $e^{- \frac{(2 n+1)}{4} y_a }$, and there is only one positive inflection point at $n = 0.0424$, this inflection point is  a \textbf{maximum} point. (\textbf{Result E.5.2})\\ 

Hence the \textbf{maximum} value of $J(y,n)$ in the interval $y_a \leq y < \infty$, is at $y=y_a$ and $n=1$  given by $J(y_a,1)=  e^{- \frac{3}{4} y_a } = 0.0563 = D_3$. \\

The fourth term in Eq.~\ref{sec_2a_2_eq_2} given by $e^{- \frac{3 n^{2}}{4} y }  e^{- (2 n+1) y}   \frac{1}{4y}$ has a maximum at $n=1$ and $y=y_a$ given by $e^{- \frac{3 }{4} y_a }  e^{- 3 y_a }   \frac{1}{4 y_a} = 3.6706 * 10^{-8} < 10^{-7} = D_4$. \\

Hence the maximum value of the terms in square bracket in Eq.~\ref{sec_2a_2_eq_2}  for $y_a \leq y < \infty$ and for $n=1$, is given by $D_1+D_2+D_3+D_4=-0.1848 + 0.0563+ 0.0563 + 10^{-7} \approx -0.0722 < 0 $. This  summation is negative for $n>1$, given Result E.5.1 and Result E.5.2 and $D_2+D_3+D_4$ is a smaller positive value and $D_1$ is more negative than the case for $n=1$. Hence $\frac{d E_0(y)}{dy} <0$ for $y_a \leq y < \infty$, given summation of negative terms for each odd $n$ and given that $e^{- \frac{n^{2}}{4} y } y^{\frac{1}{4}}  >0$ for all finite $n$ and $y$.



.\\  \tocless\subsubsection{\label{sec:Section_2a_2c} \textbf{ $\frac{d^2 E_0(y)}{dy^2} <0$ for  $\pi \leq y < y_a $ and hence $\frac{d E_0(y)}{dy} <0$ for  $\pi < y < y_a $  }  \protect\\  \lowercase{} }

We compute the second derivative $\frac{d^2 E_0(-t)}{dt^2}$ from Eq.~\ref{sec_2a_2_eq_1} as follows.\\

We copy $\frac{d E_0(y)}{dy}$ in Eq.~\ref{sec_2a_2_eq_1} as follows.

\begin{eqnarray*}\label{sec_2a_2_eq_3}  
E_0(y)= (\pi)^{-\frac{1}{4}} \displaystyle\sum\limits_{n=odd}^{\infty}  e^{- \frac{n^{2}}{4} y } y^{\frac{1}{4}} - e^{- n^{2} y }  y^{\frac{1}{4}} - e^{- \frac{(n+1)^{2}}{4} y }  y^{\frac{1}{4}} + e^{- (n+1)^{2} y}  y^{\frac{1}{4}} \\
\frac{d E_0(y)}{dy} = (\pi)^{-\frac{1}{4}} \displaystyle\sum\limits_{n=odd}^{\infty}  e^{- \frac{n^{2}}{4} y } y^{\frac{1}{4}} (\frac{1}{4y} - \frac{n^{2}}{4}   )  - e^{- n^{2} y }  y^{\frac{1}{4}} (\frac{1}{4y} - n^{2}   )  \\
- e^{- \frac{(n+1)^{2}}{4} y }  y^{\frac{1}{4}} (\frac{1}{4y} - \frac{(n+1)^{2}}{4}   )  + e^{- (n+1)^{2} y}  y^{\frac{1}{4}} (\frac{1}{4y} - (n+1)^{2}   )  \\
\end{eqnarray*}
\begin{equation} \end{equation}

We compute the second derivative $\frac{d^2 E_0(y)}{dy^2}$ as follows.

\begin{eqnarray*}\label{sec_2a_2_eq_4}  
\frac{d^2 E_0(y)}{dy^2} = (\pi)^{-\frac{1}{4}} \displaystyle\sum\limits_{n=odd}^{\infty}  e^{- \frac{n^{2}}{4} y } y^{\frac{1}{4}} ( -\frac{1}{4y^{2}} + (\frac{1}{4y} - \frac{n^{2}}{4})^{2} )  - e^{- n^{2} y }  y^{\frac{1}{4}} (  -\frac{1}{4y^{2}} + (\frac{1}{4y} - n^{2})^{2} )  \\
- e^{- \frac{(n+1)^{2}}{4} y }  y^{\frac{1}{4}} ( -\frac{1}{4y^{2}} + (\frac{1}{4y} - \frac{(n+1)^{2}}{4})^{2})  + e^{- (n+1)^{2} y}  y^{\frac{1}{4}} ( -\frac{1}{4y^{2}} + (\frac{1}{4y} - (n+1)^{2})^{2})  \\
\end{eqnarray*}
\begin{equation} \end{equation}

We simplify it as follows.

\begin{eqnarray*}\label{sec_2a_2_eq_5}  
\frac{d^2 E_0(y)}{dy^2} = (\pi)^{-\frac{1}{4}} \displaystyle\sum\limits_{n=odd}^{\infty}  e^{- \frac{n^{2}}{4} y } y^{\frac{1}{4}} ( -\frac{1}{4y^{2}} + \frac{1}{16 y^{2}} - \frac{n^{2}}{8y} + \frac{n^{4}}{16} )  \\
- e^{- n^{2} y }  y^{\frac{1}{4}} ( -\frac{1}{4y^{2}} + \frac{1}{16 y^{2}} - \frac{n^{2}}{2y} + n^{4} )  \\
- e^{- \frac{(n+1)^{2}}{4} y }  y^{\frac{1}{4}} ( -\frac{1}{4y^{2}} + \frac{1}{16 y^{2}} - \frac{(n+1)^{2}}{8y} + \frac{(n+1)^{4}}{16} ) \\
  + e^{- (n+1)^{2} y}  y^{\frac{1}{4}} ( -\frac{1}{4y^{2}} + \frac{1}{16 y^{2}} - \frac{(n+1)^{2}}{2y} + (n+1)^{4} ) \\
\end{eqnarray*}
\begin{equation} \end{equation}

We compute the \textbf{maximum} value of $\frac{d^2 E_0(y)}{dy^2}$  with $y=\pi e^{2t}$ in the range $y=[\pi, y_a)$ for $0 \leq t < t_a=0.1$, where $y_a=3.8371$, by computing the maximum value of positive terms and minimum value of absolute value of negative terms. Let the maximum value of $y$ be $y_{max}=y_a=\pi e^{2t_a}$ and minimum value of $y$ be $y_{min}=\pi$ in the interval $y=[\pi, y_a)$.\\


The first term in curved brackets in Eq.~\ref{sec_2a_2_eq_5} is given by $ -\frac{1}{4y^{2}} + \frac{1}{16 y^{2}} - \frac{n^{2}}{8y} + \frac{n^{4}}{16}  =   -\frac{3}{16 y^{2}} - \frac{n^{2}}{8y} + \frac{n^{4}}{16}    $ and the \textbf{maximum value} of the whole first term in the interval $y=[y_{min},y_{max})$ is given by $e^{- \frac{1}{4} y_{min} } (y_{max})^{\frac{1}{4}} \frac{n^{4}}{16} - e^{- \frac{1}{4} y_{max} } (y_{min})^{\frac{1}{4}} (\frac{3}{16 y_{max}^{2}} + \frac{n^{2}}{8y_{max}})$ and similarly we compute the other 3 terms  at $n=1,3,5,7,9$. The \textbf{maximum} value of $\frac{d^2 E_0(y)}{dy^2}$  in Eq.~\ref{sec_2a_2_eq_5} at $n=1,3,5,7,9$ in the interval $y=[y_{min},y_{max})$ is given by $-0.0097$ which is \textbf{negative}. (\textbf{Result E.5.5}) \href{https://www.ocf.berkeley.edu/~araman/files/math_z/test_Phi_inv_all_n.m}{Matlab simulation}) \\

\textbf{Case $ n \geq 11$}: We note that $-\frac{1}{4y^{2}} + \frac{1}{16 y^{2}} =  -\frac{3}{16 y^{2}} $ and ignore the negative terms in Eq.~\ref{sec_2a_2_eq_5} because we are computing the maximum value  of $\frac{d^2 E_0(y)}{dy^2}$ for $ n \geq 11$ given by $[\frac{d^2 E_0(y)}{dy^2}]_2$ in Eq.~\ref{sec_2a_2_eq_6} below.

\begin{eqnarray*}\label{sec_2a_2_eq_6}  
[\frac{d^2 E_0(y)}{dy^2}]_2 < (\pi)^{-\frac{1}{4}}  \displaystyle\sum\limits_{n=11,13,...}^{\infty}  e^{- \frac{n^{2}}{4} y } y^{\frac{1}{4}}  \frac{n^{4}}{16}  
+ e^{- n^{2} y } y^{\frac{1}{4}}  ( \frac{3}{16 y^{2}} + \frac{n^{2}}{2y}  )  \\
+ e^{- \frac{(n+1)^{2}}{4} y }  y^{\frac{1}{4}}  ( \frac{3}{16 y^{2}} + \frac{(n+1)^{2}}{8y}  )  + e^{- (n+1)^{2} y}  y^{\frac{1}{4}}   (n+1)^{4}  \\
\end{eqnarray*}
\begin{equation} \end{equation}

We compute the maximum value of $[\frac{d^2 E_0(y)}{dy^2}]_2$ in Eq.~\ref{sec_2a_2_eq_6} for $ n \geq 11$ by setting first term as $e^{- \frac{n^{2}}{4} y_{min} } (y_{max})^{\frac{1}{4}}$ and using $n+1 < 1.1 n$, $n^{2} <  10 e^{0.1 n^{2}} =10 [ 1 + 0.1 n^{2} + \frac{0.01}{2} n^{4}+....]$ and \\
$n^{4} < 200 e^{0.1 n^{2}}=200 [ 1 + 0.1 n^{2} + \frac{0.01}{2} n^{4}+....]$ and take the common term $y_{max}^{\frac{1}{4}}$ outside the summation, as follows.

\begin{eqnarray*}\label{sec_2a_2_eq_7}  
[\frac{d^2 E_0(y)}{dy^2}]_2  < (\pi)^{-\frac{1}{4}}  y_{max}^{\frac{1}{4}} \displaystyle\sum\limits_{n=11,13,...}^{\infty}  e^{- \frac{n^{2}}{4} y_{min} } 200 e^{0.1 n^{2}}\frac{1}{16}  
+ e^{- n^{2}  y_{min} }  ( \frac{3}{16  y_{min}^{2}} +  10 e^{0.1 n^{2}} \frac{1}{2  y_{min}}  )  \\
+ e^{- \frac{(n+1)^{2}}{4}  y_{min} }  ( \frac{3}{16  y_{min}^{2}} +  10 e^{0.1 n^{2}} \frac{(1.1)^{2}}{8  y_{min}}  ) 
  + e^{- (n+1)^{2}  y_{min}} 200 e^{0.1 n^{2}}   (1.1)^{4}  
\end{eqnarray*}
\begin{equation} \end{equation}

We use $n+1 > n$ and use $e^{- (n+1)^{2}  y_{min} K } < e^{- (n)^{2}  y_{min} K}$ for $K>0$ in the exponent terms and  simplify above equation as follows.

\begin{eqnarray*}\label{sec_2a_2_eq_8}  
[\frac{d^2 E_0(y)}{dy^2}]_2 < (\pi)^{-\frac{1}{4}}  y_{max}^{\frac{1}{4}} \displaystyle\sum\limits_{n=11,13,...}^{\infty}  e^{- n^{2} (\frac{1}{4} y_{min} -0.1) }  \frac{200}{16}  
+ e^{- n^{2} y_{min} } \frac{3}{16 y_{min}^{2}}  +  e^{- n^{2} (y_{min}-0.1) } \frac{10}{2y_{min}}    \\
+ e^{- n^{2} \frac{1}{4} y_{min} }   \frac{3}{16 y_{min}^{2}} + e^{- n^{2} ( \frac{1}{4} y_{min} - 0.1) }   \frac{(1.1)^{2}*10}{8y_{min}}  
  + e^{- n^{2} ( y_{min} - 0.1) } 200   (1.1)^{4}
\end{eqnarray*}
\begin{equation} \end{equation}

We use the complementary error function given by $erfc(z)= \frac{2}{\sqrt{\pi}} \int_{z}^{\infty} e^{-u^{2}} du $ \href{https://mathworld.wolfram.com/Erfc.html}{link}) and the fact that $\displaystyle\sum\limits_{n=11,13,...}^{\infty}  e^{- n^{2} K} < \int_{11}^{\infty} e^{-t^{2} K} dt = \frac{1}{ \sqrt{K}} \int_{11  \sqrt{K}}^{\infty} e^{-u^{2}} du = \frac{\sqrt{\pi}}{ 2 \sqrt{K}}  erfc(11  \sqrt{K})$ using the substitution $t \sqrt{K} = u$ and $dt \sqrt{K} = du$ for $K>0$ and write Eq.~\ref{sec_2a_2_eq_8} as follows.  We note that $y_{min}=\pi$ and hence $\frac{1}{4} y_{min} -0.1 >0$ and $y_{min} -0.1 >0$.

\begin{eqnarray*}\label{sec_2a_2_eq_9}  
[\frac{d^2 E_0(y)}{dy^2}]_2 < (\pi)^{-\frac{1}{4}}  y_{max}^{\frac{1}{4}} [  \frac{200}{16}  \frac{\sqrt{\pi}}{ 2 \sqrt{(\frac{1}{4} y_{min} -0.1)}}  erfc(11  \sqrt{(\frac{1}{4} y_{min} -0.1)})  + \frac{3}{16 y_{min}^{2}}   \frac{\sqrt{\pi}}{ 2 \sqrt{ y_{min} }}  erfc(11  \sqrt{ y_{min} }) \\
+  \frac{10}{2y_{min}} \frac{\sqrt{\pi}}{ 2 \sqrt{(y_{min}-0.1)}}  erfc(11  \sqrt{(y_{min}-0.1)}) 
+ \frac{3}{16 y_{min}^{2}}   \frac{\sqrt{\pi}}{ 2 \sqrt{\frac{1}{4} y_{min}}}  erfc(11  \sqrt{\frac{1}{4} y_{min}} \\
+  \frac{(1.1)^{2}*10}{8y_{min}}  \frac{\sqrt{\pi}}{ 2 \sqrt{( \frac{1}{4} y_{min} - 0.1)}}  erfc(11  \sqrt{( \frac{1}{4} y_{min} - 0.1)}) \\+ 200   (1.1)^{4} \frac{\sqrt{\pi}}{ 2 \sqrt{( y_{min} - 0.1)}}  erfc(11  \sqrt{( y_{min} - 0.1)})]  
\end{eqnarray*}
\begin{equation} \end{equation}

We compute Eq.~\ref{sec_2a_2_eq_9} numerically and get $[\frac{d^2 E_0(y)}{dy^2}]_2 < 8.65 * 10^{-37}$.  The \textbf{maximum} value of $[\frac{d^2 E_0(y)}{dy^2}]_2$  in Eq.~\ref{sec_2a_2_eq_8} at $n=11,13,..$ in the interval $y=[y_{min},y_{max})$ is given by $8.65 * 10^{-37}$ which is \textbf{positive}. (\textbf{Result E.5.6}) \href{https://www.ocf.berkeley.edu/~araman/files/math_z/test_Phi_inv_all_n.m}{Matlab simulation})\\

Using Result E.5.5 and E.5.6, we get the \textbf{maximum} value of $\frac{d^2 E_0(y)}{dy^2}$  in Eq.~\ref{sec_2a_2_eq_5} at $n=1,3,5,..$ in the interval $y=[y_{min},y_{max})$ is given by $-0.0097 + 8.65 * 10^{-37} \approx -0.0097$ which is \textbf{negative}. (\textbf{Result E.5.7})\\

Hence we have shown that  $\frac{d^2 E_0(y)}{dy^2} <0$ , for $\pi \leq y < y_a $ and hence $\frac{d E_0(y)}{dy} <0$ for $\pi < y < y_a $ given that $\frac{d E_0(y)}{dy} = 0$ at $y=\pi$, using Result E.5.a. \\

It is shown in Section~\ref{sec:Section_2a_2a} that $\frac{d E_0(y)}{dy} <0$ for $y_a \leq y < \infty$ for all finite $n$.\\

Hence  $\frac{d E_0(y)}{dy} <0$ for $\pi < y < \infty$. We see that $y= \pi   e^{2t} $ and $ \frac{dy}{dt} = 2 \pi e^{2t}= 2y $  and $\frac{d E_0(-t)}{dt} = \frac{d E_0(-t)}{dy} \frac{dy}{dt} =  \frac{d E_0(y)}{dy} 2 y$ and hence  $\frac{d E_0(-t)}{dt}  < 0$ for $t > 0$. Hence $E_0(t)=E_0(-t)$ is strictly decreasing for $t > 0$.


.\\  \tocless\subsection{\label{sec:Section_2a_3} \textbf{Proof of Riemann's Hypothesis for Dirichlet Eta function}  \protect\\  \lowercase{} }

The proof of Riemann's Hypothesis presented in Section~\ref{sec:Section_2} to Section~\ref{sec:Section_3} can be used to prove Riemann's Hypothesis for Dirichlet Eta function, as detailed below.\\

$\bullet$ In Section~\ref{sec:Section_2} to Section~\ref{sec:Section_3} and ~\ref{sec:appendix_A} to ~\ref{sec:appendix_I}, we replace $\xi(s)$ with $E(s)$ which is a holomorphic function in the region $0 < Re[s] < 1$ and replace $E_0(t)$ for Riemann's Xi function with $E_0(t)= \displaystyle\sum\limits_{n=1}^{\infty} (-1)^{n-1}  ( e^{-\pi \frac{n^{2}}{4} e^{ -2 t} } - e^{-\pi n^{2} e^{ -2 t} } ) e^{ -\frac{t}{2} } $ for Dirichlet Eta function detailed in ~\ref{sec:Sec_A_1_7a}, which is a real and \textbf{even function} of variable $t$.\\

We use $\xi(s) = \xi(1-s) =  \frac{1}{2} s (s-1) \pi^{-\frac{s}{2}} \Gamma(\frac{s}{2}) \zeta(s) $ where $\zeta(s)=\displaystyle\frac{\eta(s)}{1-2^{1-s}}$. \href{https://www.ocf.berkeley.edu/~araman/files/math_z/Titchmarsh_pp16_17.pdf}{Titchmarsh pp16-17}) Using Result E.0 and Eq.~\ref{sec:sec_a_1_7_eq_5}, we see that, if Dirichlet Eta function has a zero in the critical strip at $\omega=\omega_0$, then the Fourier transform of $E_p(t)= E_0(t)  e^{-\sigma t}$ \textbf{also} has a zero at $\omega=\omega_0$ for $0 < |\sigma| < \frac{1}{2}$.\\

$\bullet$ We replace Section~\ref{sec:Section_A_1_6} with ~\ref{sec:Section_2a_2}  which shows that  $E_0(t)$ is a \textbf{strictly decreasing} function for $t>0$.\\

$\bullet$ We use the result in ~\ref{sec:appendix_C_7} with $\xi(s)$ replaced by $E(s)$ and use $E(s)=E(1-s)$ using Eq.~\ref{sec:sec_a_1_7_eq_1} and show that $E_0(t)$ which is a real and \textbf{even} function of $t$.\\

$\bullet$ We replace ~\ref{sec:appendix_C_5} with ~\ref{sec:appendix_C_5b} which shows that $E_0(t), E_p(t), x(t)=E_0(t) e^{-2 \sigma t}$ have exponential fall off rate.\\

$\bullet$ We replace ~\ref{sec:appendix_C_5a} with ~\ref{sec:Section_2a_1} which shows that $E_0(t) > 0$ for $-\infty < t < \infty$.\\

$\bullet$ We show the above result for $0 < \sigma < \frac{1}{2}$. Given that $E_p(t) =  E_0(t) e^{-\sigma t}$  is real, its Fourier transform $E_{p\omega}(\omega)=E_{pR\omega}(\omega) + i E_{pI\omega}(\omega) $ has symmetry properties and hence $E_{pR\omega}(-\omega)=E_{pR\omega}(\omega)$  and $E_{pI\omega}(-\omega)= - E_{pI\omega}(\omega)$ (\href{https://lpsa.swarthmore.edu/Fourier/Xforms/FXProps.html}{Symmetry property of Fourier Transform}) also have a zero at $\omega=\omega_0$ and hence $E_{p\omega}(-\omega)=E_{pR\omega}(\omega) - i E_{pI\omega}(\omega) $ \textbf{also} has a zero at $\omega=\omega_0$ to satisfy Statement A. \\

If $E_{p\omega}(\omega)$ and $\eta(\frac{1}{2}+\sigma + i \omega)$ has a zero at $\omega=\omega_0$ to satisfy Statement A, then $E_{p\omega}(-\omega)$ and $\eta(\frac{1}{2}+\sigma - i \omega)$ also has a zero at $\omega=\omega_0$(using last paragraph)  and $\eta(\frac{1}{2}-\sigma + i \omega)$ also has a zero at $\omega=\omega_0$ using the functional equation for Dirichlet Eta function derived in ~\ref{sec:Section_2a_2z} which relates $\eta(s)$ and $\eta(1-s)$. Hence the results in above sections hold for $-\frac{1}{2} < \sigma < 0$ and for  $0 < |\sigma| < \frac{1}{2}$.\\

$\bullet$ Complete details in \href{https://zenodo.org/record/8383224}{(link)}. 

.\\  \tocless\subsection{\label{sec:appendix_C_5c} \textbf{Exponential fall-off rate of Dirichlet Eta function}  \protect\\  \lowercase{} }

The integrand in Eq.~\ref{sec_2a_eq_3} given by $\displaystyle\sum\limits_{n=1}^{\infty} (-1)^{n-1}   e^{-\pi n^{2} e^{ -2 t} }  e^{ -\frac{t}{2} } e^{-\sigma t}$ goes to zero with \textbf{exponential} fall-off rate, as $t \to -\infty$ because the term $ e^{-\pi n^{2} e^{ -2 t} } $ has a faster fall-off rate than the term $ e^{ -\frac{t}{2} } e^{-\sigma t}$, for each $n$.\\

The integrand in Eq.~\ref{sec_2a_eq_3} given by $\displaystyle\sum\limits_{n=1}^{\infty} (-1)^{n-1}   e^{-\pi n^{2} e^{ -2 t} }  e^{ -\frac{t}{2} } e^{-\sigma t}$ goes to zero with \textbf{exponential} fall-off rate, as $t \to +\infty$ because the term $ \lim_{t \to \infty}  e^{-\pi n^{2} e^{ -2 t} }  = 1$ for each $n$ and hence  \\$ \lim_{t \to \infty} \displaystyle\sum\limits_{n=1}^{\infty} (-1)^{n-1}   e^{-\pi n^{2} e^{ -2 t} } = 1-1+1-1... = \frac{1}{2} $ \href{https://www.ocf.berkeley.edu/~araman/files/math_z/Hardy_Divergent_Series.pdf#page=19}{(Eq.1.2.7 in page 2)}  and the term  $  \lim_{t \to \infty}  e^{ -\frac{t}{2} } e^{-\sigma t} =0$ for $0 < \sigma < \frac{1}{2}$.\\


.\\  \tocless\subsection{\label{sec:Section_2a_2z} \textbf{Functional equation for Dirichlet Eta function}  \protect\\  \lowercase{} }

We use the \textbf{functional equation} for Riemann's zeta function given by $\zeta(s)= \zeta(1-s) \Gamma(1-s) \sin{(\frac{s \pi}{2})} \pi^{(s-1)} 2^s$ and use  $\zeta(s)=\displaystyle\frac{\eta(s)}{1-2^{1-s}}$ and $s=\frac{1}{2}+\sigma+i \omega$ and  $1-s=\frac{1}{2}-\sigma-i \omega$. 

\begin{eqnarray*}\label{sec_2b_2_eq_1}  
\zeta(s)= \zeta(1-s) \Gamma(1-s) \sin{(\frac{s \pi}{2})} \pi^{(s-1)} 2^s \\
\displaystyle\frac{\eta(s)}{1-2^{1-s}}= \displaystyle\frac{\eta(1-s)}{1-2^{s}} \Gamma(1-s) \sin{(\frac{s \pi}{2})} \pi^{(s-1)} 2^s
\end{eqnarray*}
\begin{equation} \end{equation}

We use well known properties of Gamma function $\Gamma(s) \Gamma(1-s) = \frac{\pi}{\sin{(s \pi)}}=  \frac{\pi}{2 \sin{( \frac{s \pi}{2})} \cos{( \frac{s \pi}{2})}}$ in Eq.~\ref{sec_2b_2_eq_1} as follows. \href{https://www.ocf.berkeley.edu/~araman/files/math_z/Ellison_p147-152.pdf#page=4}{(link)}

\begin{equation} \label{sec_2b_2_eq_2}  
\displaystyle\frac{\eta(s)}{1-2^{1-s}}= \displaystyle\frac{\eta(1-s)}{1-2^{s}} \frac{\pi}{2 \sin{( \frac{s \pi}{2})} \cos{( \frac{s \pi}{2})}\Gamma(s)} \sin{(\frac{s \pi}{2})} \pi^{(s-1)} 2^s
\end{equation}

We cancel the common term $ \sin{(\frac{s \pi}{2})}$ in Eq.~\ref{sec_2b_2_eq_2} for $0 < Re[s] < 1$ and rearrange the terms as follows.

\begin{equation}  \label{sec_2b_2_eq_3}  
\eta(1-s) = \eta(s)   \Gamma(s) \cos{( \frac{s \pi}{2})} \frac{(1-2^{s}) }{(1-2^{1-s}) \pi^{s} 2^{s-1}} 
\end{equation}

In the modified functional equation in Eq.~\ref{sec_2b_2_eq_3} , we see that,  \textbf{if} Dirichlet Eta function $\eta(s)$ has a zero in the region $0 < Re[s] < 1$ at $s=s_0$, \textbf{then} $\eta(s)$ also has a zero at  $s=1-s_0$, due to the term $\eta(1-s)$, given that for $Re[s]>0$, the gamma function is analytic in the complex plane \href{https://proofwiki.org/wiki/Poles_of_Gamma_Function}{(link)}.

.\\  \tocless\subsection{\label{sec:appendix_C_5b} \textbf{ Exponential Fall off rate of  $E_0(t)$ and $E_p(t)=E_0(t) e^{-\sigma t}$ and $x(t)=E_0(t) e^{-2 \sigma t}$ } \protect\\  \lowercase{} }

Given that $E_0(t)=E_0(-t)$ (using Result E.1 in ~\ref{sec:Sec_A_1_7a}), we write $E_0(t)$ in Eq.~\ref{sec:sec_a_1_7_eq_5} as follows.

\begin{eqnarray*}\label{sec_app_C_5b_eq_1}
E_0(t) =   \displaystyle\sum\limits_{n=1}^{\infty} (-1)^{n-1}  ( e^{-\pi \frac{n^{2}}{4} e^{ 2 t} } - e^{-\pi n^{2} e^{2 t} } ) e^{ \frac{t}{2} } =  \displaystyle\sum\limits_{n=1}^{\infty} (-1)^{n-1}  e^{-\pi \frac{n^{2}}{4} e^{ 2 t} } (1 - e^{-\pi  \frac{3 n^{2}}{4} e^{2 t} }  ) e^{ \frac{t}{2} } 
\end{eqnarray*}
\begin{equation} \end{equation}

We use Taylor series expansion around $t=0$ for $e^{2t}= \displaystyle\sum_{r=0}^{\infty} \frac{(2t)^{r}}{!r}$, given that $e^{2t}$ is an analytic function for real $t$.

\begin{eqnarray*}\label{sec_app_C_5b_eq_1}
E_0(t) = \sum_{n=1}^{\infty}   (-1)^{n-1}  e^{- \pi \frac{n^{2}}{4} (1+ 2t)} e^{- \pi\frac{n^{2}}{4} ( \frac{(2t)^2}{!2} +  \frac{(2t)^3}{!3} ....)} ( 1 - e^{-\pi  \frac{3 n^{2}}{4} e^{2 t} } ) e^{\frac{t}{2}} 
\end{eqnarray*}
\begin{equation} \end{equation}

We take the term $ e^{- \frac{\pi}{2} t} e^{\frac{t}{2}} = e^{-1.0708 t}$ out of the summation, corresponding to $n=1$ and write Eq.~\ref{sec_app_C_5b_eq_1} as follows.

\begin{equation} \label{sec_app_C_5b_eq_2}
E_0(t) = e^{- \frac{\pi}{2} t} e^{\frac{t}{2}} \sum_{n=1}^{\infty}   (-1)^{n-1}  e^{- \pi \frac{n^{2}}{4}}  e^{- \frac{\pi}{2} (n^{2}-1) t} e^{- \pi\frac{n^{2}}{4} ( \frac{(2t)^2}{!2} +  \frac{(2t)^3}{!3} ....)} ( 1 - e^{-\pi  \frac{3 n^{2}}{4} e^{2 t} } ) 
 \end{equation}

For $t > 0$, we see that the term corresponding to $n=1$ in Eq.~\ref{sec_app_C_5b_eq_2} has an asymptotic fall-off rate of $o[e^{- t} ]$. The terms corresponding to $n > 1$ have fall-off rates \textbf{higher} than $o[e^{- t} ]$, due to the term $e^{- \frac{\pi}{2} (n^{2}-1) t}$.\\

Hence we see that $E_0(t)$ has an asymptotic fall-off rate of $o[e^{-t} ]$, for $t>0$. Given that $E_0(t)=E_0(-t)$(~\ref{sec:appendix_C_7}), we see that $E_0(t)$ has an \textbf{exponential} asymptotic fall-off rate of $o[e^{- |t|} ]$.\\

Similarly, $E_p(t)= E_0(t) e^{- \sigma t}$ has an asymptotic \textbf{exponential} fall-off rate of  $ o[e^{-0.5 |t|} ] $ (using $o[e^{-(1-  \sigma) |t| } ]$ ), for $0 \leq |\sigma| < \frac{1}{2}$.\\

Similarly, $x(t)=E_0(t) e^{-2 \sigma t}$ has an asymptotic \textbf{exponential} fall-off rate of   $ o[e^{-\delta |t|} ] $ (using $o[e^{-(1- 2 \sigma) |t| } ]$), for $0 \leq |\sigma| < \frac{1}{2}$ and $0 < \delta << 1$.\\

\end{document}